\documentclass[a4paper,10pt]{article}
\usepackage[utf8x]{inputenc}
\usepackage{graphicx}
\usepackage{amsmath,bbm} 
\usepackage{amsfonts} 
\usepackage{amssymb}
\usepackage[ruled,linesnumbered]{algorithm2e}

\graphicspath{{Figures/}}

\title{2D Empirical Transforms.\\ Wavelets, Ridgelets and Curvelets revisited}
\author{J\'er\^ome Gilles and Giang Tran and Stanley Osher \thanks{University of California, Los Angeles (UCLA), Department of Mathematics, 520 Portola Plaza, Los Angeles, CA 90095 ({\tt jegilles@math.ucla.edu, 
giangtran@math.ucla.edu, sjo@math.ucla.edu})}}

\def\F{\mathcal{F}}
\def\W{\mathcal{W}}
\def\E{\mathcal{E}}
\def\B{\mathcal{B}}
\def\x{\mathbf{x}}
\def\o{\mathbf{\omega}}

\newtheorem{theorem}{Theorem}
\newtheorem{proof}{Proof}

\begin{document}

\maketitle

\newcommand{\slugmaster}{
\slugger{siims}{xxxx}{xx}{x}{x--x}}

\begin{abstract}
A recently developed new approach, called ``Empirical Wavelet Transform'', aims to build 1D adaptive wavelet frames accordingly to the analyzed signal. In this paper, we present 
several extensions of this approach to 2D signals (images). We revisit some well-known transforms (tensor wavelets, Littlewood-Paley wavelets, ridgelets and curvelets) and show that it is possible 
to build their empirical counterpart. We prove that such constructions lead to different adaptive frames which show some promising properties for image analysis and processing.
\end{abstract}

\section{Introduction}
Wavelets and their geometric extensions (framelets, ridgelets, curve\-lets,\ldots) are not only a useful mathematical tool arising in harmonic analysis and used to study function spaces, but they are also very efficient in image processing. 
For instance, in texture analysis, such transforms permit us to extract relevant characteristics which can be used as inputs of a classifier. Since the advent of compressive sensing theory, many variational models have been proposed in the 
literature, using the sparsity properties in these transform domains \cite{Cai2009,Cai2009a,Cai2009b,Shen2010,Starck2005}, to perform efficient denoising or deconvolution. While wavelets are well established in the scientific 
community, their foundation has not really evolved since their origin and almost all discrete wavelet bases are built on a prescribed scheme corresponding to the multiresolution analysis definition. The main consequence of using 
such a definition is that those basis are based on dyadic scale decomposition which does not ensure that the corresponding filters are the optimal ones to represent an image. For instance, it is easy to build a texture 
based on compact information lying between two scales and which will be separated by a usual wavelet transform. A better approach is to have an adaptive representation where the basis is generated accordingly to the 
information contained in the image itself. Surprisingly, there are only a few attempts to build such adaptive representations. One of the first ones was the Brushlet transform \cite{Borup2003,Meyer1997} which aims to build some 
wavelet filters based on adaptive supports in the Fourier domain. While this construction is difficult to implement and did not get a real success, the corresponding idea is interesting and, as we will present later, 
will inspire our work. Mallat and his students investigated the opportunity to extract some geometric information from the image itself to build adaptive wavelets representations like the bandlet transform \cite{LePennec2005a, LePennec2005b,
peyre2004} or the geometrical grouplet transform \cite{Mallat2009}. A completely different approach to design an adaptive representation exists in the literature: the empirical mode decomposition (EMD). In the 1D case, the EMD was proposed 
by Huang et al. \cite{Huang1998} where the aim of the method is to extract amplitude modulated-frequency modulated (AM-FM) components of an input signal. In practice, it is based on the detection of the upper and lower envelopes 
which provide the global trend of the signal. This low varying signal is then subtracted to the original one to get the highest oscillatory component. This procedure is repeated until all oscillating parts are extracted.
Bidimensional extensions \cite{Damerval2005,Numes2003} of this algorithm were proposed to process images. The main drawback of the EMD approach is that it is a pure nonlinear algorithmic procedure and the obtained 
representation is implementation dependent (e.g it depends on how the envelopes are detected, which interpolation process is used and the chosen stopping criteria). Moreover due to the nonlinear aspects, no theoretical
background supports this method. Recently, in \cite{Gilles2013} the author proposed the construction of 1D empirical wavelets where the aim is, like the EMD, to extract AM-FM components from a signal. The definition of the 
empirical wavelets is based on the formulation of a Littlewood-Paley \cite{grafakos2008classical} wavelet where the choice of their supports in the Fourier domain is not prescribed to a dyadic tiling but chosen accordingly to the analyzed signal. The author showed 
that it is possible to build a wavelet tight frame which corresponds to an adaptive filter bank. The empirical wavelet transform (EWT) consists of two major steps: detect the Fourier supports and build the corresponding 
wavelet accordingly to those supports; filter the input signal with the obtained filter bank to get the different components.\\
In \cite{Gilles2013}, a straightforward tensor extension was proposed to do some image analysis. In this paper, we propose to investigate the possibility of generalizing this empirical approach to different existing 2D 
wavelet transforms. The remainder of the paper is as follows. In section \ref{sec:background}, we define some notations which will be used throughout the paper and we recall the definition of the 1D EWT. In section 
\ref{sec:bound}, we investigate an important aspect of the EWT: the detection of the Fourier supports. Section \ref{sec:tensor} recalls the 2D tensor extension while sections \ref{sec:lpewt}, \ref{sec:eridge} and 
\ref{sec:ecurv} present the empirical extensions of the 2D Littlewood-Paley wavelet transform, the ridgelet transform and the curvelet transform, respectively. Some experiments will be given in section \ref{sec:expe} and
we conclude in section \ref{sec:conc}.

\section{Background}\label{sec:background}
\subsection{Notation}
Let us start to fix some notations which will be used throughout the paper. We denote $\x=(x_1,x_2)$ a spatial position in the 2D plane, $\o=(\omega_1,\omega_2)$ 
the coordinates in the 2D frequency plane (in the case of a 1D signal we will use $\omega$ to represent the 1D frequency). 
The convolution product will be denoted by $\star$. We define the following operators:
\begin{itemize}
 \item $\F_{1,y}$, $\F_{1,y}^*$ the usual 1D Fourier transform and its inverse with respect to the $y$ variable,
 \item $\F_2$, $\F_2^*$ the usual 2D Fourier transform and its inverse,
 \item $\F_P$, $\F_P^*$ the 2D Pseudo-Polar Fourier transform and its adjoint \cite{PseudoPolar} (their definitions are recalled in appendix~\ref{app:ppfft}),
 \item $\W_{1,y}$, $\W_{1,y}^*$ the standard dyadic 1D wavelet transform and its inverse with respect to the $y$ variable,
\end{itemize}

\subsection{The 1D EWT}\label{subsec:1dewt}
In \cite{Gilles2013}, the author proposed to build an empirical wavelet transform (EWT). The idea consists of defining a bank of $N$ wavelet filters 
(one lowpass and $N-1$ bandpass filters corresponding to the approximation and details components, respectively) based on 
``well chosen'' Fourier supports (meaning we select relevant modes in the signal spectrum). If we denote $f(t)$ a 1D signal, then we first detect the 
boundaries of each Fourier support from $|\F_{1,t}(f)|(\omega)$ (we will give more details about the procedure to detect such supports in section \ref{sec:bound}).
This operation provides us with a set of boundaries $\Omega=\{\omega^n\}_{n=0,\ldots,N}$ (we restrict our discussion to the interval $[0,\pi]$ and take the 
convention $\omega^0=0$ and $\omega^N=\pi$). Based on $\Omega$ (see Figure.~\ref{fig:perio}), we can define a wavelet tight frame, $\B=\left\{\phi_1(t),\{\psi_n(t)\}_{n=1}^{N-1}\right\}$, 
inspired from Meyer's and Littlewood-Paley wavelets, their Fourier transforms are given by:

\begin{figure}
\centering
\includegraphics{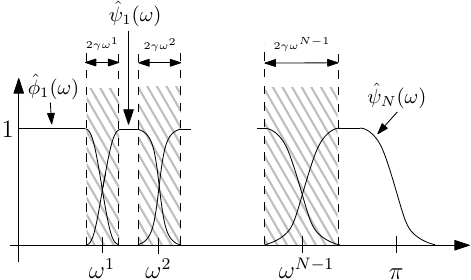}
\caption{Fourier line decomposition principle and EWT basis construction.}
\label{fig:perio}
\end{figure}

\begin{equation}\label{eq:phi2}
\F_{1,t}(\phi_1)(\omega)=
\begin{cases}
1  \; &\text{if}\;|\omega|\leq (1-\gamma)\omega^1\\
\cos\left[\frac{\pi}{2}\beta\left(\frac{1}{2\gamma\omega^1}(|\omega|-(1-\gamma)\omega^1)\right)\right] \; &\text{if}\; (1-\gamma)\omega^1\leq |\omega|\leq (1+\gamma)\omega^1 \\
0  \; &\text{otherwise,}
\end{cases}
\end{equation}
and
\begin{equation}\label{eq:psi2}
\F_{1,t}(\psi_n)(\omega)=
\begin{cases}
1 \qquad\qquad\qquad\qquad\qquad\quad\;\;\; \text{if}\; (1+\gamma)\omega^n\leq |\omega|\leq (1-\gamma)\omega^{n+1} \\
\cos\left[\frac{\pi}{2}\beta\left(\frac{1}{2\gamma\omega^{n+1}}(|\omega|-(1-\gamma) \omega^{n+1})\right)\right]  \\
\qquad\qquad\qquad\qquad\qquad\qquad\; \text{if}\; (1-\gamma)\omega^{n+1}\leq |\omega|\leq (1+\gamma)\omega^{n+1}\\
\sin\left[\frac{\pi}{2}\beta\left(\frac{1}{2\gamma\omega^n}(|\omega|-(1-\gamma)\omega^n)\right)\right] \\
\qquad\qquad\qquad\qquad\qquad\qquad\;  \text{if}\; (1-\gamma)\omega^n\leq |\omega|\leq (1+\gamma)\omega^n\\
0 \qquad\qquad\qquad\qquad\qquad\quad\;\;\; \text{otherwise.}
\end{cases}
\end{equation}
Where $\beta$ is an arbitrary $\mathcal{C}^k([0,1])$ function, fulfilling the properties $\beta(x)=0$ if $x\leq 0$, $\beta(x)=1$ if $x\geq 1$ 
and $\beta(x)+\beta(1-x)=1, \forall x\in[0,1]$. For instance, in the construction of Meyer's wavelet, Daubechies \cite{Daubechies1992}
proposed to use
\begin{equation}\label{eq2}
\beta(x)=
\begin{cases}
0 & \text{if}\; x\leq 0\\
1 & \text{if}\; x\geq 1\\
x^4(35-84x+70x^2-20x^3) & \text{if}\; \forall x\in [0,1].
\end{cases}
\end{equation}
The parameter $\gamma$ allows us to ensure that two consecutive transitions areas (dashed regions in Figure.~\ref{fig:perio}) do not overlap. A necessary condition on $\gamma$ is given in \cite{Gilles2013} in order to 
have the tight frame property of $\B$ and allows us to automatically 
choose this parameter. Equipped with this set of filters, the 1D empirical wavelet transform (1D-EWT) is defined by
\begin{equation}
\mathcal{W}^{\E}(f)(n,t)=\F_{1,\omega}^*\left(\F_{1,t}(f)(\omega)\overline{\F_{1,t}(\psi_n)(\omega)}\right),
\end{equation}
for the detail coefficients and the approximation coefficients (we adopt the convention $\mathcal{W}_f^{\E}(0,t)$ to denote them) by:
\begin{equation}
\mathcal{W}^{\E}(f)(0,t)=\F_{1,\omega}^*\left(\F_{1,t}(f)(\omega)\overline{\F_{1,t}(\phi_1)(\omega)}\right).
\end{equation}
The inverse transform is straightforward by using inverse/adjoint operators:
\begin{align}
f(t)&=\mathcal{W}^{\mathcal{E}}(f)(0,t)\star\phi_1(t)+\sum_{n=1}^{N-1}\mathcal{W}^{\mathcal{E}}(f)(n,t)\star\psi_n(t)\\
&=\F_{1,\omega}^*\left(\F_{1,t}(\mathcal{W}^{\E})(0,\omega)\F_{1,t}(\phi_1)(\omega)+\sum_{n=1}^{N-1}\F_{1,t}(\mathcal{W}^{\E})(n,\omega)\F_{1,t}(\psi_n)(\omega)\right).
\end{align}
This EWT shows interesting results in extracting different information contained in the signal and permits us to obtain, with the help 
of the Hilbert transform, a precise time-frequency representation for 1D signals \cite{Gilles2013}. It is clear that the process detecting the Fourier boundaries is an important step. The next section recalls 
the approach used in \cite{Gilles2013} and proposes new ones.

\section{Fourier boundaries detection}\label{sec:bound}
The original method proposed in \cite{Gilles2013} to detect the set of Fourier boundaries, $\Omega$, is based on the detection of local maxima in the spectrum magnitude, assuming that the spectrum is composed 
of sufficiently separated modes. However, in the case of images, except for highly textured 
images, their spectra do not exhibit obvious distinct modes and other criteria must be used to detect useful Fourier supports. After recalling the method based on the local maxima, we will introduce some improvements as 
well as other options to perform this detection task on the Fourier spectrum. In the following, we use the notation $H(\omega)=|\F_{1,y}|(\omega)$, restricted to $\omega\in[0,\pi]$ and we keep the convention given in 
section \ref{subsec:1dewt}: $\omega^0=0$ and $\omega^N=\pi$.

\subsection{Original detection of local maxima method}
\begin{figure}[!t]
\begin{center}
\begin{tabular}{c}
\begin{tabular}{cc}
\includegraphics[width=0.34\textwidth]{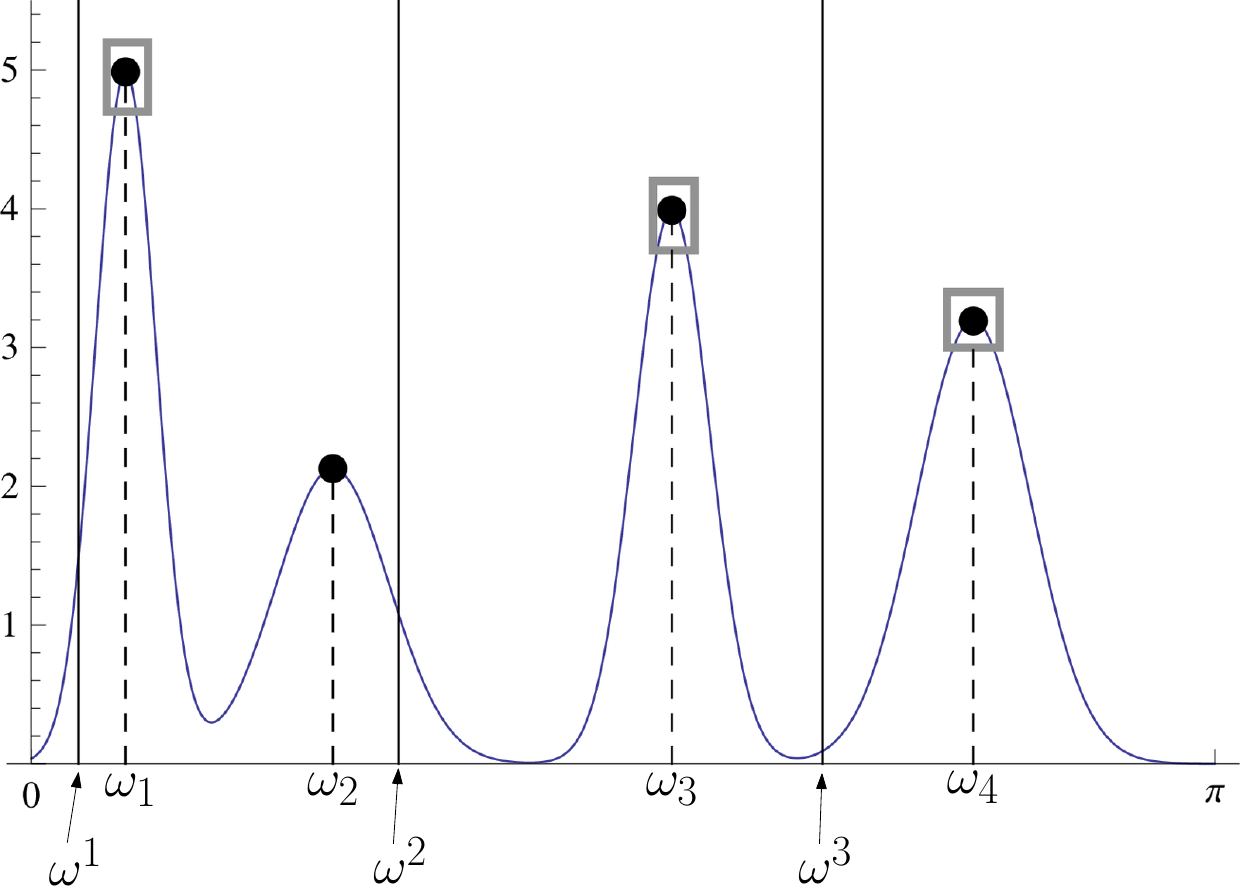} & \includegraphics[width=0.34\textwidth]{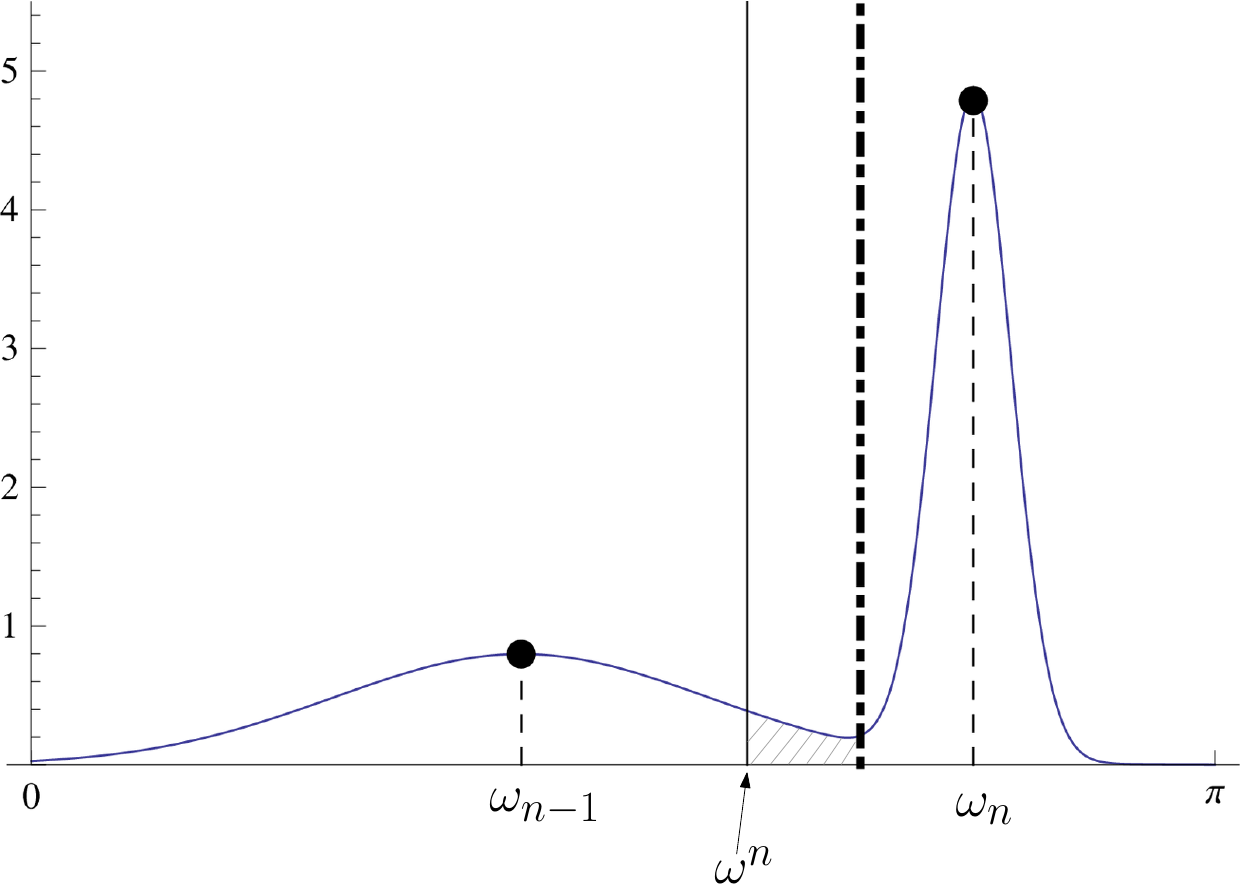}\\
(a) Local maxima detection principle. & (b) The Flat-Picked modes issue.\vspace{5mm}
\end{tabular}\\
\includegraphics[width=0.34\textwidth]{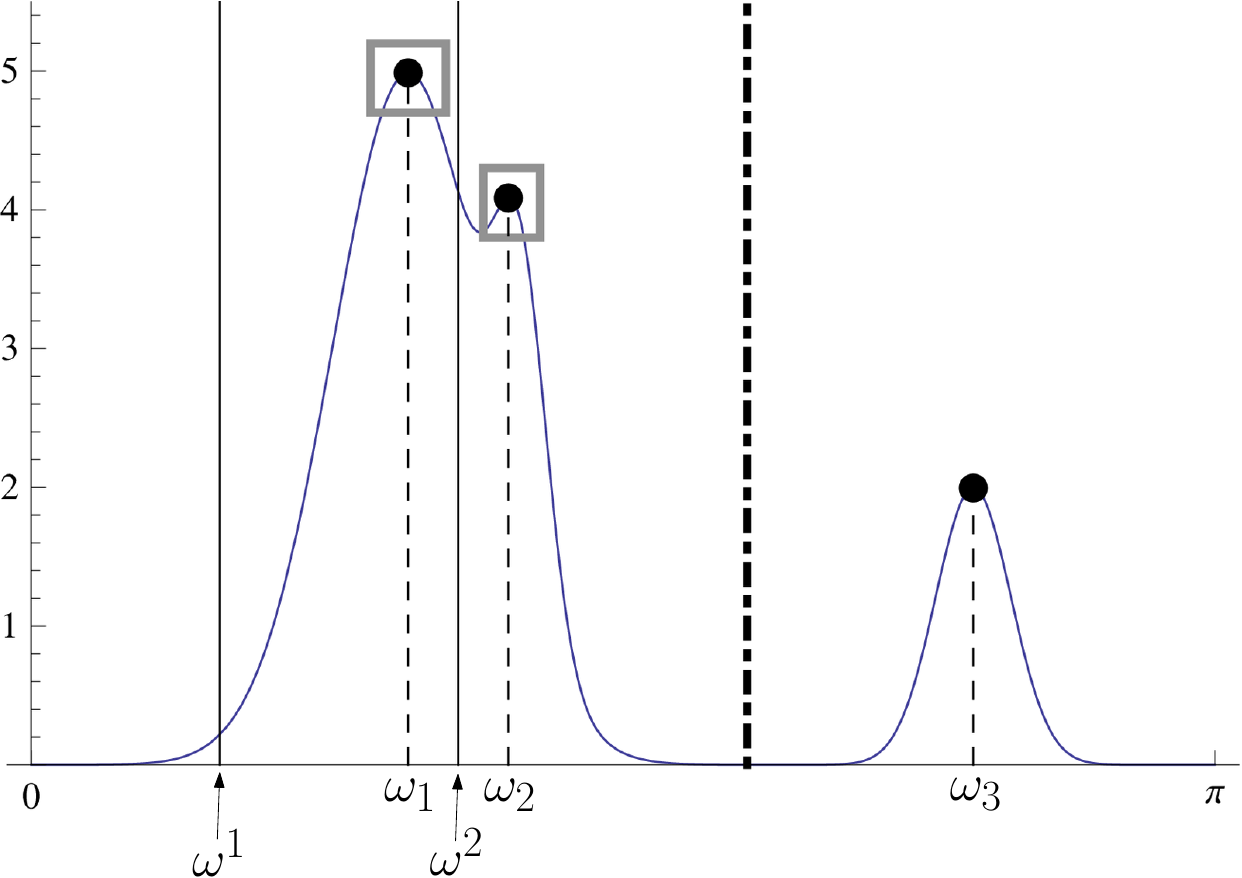}\\
(c) Global v.s local modes
\end{tabular}
\end{center}
\caption{Fourier supports definition based on local maxima detection and two of its main issues. See the text for full explanation.}
\label{fig:locmax}
\end{figure}

In this approach, see Figure.~\ref{fig:locmax}.a, we first compute all local maxima $M_i$ of $H$ (the black dots on the figure) and deduce their corresponding position $\omega_i$. Next we keep the set of all $\omega_i$ 
corresponding to the $N-1$ largest maxima (the square boxes on Figure.~\ref{fig:locmax} illustrate the selected maxima when $N=4$) and re-index them as $\omega_n$ where $1\leq n\leq N-1$.
Finally we deduce the set of Fourier boundaries $\Omega=\{\omega^n\}_{n=0,\ldots,N}$ (the solid vertical lines in the figure) by ($\omega_0=0$):
\begin{equation}\label{eq:middle}
\omega^0=0\quad ;\quad\omega^N=\pi ; \qquad\omega^n=\frac{\omega_{n}+\omega_{n-1}}{2} \quad\text{for}\;\;1\leq n\leq N-1
\end{equation}
This simple method is efficient when the spectrum is composed of relatively well separated modes but can provide unexpected segmentation in some specific, but frequent in practice, cases. For instance, let us consider two 
close consecutive modes where one has a wide support while the other has a narrow support (see Figure.~\ref{fig:locmax}.b). It is straightforward to see that the corresponding boundary $\omega^n$ obtained from 
the above process will fall in the largest support of the first mode. This implies that the dashed area on the figure will be considered in the second mode while it is obviously part of the first mode. In fact a ``good'' 
segmentation should correspond the dot-dashed vertical line instead of the solid vertical line. A second example where this method does not provide the best segmentation is when several local maxima belong 
to the same mode and are larger than other modes. Figure~\ref{fig:locmax}.c illustrates this case. Imagine we want three bands, the square boxes show the selected local maxima and the solid vertical lines represent the found 
Fourier boundaries. Conceptually, it will be more reasonable to consider that $\omega_1$ and $\omega_2$ are part of the same mode and then keep $\omega_3$ instead of $\omega_2$ (giving the dot-dashed vertical line as 
the second Fourier boundary and hence well separating to two modes). In this case, the problem comes from the fact that the segmentation method considers only a local information but it should be better to also take into  
account the spectrum global trend to avoid such issues.

\subsubsection{Finding the lowest minima}
The issue depicted on Figure~\ref{fig:locmax}.b is due to the fact that the boundary between two consecutive detected maxima is computed as the point at equal distance from the two maxima (Eq.~\ref{eq:middle}). A 
simple way to avoid such situation is to retain the position of the lowest minima in the segment defined by consecutive maxima (this correspond to find the global minimum in this segment). If we denote $\mho_n$ the set 
of all local minima located between $\omega_{n-1}$ and $\omega_n$, then
\begin{equation}\label{eq:lowest}
\omega^0=0\quad ;\quad\omega^N=\pi ; \qquad\omega^n=\underset{\omega}{\arg}\min\mho_n \quad\text{for}\;\;1\leq n\leq N-1
\end{equation}
For instance, the use of this method will automatically provide the dot-dashed boundary of Figure.~\ref{fig:locmax}.b.

\subsubsection{Global trend removing}
\begin{figure}
\centering\includegraphics[width=0.46\textwidth]{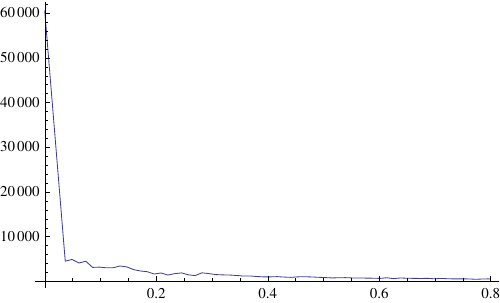}
\caption{Example of a profile extracted from the spectrum of a natural image.}
\label{fig:zoomlenaspect}
\end{figure}

The example showed on Figure.~\ref{fig:locmax}.c illustrates that many local maxima can belong to a common mode and may induce many unexpected boundaries. Such situation is quite common when we look at the spectrum profile 
of an image where the low frequencies concentrate a lot of the energy (see Figure.~\ref{fig:zoomlenaspect}) with many local maxima but belonging to the same mode. In order to get rid of this situation, it should be 
interesting to work with the logarithm of $H$ and/or to remove the global trend (denoted $T(\omega)$) of the analyzed spectrum before the detection task. Different approaches can be used to find the global trend of $H$. In this paper 
we propose different options:
\begin{itemize}
\item ``\textit{plaw}'': we approximate $H$ with a power law of the form $T(\omega)=\omega^{-s}$. The exponent $s$ can be estimated by using a least mean square criteria. In its discrete version it is equivalent to: 
\begin{equation}
s=\arg\min \|H(w)-w^{-s}\|_2=-\frac{\sum_n\ln\omega_n\ln H(\omega_n)}{\sum_n(\ln\omega_n)^2}.
\end{equation}
\item ``\textit{poly}'': we approximate $H$ with a polynomial of degree $d$. Usual methods (mean square error) can be used to compute the polynomial coefficients. Concerning the polynomial degree, there is no automatic 
way to choose it and we will consider it as a parameter of this method. In practice, experiments seem to indicate that $d\approx 5$ works for usual images.
\item ``\textit{morpho}'': this approach is inspired from the mathematical morphology \cite{serra1982,Soille1999}, see appendix~\ref{app:morpho} for a review of the mathematical morphology operators used in this paper. 
The opening operator ($Op$) provides a lower envelope of $H$ while the closing operator ($Cl$) provides an upper 
envelope. Then the global trend can be obtained by 
\begin{equation}
T(\omega)=\frac{Op(H,b)(\omega)+Cl(H,b)(\omega)}{2}.
\end{equation}
Concerning the size of the structuring function $b$, because we want to remove the quickly changes in $H$, we can choose a size equal to the smallest distance between two consecutive 
maxima.
\item ``\textit{tophat}'': this method, also based on mathematical morphology, consists to choose $T(\omega)=Op(H,b)(\omega)$. The size of $b$ is chosen like in the \textit{morpho} case.

\end{itemize}

\subsection{Fine to Coarse histogram segmentation algorithm}
In \cite{Delon2007}, the authors propose a fully automatic algorithm to detect the modes in an histogram. In our context, $H$ can be viewed as an histogram 
representing how the energy is distributed among the frequencies. Then we can use this algorithm on $H$. The main advantage of this method is in that it automatically selects the number of modes. It is based on a fine to coarse 
segmentation of $H$. The algorithm is initialized with all local minima of $H$. The authors defined a statistical criteria to decide if two consecutive supports correspond to a common global trend or if they are 
part of true separated modes; this criteria is based on the $\epsilon-$meaningful events theory \cite{Delon2007,Desolneux2008}. The (parameterless) algorithm can be resumed as following: start from the fine segmentation 
given by all local minima, choose one minimum and check if adjacent supports are part of a same trend or not. 
If yes then merge these supports by removing this local minima from the list. Repeat until no merging are possible.

\subsection{Comparison of the different approaches}
\begin{figure}[!h]
\begin{center}
\begin{tabular}{c|c|c}
 & ``middle point'' detection & ``lowest minima'' detection \\ \hline\hline
No preprocessing & \includegraphics[width=0.3\textwidth]{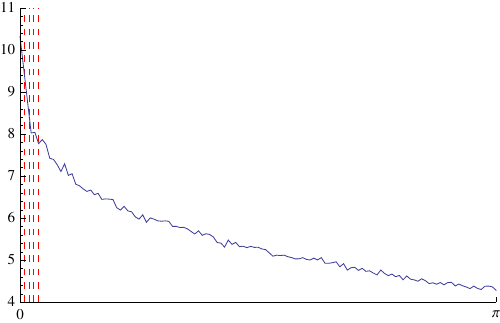} & \includegraphics[width=0.3\textwidth]{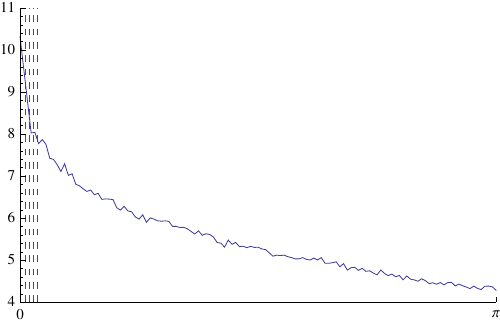} \\
%No preprocessing & No preprocessing \\ \hline
\textit{plaw} & \includegraphics[width=0.3\textwidth]{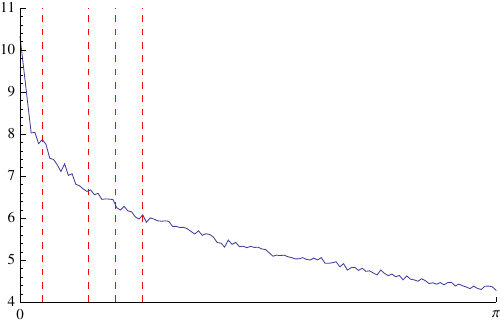} & \includegraphics[width=0.3\textwidth]{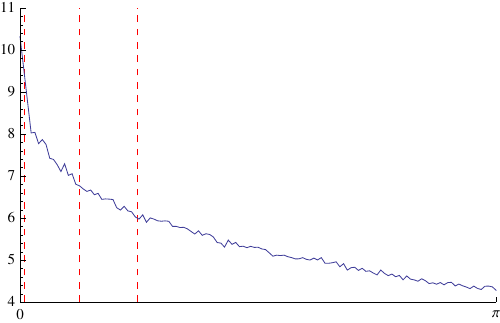} \\
%\textit{plaw} & \textit{plaw} \\ \hline
\textit{poly} & \includegraphics[width=0.3\textwidth]{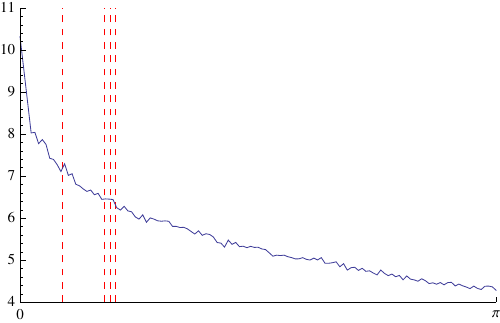} & \includegraphics[width=0.3\textwidth]{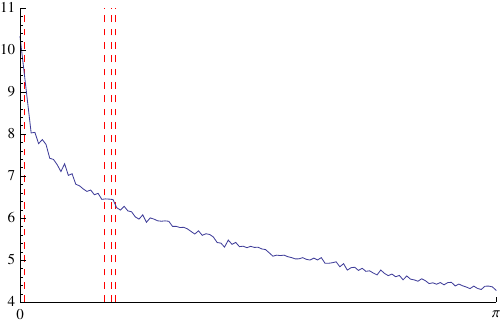} \\
%\textit{poly} & \textit{poly} \\ \hline
\textit{morpho} & \includegraphics[width=0.3\textwidth]{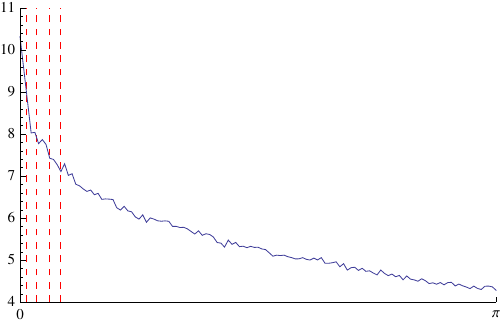} & \includegraphics[width=0.3\textwidth]{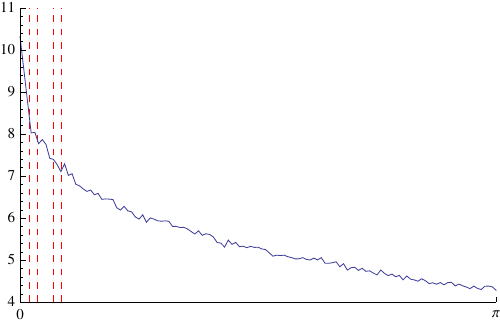} \\
%\textit{morpho} & \textit{morpho} \\  \hline
\textit{tophat} & \includegraphics[width=0.3\textwidth]{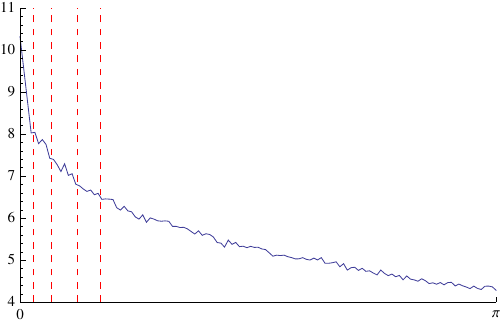} & \includegraphics[width=0.3\textwidth]{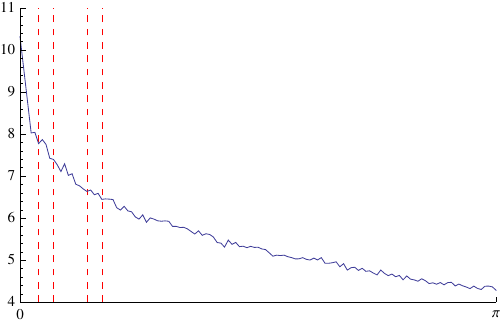} \\
%\textit{tophat} & \textit{tophat}
\end{tabular}
\end{center}
\caption{Fourier support detection on a spectrum with a global trend (without any logarithm and $N=5$). Each row correspond to different preprocessing: none, \textit{plaw}, \textit{poly} and \textit{morpho}, respectively. On the left column the 
detection based on the middle point between consecutive maxima is shown while the right column show the results when we keep the lowest minimum between consecutive maxima.}
\label{fig:boundlenanolog}
\end{figure}

\begin{figure}[!h]
\begin{center}
\begin{tabular}{c|c|c}
& ``middle point'' detection & ``lowest minima'' detection \\ \hline\hline
No preprocessing & \includegraphics[width=0.3\textwidth]{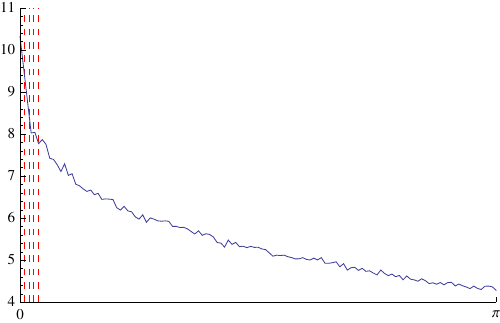} & \includegraphics[width=0.3\textwidth]{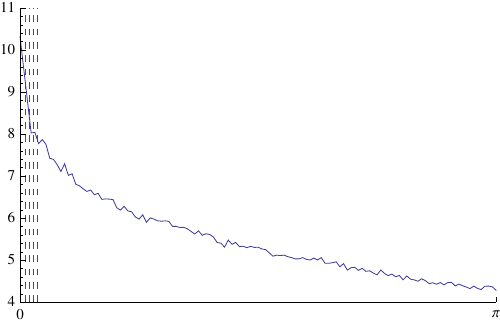} \\
%No preprocessing & No preprocessing \\ \hline
\textit{plaw} & \includegraphics[width=0.3\textwidth]{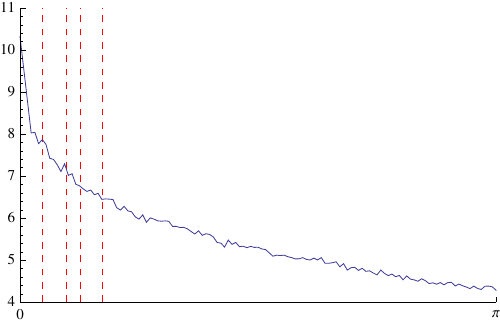} & \includegraphics[width=0.3\textwidth]{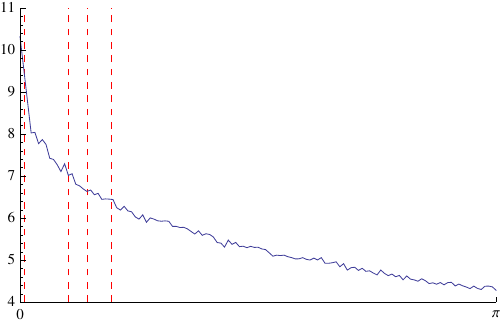} \\
%\textit{plaw} & \textit{plaw} \\ \hline
\textit{poly} & \includegraphics[width=0.3\textwidth]{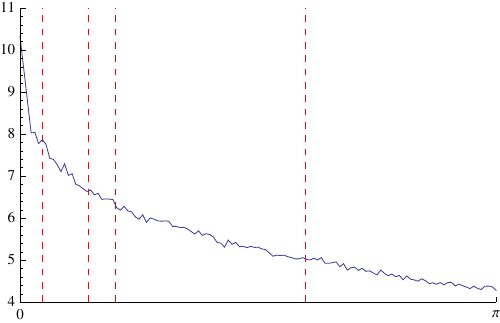} & \includegraphics[width=0.3\textwidth]{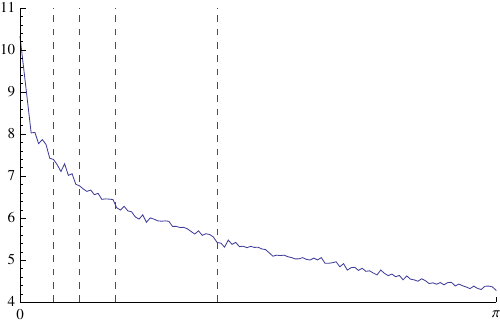} \\
%\textit{poly} & \textit{poly} \\ \hline
\textit{morpho} & \includegraphics[width=0.3\textwidth]{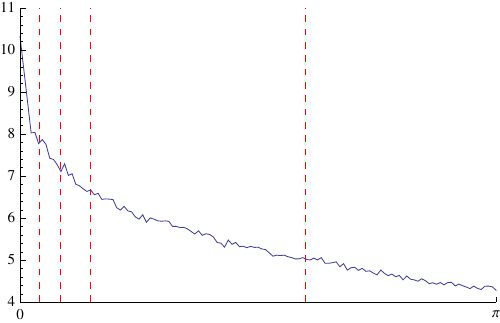} & \includegraphics[width=0.3\textwidth]{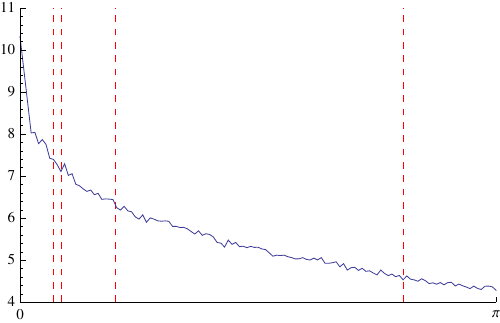} \\
%\textit{morpho} & \textit{morpho}  \\ \hline
\textit{tophat} & \includegraphics[width=0.3\textwidth]{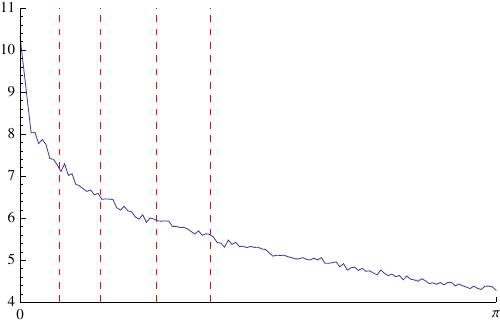} & \includegraphics[width=0.3\textwidth]{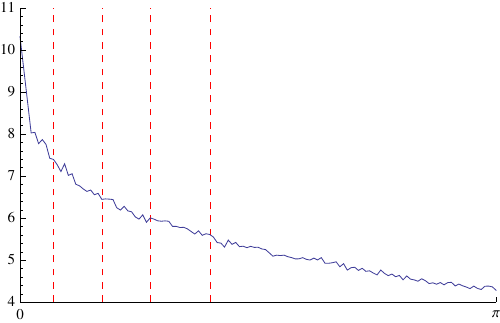} \\
%\textit{tophat} & \textit{tophat}
\end{tabular}
\end{center}
\caption{Fourier support detection on a spectrum with a global trend (with applying the logarithm and $N=5$). Each row correspond to different preprocessing: none, \textit{plaw}, \textit{poly} and \textit{morpho}, respectively. On the left column the 
detection based on the middle point between consecutive maxima is shown while the right column show the results when we keep the lowest minimum between consecutive maxima.}
\label{fig:boundlenalog}
\end{figure}

\begin{figure}[!h]
\begin{center}
\begin{tabular}{c|c|c}
& ``middle point'' detection & ``lowest minima'' detection \\ \hline\hline
No preprocessing & \includegraphics[width=0.28\textwidth]{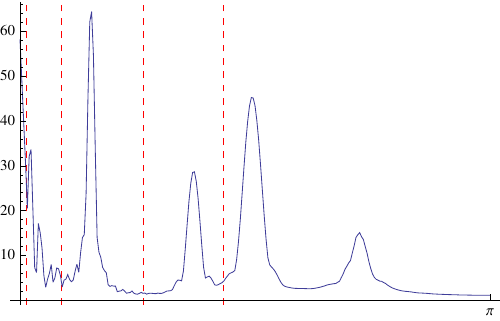} & \includegraphics[width=0.28\textwidth]{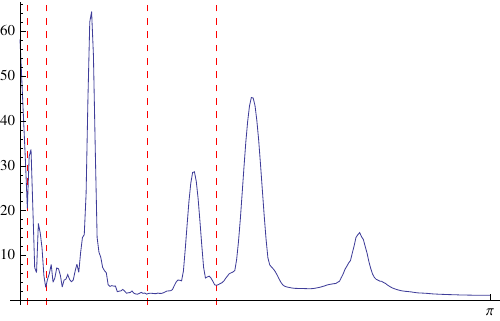} \\
%No preprocessing & No preprocessing \\ \hline
\textit{plaw} & \includegraphics[width=0.28\textwidth]{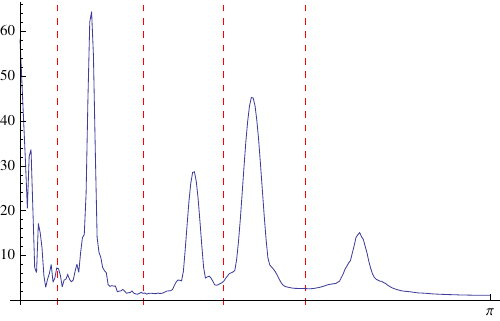} & \includegraphics[width=0.28\textwidth]{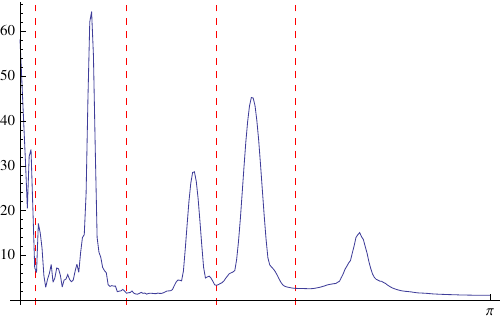} \\
%\textit{plaw} & \textit{plaw} \\ \hline
\textit{poly} & \includegraphics[width=0.28\textwidth]{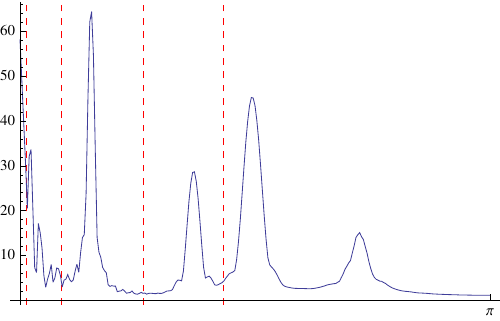} & \includegraphics[width=0.28\textwidth]{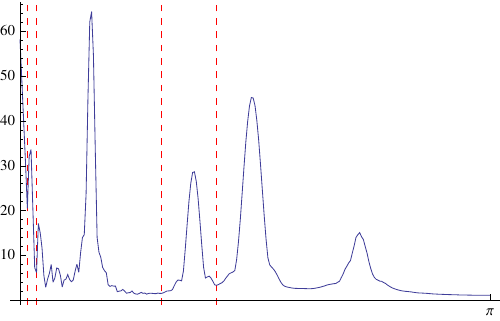} \\
%\textit{poly} & \textit{poly} \\ \hline
\textit{morpho} & \includegraphics[width=0.28\textwidth]{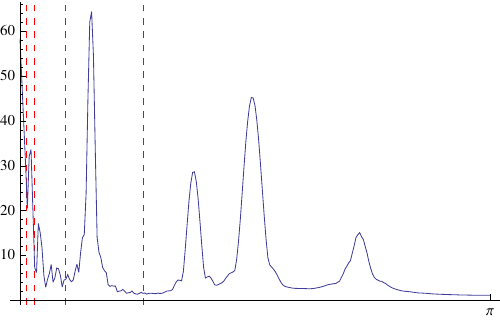} & \includegraphics[width=0.28\textwidth]{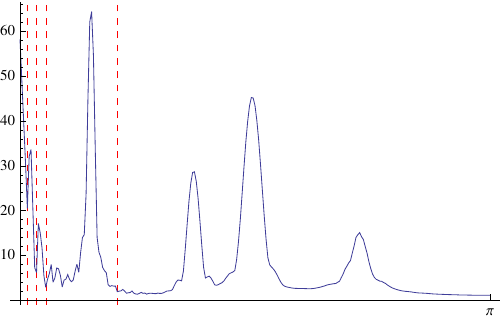} \\
%\textit{morpho} & \textit{morpho} \\ \hline
\textit{tophat} & \includegraphics[width=0.28\textwidth]{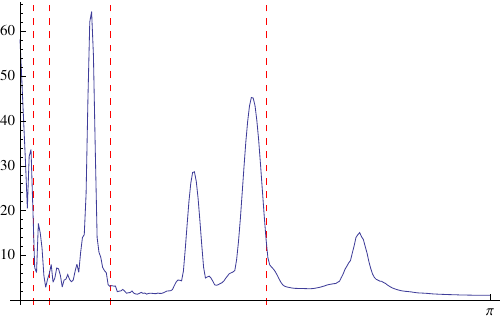} & \includegraphics[width=0.28\textwidth]{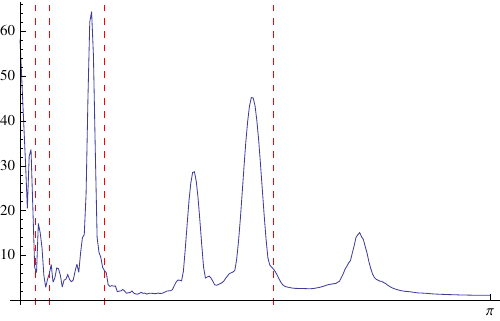} \\
%\textit{tophat} & \textit{tophat}
\end{tabular}
\end{center}
\caption{Fourier support detection on a spectrum with distinct modes (without any logarithm and $N=5$). Each row correspond to different preprocessing: none, \textit{plaw}, \textit{poly} and \textit{morpho}, respectively. On the left column the 
detection based on the middle point between consecutive maxima is shown while the right column show the results when we keep the lowest minimum between consecutive maxima.}
\label{fig:boundtexturesnolog}
\end{figure}

\begin{figure}[!h]
\begin{center}
\begin{tabular}{c|c|c}
& ``middle point'' detection & ``lowest minima'' detection \\ \hline\hline
No preprocessing & \includegraphics[width=0.28\textwidth]{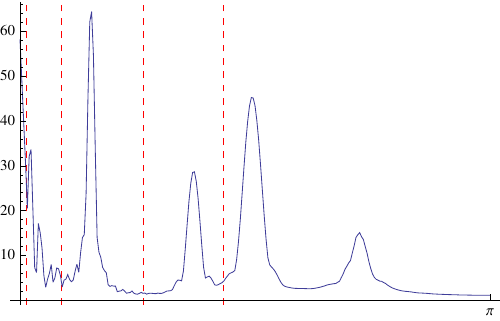} & \includegraphics[width=0.28\textwidth]{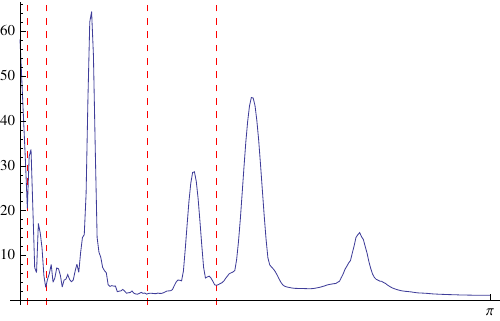} \\
%No preprocessing & No preprocessing \\ \hline
\textit{plaw} & \includegraphics[width=0.28\textwidth]{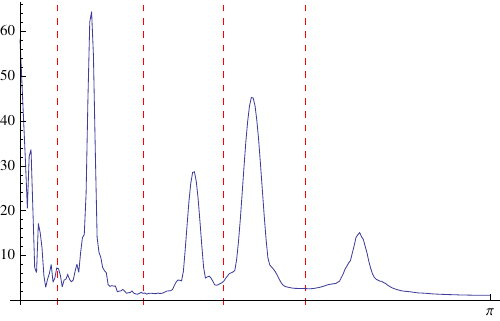} & \includegraphics[width=0.28\textwidth]{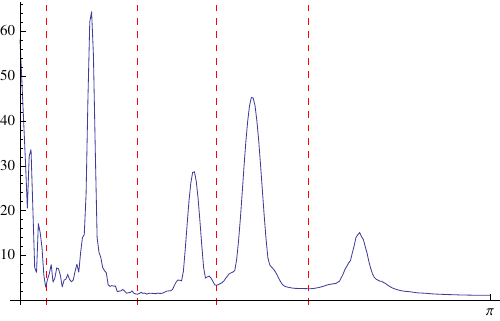} \\
%\textit{plaw} & \textit{plaw} \\ \hline
\textit{poly} & \includegraphics[width=0.28\textwidth]{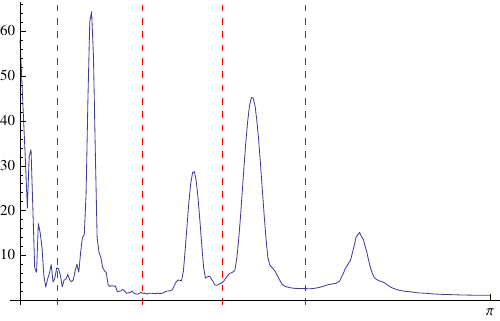} & \includegraphics[width=0.28\textwidth]{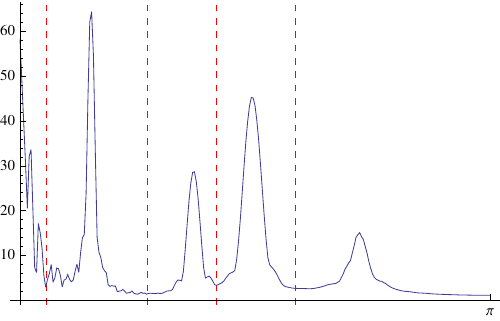} \\
%\textit{poly} & \textit{poly} \\ \hline
\textit{morpho} & \includegraphics[width=0.28\textwidth]{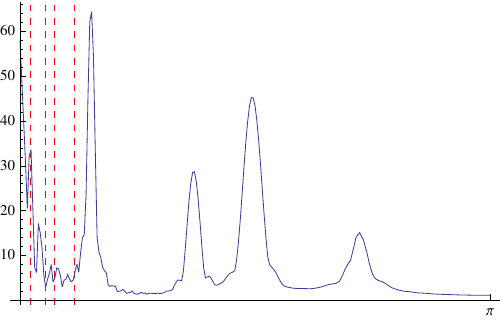} & \includegraphics[width=0.28\textwidth]{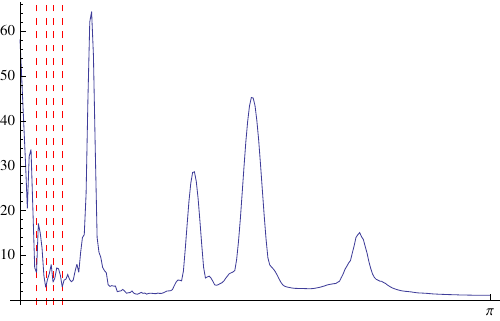} \\
%\textit{morpho} & \textit{morpho} \\ \hline
\textit{tophat} & \includegraphics[width=0.28\textwidth]{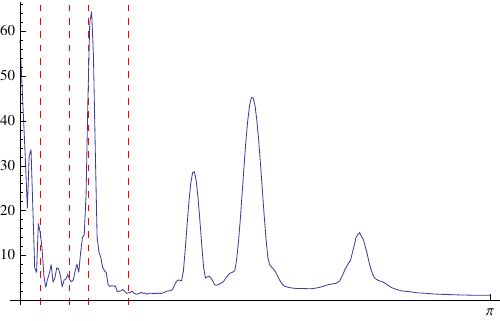} & \includegraphics[width=0.28\textwidth]{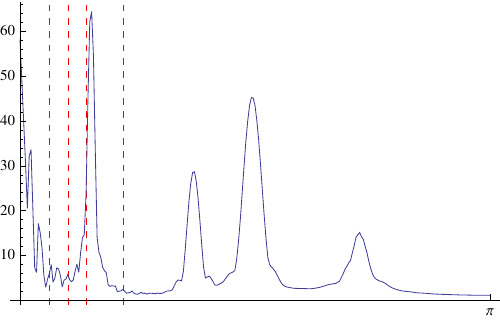} \\
%\textit{tophat} & \textit{tophat}
\end{tabular}
\end{center}
\caption{Fourier support detection on a spectrum with distinct modes (with applying the logarithm and $N=5$). Each row correspond to different preprocessing: none, \textit{plaw}, \textit{poly} and \textit{morpho}, respectively. On the left column the 
detection based on the middle point between consecutive maxima is shown while the right column show the results when we keep the lowest minimum between consecutive maxima.}
\label{fig:boundtextureslog}
\end{figure}

\begin{figure}[!h]
\begin{center}
\begin{tabular}{c|c}
\includegraphics[width=0.46\textwidth]{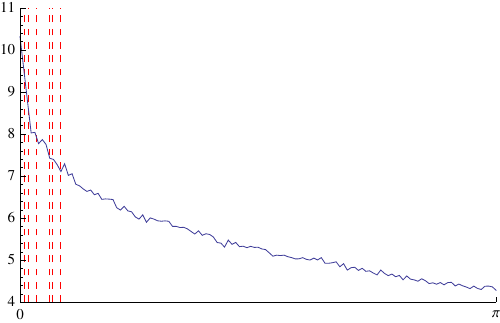} & \includegraphics[width=0.46\textwidth]{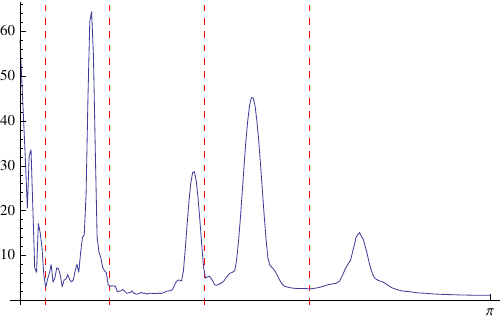}
\end{tabular}
\end{center}
\caption{Fourier support detection using the Fine to Coarse algorithm.}
\label{fig:boundftc}
\end{figure}

In order to understand the advantages and disadvantages of each method, we test each combination on two types of spectrum. The first one has a global trend with a huge magnitude difference between the low frequencies 
and the rest of the spectrum. The second one is mainly composed of distinct modes with a fast decaying global trend localized in the low frequencies. The respective detected Fourier boundaries are 
presented in Figures.~\ref{fig:boundlenanolog}, \ref{fig:boundlenalog}, \ref{fig:boundtexturesnolog} and \ref{fig:boundtextureslog}, respectively. We can observe in the case of the spectrum with a global 
trend that it is necessary to do a preprocessing in order to magnify some ``small'' modes. Moreover, we see that by using the logarithm of the spectrum, all preprocessing methods provide interesting detections (notice 
that the boundary closest to $\pi$ detected by the \textit{morpho} method is not an artifact but a very small mode, too small to be visible on the provided graphs). 
In the case of the second spectrum with distinct modes, there are no specific advantages to take the logarithm. If the detection without any preprocessing gives interesting results, we can see that the use of the \textit{plaw} 
method permits to deal with the fast decaying trend. The results of the Fine to Coarse algorithm are given in Figure.~\ref{fig:boundftc}. If this algorithm performs well in the case of distinct modes, it looks like it 
over segments the spectrum with the global trend.

\section{Tensor 2D EWT}\label{sec:tensor}
In \cite{Gilles2013}, the author proposed to extend the 1D EWT to images by considering a tensor product approach as with the usual 2D discrete wavelet transform. The idea is to use separately the 1D EWT on rows and 
columns, respectively. It is easy to see that if we consider each row (or column) independently we can obtain different sets of filters defined on very different Fourier supports. Indeed, for instance, let us consider two 
different rows in an image where one mainly contains edges while the other reflects the presence of textures. Their corresponding spectrum will be very different, hence providing different Fourier supports. It is possible 
that a scale $n$ in one case corresponds to a completely different spectral information as opposed to the same scale from the other row. Such ``discontinuity''in the filter bank leads to discontinuities in the EWT components. 
In order to avoid such behavior, 
the author proposed to use the same filter bank to process all rows and a second one to process all columns. These filters are generated by using the following procedure: for the rows, we compute the 1D FFT of each row and evaluate the average magnitude 
spectrum; then the Fourier supports detection is made on this average spectrum. Those supports permit us to build a unique filter bank which is used to process each 
row. The column case is completely equivalent by computing a column average magnitude spectrum. If we denote $N_{row}$ and $N_{col}$ the number of rows and columns, respectively, the 2D Tensor EWT can be resumed by 
Algorithm.~\ref{algo:tensor}.

\begin{algorithm}
{\bf Inputs}: image $f(\x)$, $N_R$ (number of filters for the rows) and $N_C$ (number of filters for the columns)\;
For each row $i$, compute its 1D FFT; $\F_{1,x_1}(f)(i,\omega_2)$; then compute the average row spectrum magnitude: $$\widetilde{\F}_{row}(\omega_2)=\frac{1}{N_{row}}\sum_{i=0}^{N_{row}-1}|\F_{1,x_1}(f)(i,\omega_2)|;$$\\
For each column $j$, compute its 1D FFT; $\F_{1,x_2}(f)(\omega_1,j)$; then compute the average row spectrum magnitude: $$\widetilde{\F}_{col}(\omega_1)=\frac{1}{N_{col}}\sum_{j=0}^{N_{col}-1}|\F_{1,x_2}(f)(\omega_1,j)|;$$\\
Perform the Fourier boundaries detection on $\widetilde{\F}_{row}$ to get $\Omega_{row}$ and build $\B_{row}=\{\phi_1^{row},\{\psi_n^{row}\}_{n=1}^{N_R-1}\}$ accordingly to equations \eqref{eq:phi2} and \eqref{eq:psi2}\;
Perform the Fourier boundaries detection on $\widetilde{\F}_{col}$ to get $\Omega_{col}$ and build $\B_{col}=\{\phi_1^{col},\{\psi_m^{col}\}_{m=1}^{N_C-1}\}$ accordingly to equations \eqref{eq:phi2} and \eqref{eq:psi2}\;
Filter $f$ along the rows with $\B_{row}$, this provides $N_R$ output images\;
Filter each previous output images along the columns with $\B_{col}$\;
{\bf Outputs}: $\B_{row},\B_{col},\W^{\E\mathcal{T}}(n,m,\x)$
\caption{2D Tensor Empirical Wavelet Transform algorithm.}
\label{algo:tensor}
\end{algorithm}
At the end, this transform $\W^{\E\mathcal{T}}(n,m,\x)$ provides a set of $N_R\times N_C$ subband images corresponding to a projection on the 2D frame 
(adopting the convention $\psi_0=\phi_1$):
\begin{equation}
\B^{\E\mathcal{T}}=\left\{\psi_{nm}(\x)=\psi_n^{row}(x_1)\psi_m^{col}(x_2)\right\}_{\substack{n=0,\ldots,N_R-1 \\ m=0,\ldots,N_C-1}}.
\end{equation}
Another point of view is that this procedure builds functions $\psi_{nm}(\x)$ which have adaptive rectangular supports in the Fourier plane (this will be illustrated in the experiments section). The inverse transform can be obtained by using the adjoint formulation (e.g applying 
the inverse 1D EWT with respect to the columns first and then with respect to the rows).

\section{2D Littlewood-Paley EWT}\label{sec:lpewt}
The classic 2D Littlewood-Paley wavelet transform corresponds to filter images with 2D wavelets defined in the Fourier domain on annuli supports (centered around the origin) \cite{grafakos2008classical}. The inner and outer radius 
of these supports are fixed upon a dyadic decomposition of the Fourier plane (corresponding to the usual notion of scales). It is easy to build an image for which its Fourier energy is spread into two consecutive 
supports leading to the separation of this information after the wavelet filtering. In this paper, we propose to apply the empirical approach to detect the radius of each annuli. The best point of view to 
perform this detection is to consider the Fourier plane in a polar representation because finding such boundaries is equivalent to working with the frequency modulus $|\omega|$. Some work has been made on the construction 
of a Pseudo-Polar FFT \cite{PseudoPolar,PseudoPolar2} providing an operator $\F_P(f)(\theta,|\omega|)$. For each angle $\theta$ we have a 1D Fourier spectrum but, as in the tensor transform case, if we perform the Fourier 
boundaries detection for each $\theta$ independently we will introduce some discontinuities in the output components. In order to avoid such artifacts, we adapt the idea used in the tensor transform by computing an average spectrum 
where the averaging is taken with respect to $\theta$: $\widetilde{\F}(|\omega|)=\frac{1}{N_{\theta}}\sum_{i=0}^{N_{\theta}-1}\F_P(f)(\theta_i,|\omega|)$, where $N_{\theta}$ is the number of discrete angles. Then we 
perform the Fourier boundaries detection on $\widetilde{\F}(|\omega|)$ and get the set of spectral radius, denoted $\{\omega^n\}_{n=0,\ldots,N}$ (with $\omega^0=0$ and $\omega^N=\pi$), which we can use to build a set of 2D 
empirical Littlewood-Paley wavelets $\B^{\E\mathcal{LP}}=\left\{\phi_1(\x),\{\psi_n(\x)\}_{n=1}^{N-1}\right\}$. Their definition is a straightforward extension of equations \eqref{eq:phi2} and \eqref{eq:psi2}, except for the last annuli (for $\omega^{N-1}\leq |\omega|
\leq\omega^N=\pi$ where we extend the ring in order to keep the ``corners'' of the Fourier plane):
\begin{equation}\label{eq:lpphi2}
\F_2(\phi_1)(\omega)=
\begin{cases}
1  \; &\text{if}\;|\omega|\leq (1-\gamma)\omega^1\\
\cos\left[\frac{\pi}{2}\beta\left(\frac{1}{2\gamma\omega^1}(|\omega|-(1-\gamma)\omega^1)\right)\right] \; &\text{if}\; (1-\gamma)\omega^1\leq |\omega|\leq (1+\gamma)\omega^1 \\
0  \; &\text{otherwise,}
\end{cases}
\end{equation}
and, if $n\neq N-1$:
\begin{equation}\label{eq:lppsi2a}
\F_2(\psi_n)(\omega)=
\begin{cases}
1 \qquad\qquad\qquad\qquad\qquad\quad\;\;\; \text{if}\; (1+\gamma)\omega^n\leq |\omega|\leq (1-\gamma)\omega^{n+1} \\
\cos\left[\frac{\pi}{2}\beta\left(\frac{1}{2\gamma\omega^{n+1}}(|\omega|-(1-\gamma) \omega^{n+1})\right)\right]  \\
\qquad\qquad\qquad\qquad\qquad\qquad\; \text{if}\; (1-\gamma)\omega^{n+1}\leq |\omega|\leq (1+\gamma)\omega^{n+1}\\
\sin\left[\frac{\pi}{2}\beta\left(\frac{1}{2\gamma\omega^n}(|\omega|-(1-\gamma)\omega^n)\right)\right] \\
\qquad\qquad\qquad\qquad\qquad\qquad\;  \text{if}\; (1-\gamma)\omega^n\leq |\omega|\leq (1+\gamma)\omega^n\\
0 \qquad\qquad\qquad\qquad\qquad\quad\;\;\; \text{otherwise,}
\end{cases}
\end{equation}
and, if $n=N-1$:
\begin{equation}\label{eq:lppsi2b}
\F_2(\psi_{N-1})(\omega)=
\begin{cases}
1 \qquad\qquad\qquad\qquad\qquad\quad\;\;\; \text{if}\; (1+\gamma)\omega^{N-1}\leq |\omega| \\
\sin\left[\frac{\pi}{2}\beta\left(\frac{1}{2\gamma\omega^{N-1}}(|\omega|-(1-\gamma)\omega^{N-1})\right)\right] \\
\qquad\qquad\qquad\qquad\qquad\qquad\;  \text{if}\; (1-\gamma)\omega^{N-1}\leq |\omega|\leq (1+\gamma)\omega^{N-1}\\
0 \qquad\qquad\qquad\qquad\qquad\quad\;\;\; \text{otherwise.}
\end{cases}
\end{equation}
Then the 2D empirical Littlewood-Paley transform of an input image $f$ is given by
\begin{equation}\label{eq:lpta}
\W_f^{\E\mathcal{LP}}(n,\x)=\F_2^*\left(\F_2(f)(\omega)\overline{\F_2(\psi_n)(\omega)}\right),
\end{equation}
for the detail coefficients and the approximation coefficients (the convention $\mathcal{W}_f^{\E\mathcal{LP}}(0,\x)$ is adopted to denote them) by:
\begin{equation}\label{eq:lptb}
\W_f^{\E\mathcal{LP}}(0,\x)=\F_2^*\left(\F_2(f)(\omega)\overline{\F_2(\phi_1)(\omega)}\right).
\end{equation}
The corresponding algorithm is resumed in Algorithm.~\ref{algo:lp}. The inverse transform is obtained via the adjoint formulation:
\begin{equation}
f(\x)=\F_2^*\left(\F_2\left(\W_f^{\E\mathcal{LP}}\right)(0,\omega)\F_2(\phi_1)(\omega)+\sum_{n=1}^{N-1}\F_2\left(\W_f^{\E\mathcal{LP}}\right)(n,\omega)\F_2(\psi_n)(\omega)\right).
\end{equation}

\begin{algorithm}
{\bf Inputs}: image $f(\x)$, $N$ (number of filters)\;
Compute the Pseudo-Polar FFT, $\F_P(f)(\theta,|\omega|)$, and take the average with respect to $\theta$: 
\begin{equation}\label{eq:meanppfft}
\widetilde{\F_P}(|\omega|)=\frac{1}{N_{\theta}}\sum_{i=0}^{N_{\theta}-1}|\F_P(f)(\theta_i,|\omega|)|;
\end{equation}\\
Perform the Fourier boundaries detection on $\widetilde{\F_P}(|\omega|)$ to get $\Omega$ and build the corresponding filter bank $\B=\left\{\phi_1(\x),\{\psi_n(\x)\}_{n=1}^{N-1}\right\}$ accordingly to equations \eqref{eq:lpphi2}-\eqref{eq:lppsi2b}\;
Filter $f$ by using equations \eqref{eq:lpta} and \eqref{eq:lptb}\;
{\bf Outputs}: $\B^{\E\mathcal{LP}}$ and $\W_f^{\E\mathcal{LP}}(n,\x)$
\caption{2D Empirical Littlewood-Paley Wavelet Transform algorithm.}
\label{algo:lp}
\end{algorithm}

\section{Empirical Ridgelet Transform}\label{sec:eridge}
The ridgelet transform was introduced by Cand\`es and Donoho \cite{candes1,candes2,Donoho98} and was one of the first directional 2D wavelet type transforms. The idea is to consider the usual wavelet transform on 
a collection of 1D signals. This 1D signal collection is generated by an inverse 1D FFT of the 1D spectrum profile extracted from the 2D FFT of the input image along lines going through the origin and with angle $\theta$.  
By using our notations, the classical ridgelet transform 
can be written, $0<n<N$ (where $n=0$ corresponds to the approximation subband),
\begin{equation}
\W_f^{\mathcal{R}}(n,\theta,t)=\W_{1,t}(\F_{1,|\omega|}^*(\F_P(f)(\theta,|\omega|))),
\end{equation}
and its inverse
\begin{equation}
f(\x)=\F_P^*(\F_{1,t}(\W_{1,t}^*(\W_f^{\mathcal{R}})).
\end{equation}
We can define the empirical ridgelet transform by the same algorithm except that we replace the standard 1D wavelet transform by the 1D empirical wavelet transform.
Note that in this case the number of subbands is supposed to depend on $\theta$ which means that we will face the same discontinuity problems as for the tensor or Littlewood-Paley approaches. 
One more time, in order to avoid such issue, we will use the same method as for the Littlewood-Paley algorithm: we will detect the Fourier boundaries on an average spectrum (with respect to $\theta$) obtained from 
the Pseudo-Polar FFT. Then we build the set $\B^{\E\mathcal{R}}=\left\{\phi_1(t),\{\psi_n(t)\}_{n=1}^{N-1}\right\}$ which will be used to perform the 1D EWT. Hence the empirical ridgelet transform is defined by ($\W_t^{\E}$ is the 
1D EWT based on $\B$ with respect to $t$)
\begin{equation}
\W_f^{\E\mathcal{R}}(n,\theta,t)=\W_t^{\E}(\F_{1,\omega}^*(\F_P(f))),
\end{equation}
and its inverse
\begin{equation}
f(\x)=\F_P^*(\F_{1,t}(\W_t^{\E*}(\W_f^{\E\mathcal{R}})).
\end{equation}
Because the empirical wavelet is defined in the Fourier domain, we can simplify these equations to
\begin{equation}\label{eq:erta}
\W_f^{\E\mathcal{R}}(n,\theta,t)=\F_{1,\omega}^*\left(\F_P(f)(\omega,\theta)\overline{\F_{1,t}(\psi_n)(\omega)}\right),
\end{equation}
\begin{equation}\label{eq:ertb}
\W_f^{\E\mathcal{R}}(0,\theta,t)=\F_{1,\omega}^*\left(\F_P(f)(\omega,\theta)\overline{\F_{1,t}(\phi_1)(\omega)}\right),
\end{equation}
and its inverse
\begin{equation}
f(\x)=\F_P^*\left(\F_{1,t}\left(\W_f^{\E\mathcal{R}}(0,\theta,t)\right)\F_{1,t}(\phi_1)(\omega)+\sum_{n=1}^{N-1}\F_{1,t}\left(\W_f^{\E\mathcal{R}}(n,\theta,t)\right)\F_{1,t}(\psi_n)(\omega)\right).
\end{equation}
The corresponding algorithm is resumed in Algorithm.~\ref{algo:ridge}.

\begin{algorithm}
{\bf Inputs}: image $f(\x)$, $N$ (number of filters)\;
Compute the Pseudo-Polar FFT, $\F_P(f)(\theta,|\omega|)$, and take the average with respect to $\theta$: 
\begin{equation}
\widetilde{\F_P}(|\omega|)=\frac{1}{N_{\theta}}\sum_{i=0}^{N_{\theta}-1}|\F_P(f)(\theta_i,|\omega|)|;
\end{equation}\\
Perform the Fourier boundaries detection on $\widetilde{\F_P}(|\omega|)$ to get $\Omega$ and build the corresponding filter bank $\B^{\E\mathcal{R}}=\left\{\phi_1(\x),\{\psi_n(\x)\}_{n=1}^{N-1}\right\}$ accordingly to equations \eqref{eq:phi2}-\eqref{eq:psi2}\;
Filter $f$ by using equations \eqref{eq:erta} and \eqref{eq:ertb}\;
{\bf Outputs}: $\B^{\E\mathcal{R}}$ and $\W_f^{\E\mathcal{R}}(n,\theta,t)$
\caption{Empirical Ridgelet Transform algorithm.}
\label{algo:ridge}
\end{algorithm}

\section{Empirical Curvelet Transform}\label{sec:ecurv}
\begin{figure}[!t]
\begin{center}
\begin{tabular}{ccc}
\includegraphics[width=0.3\textwidth]{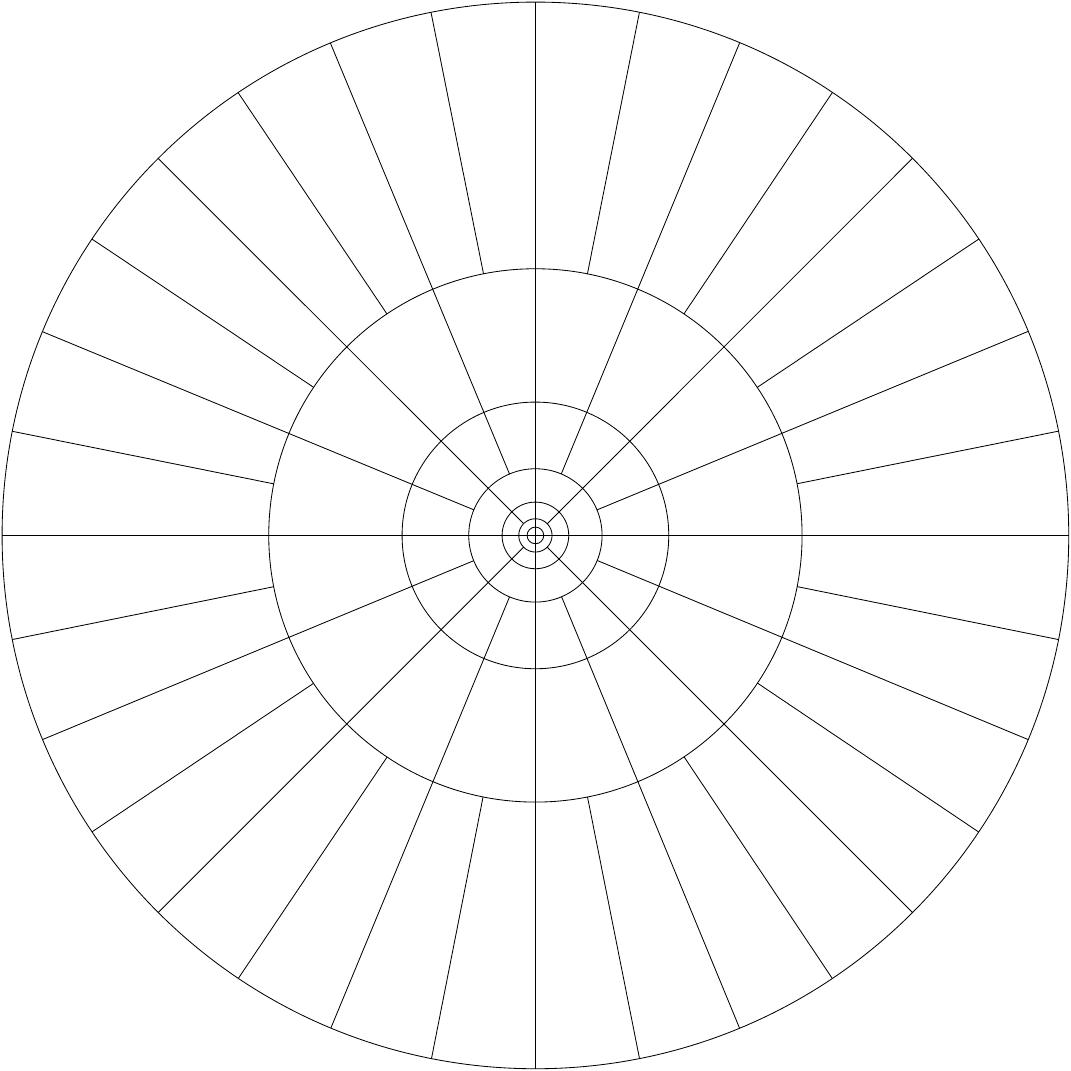} & \includegraphics[width=0.3\textwidth]{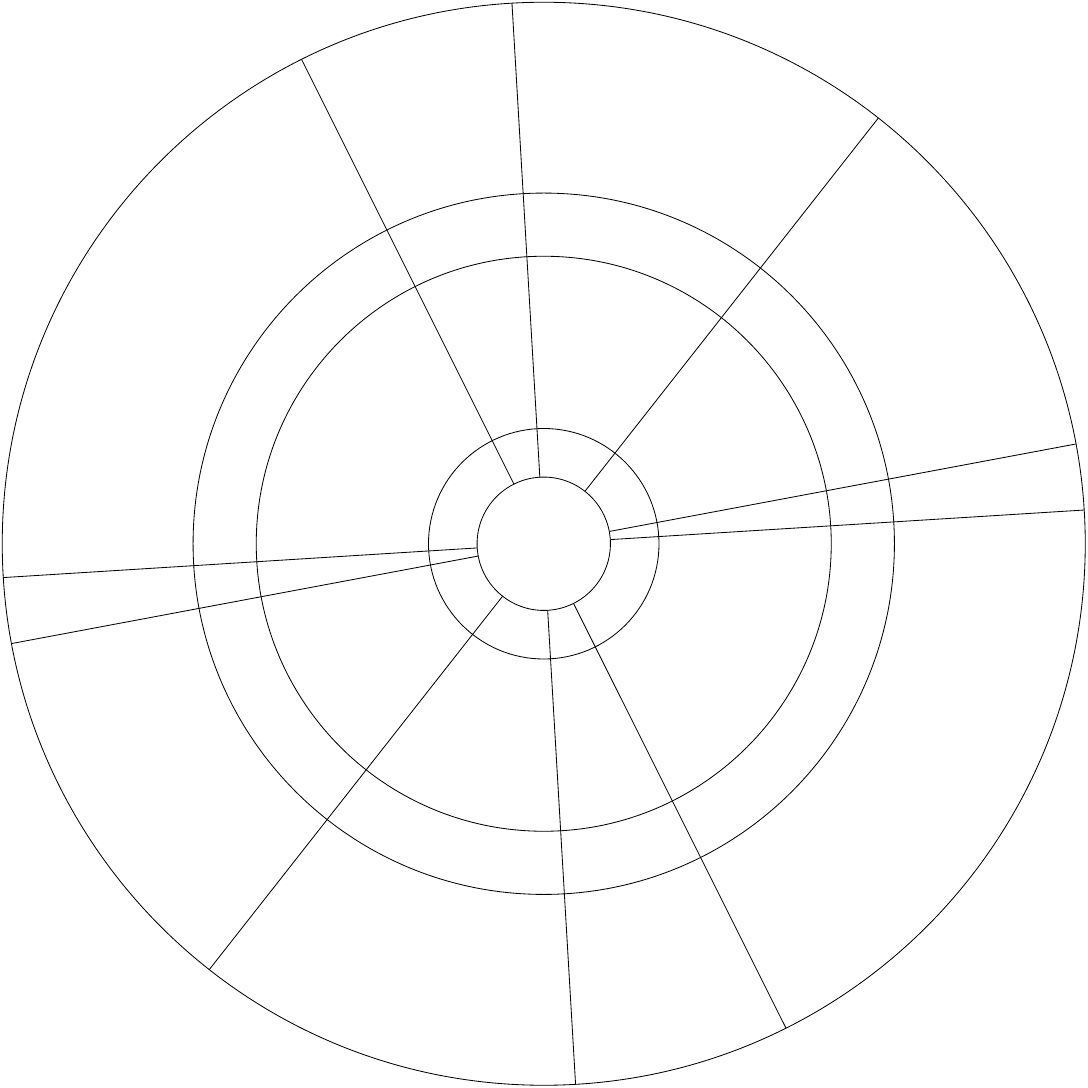} & \includegraphics[width=0.3\textwidth]{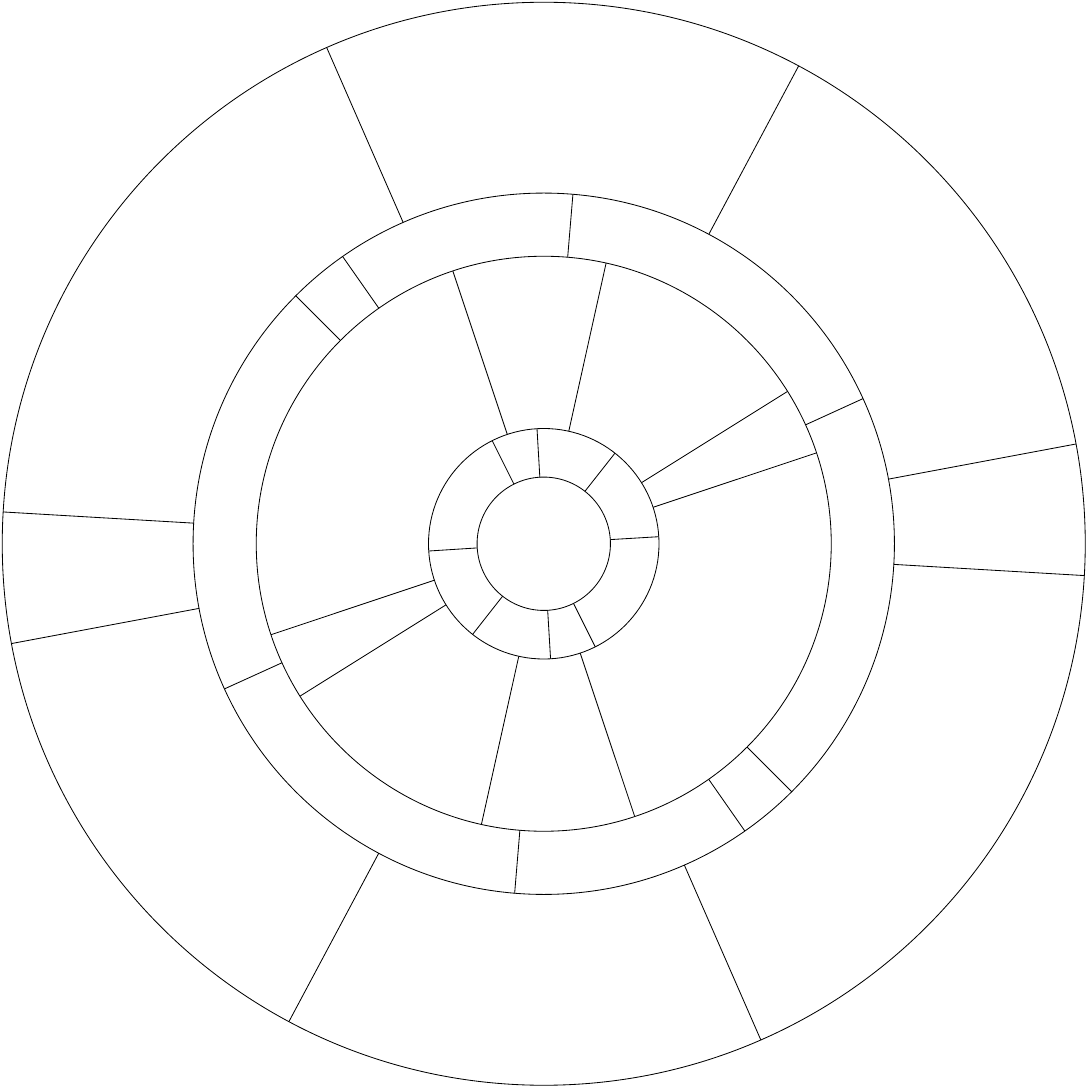} \\
Standard transform & Empirical transform-I & Empirical transform-II
\end{tabular}
\end{center}
\caption{Curvelet tiling of the Fourier plane.}
\label{fig:curvtiling}
\end{figure}

The curvelet transform was introduced by Cand\`es et al. in \cite{Candes3,Candes4,Candes5}. The concept behind this transform is to build a filter bank in the Fourier domain where each filter has its support located 
on a ``polar wedge''. The Fourier transform of one curvelet can be written ($(\omega,\theta)$ are the polar coordinates in the Fourier plane) as $\F_2(\psi_j)(\omega,\theta)=2^{-3j/4}W(2^{-j}\omega)V\left(\frac{2^{\lfloor j/2\rfloor}\theta}{2\pi}\right)$ where $W(r)$ and $V(t)$ (called 
the radial window and the angular window, respectively) are smooth, nonnegative and real-valued functions defined over compact supports (with $r\in(1/2,2)$ and $t\in[-1,1]$). A dyadic decomposition is used to partition 
the Fourier plane: the low frequencies are located on a disk centered at the zero frequency and each scale are defined on concentric annuli (in the same idea as the Littlewood-Paley transform), the angle ranges 
corresponding to each angular sector are also split into dyadic intervals and their number double every two scales (see left diagram on Figure.~\ref{fig:curvtiling}).
The empirical extension will consist in both empirically detect the scales and the angles corresponding to each polar wedge. We assume that the number of scales $N_s$ and the number of angular sectors $N_{\theta}$ are 
given. The detection process provides us with a set of scale boundaries $\Omega_{\omega}=\left\{\omega^n\right\}_{n=0,\ldots N_s}$ and a set of angular boundaries $\Omega_{\theta}=\left\{\theta^m\right\}_{m=1,\ldots,N_{\theta}}$. 
The (purely radial) lowpass filter $\phi_1$ is equivalent to the one used in the 
Littlewood-Paley transform and is given by Equation.~\ref{eq:lpphi2}. As with the standard curvelets, a polar wedge $\psi_{nm}$ ($n$ and $m$ are the scale and angular indices, respectively) in the Fourier domain is 
the product of a radial window, $W_n$, which for $n\neq N_s-1$ is given
\begin{equation}\label{eq:Wna}
W_n(\omega)=
\begin{cases}
1 \qquad\qquad\qquad\qquad\qquad\quad\;\;\; \text{if}\; (1+\gamma)\omega^n\leq |\omega|\leq (1-\gamma)\omega^{n+1} \\
\cos\left[\frac{\pi}{2}\beta\left(\frac{1}{2\gamma\omega^{n+1}}(|\omega|-(1-\gamma) \omega^{n+1})\right)\right]  \\
\qquad\qquad\qquad\qquad\qquad\qquad\; \text{if}\; (1-\gamma)\omega^{n+1}\leq |\omega|\leq (1+\gamma)\omega^{n+1}\\
\sin\left[\frac{\pi}{2}\beta\left(\frac{1}{2\gamma\omega^n}(|\omega|-(1-\gamma)\omega^n)\right)\right] \\
\qquad\qquad\qquad\qquad\qquad\qquad\;  \text{if}\; (1-\gamma)\omega^n\leq |\omega|\leq (1+\gamma)\omega^n\\
0 \qquad\qquad\qquad\qquad\qquad\quad\;\;\; \text{otherwise.}
\end{cases}
\end{equation}
and, if $n=N_s-1$:
\begin{equation}\label{eq:Wnb}
W_{N_s-1}(\omega)=
\begin{cases}
1 \qquad\qquad\qquad\qquad\qquad\quad\;\;\; \text{if}\; (1+\gamma)\omega^{N_s-1}\leq |\omega| \\
\sin\left[\frac{\pi}{2}\beta\left(\frac{1}{2\gamma\omega^{N_s-1}}(|\omega|-(1-\gamma)\omega^{N_s-1})\right)\right] \\
\qquad\qquad\qquad\qquad\qquad\qquad\;  \text{if}\; (1-\gamma)\omega^{N_s-1}\leq |\omega|\leq (1+\gamma)\omega^{N_s-1}\\
0 \qquad\qquad\qquad\qquad\qquad\quad\;\;\; \text{otherwise.}
\end{cases}
\end{equation}
and a polar window, $V_m$,
\begin{equation}\label{eq:Vnb}
V_m(\theta)=
\begin{cases}
1 \qquad\qquad\qquad\qquad\qquad\quad\;\;\; \text{if}\quad \theta^m+\Delta\theta\leq \theta\leq \theta^{m+1}-\Delta\theta \\
\cos\left[\frac{\pi}{2}\beta\left(\frac{1}{2\Delta\theta}(\theta-\theta^{m+1}+\Delta\theta)\right)\right]  \\
\qquad\qquad\qquad\qquad\qquad\qquad\; \text{if}\quad \theta^{m+1}-\Delta\theta\leq \theta\leq \theta^{m+1}+\Delta\theta\\
\sin\left[\frac{\pi}{2}\beta\left(\frac{1}{2\Delta\theta}(\theta-\theta^m+\Delta\theta)\right)\right] \\
\qquad\qquad\qquad\qquad\qquad\qquad\;  \text{if}\quad \theta^m-\Delta\theta\leq \theta\leq \theta^m+\Delta\theta\\
0 \qquad\qquad\qquad\qquad\qquad\quad\;\;\; \text{otherwise.}
\end{cases}
\end{equation}
where $\theta^{N_{\theta}+1}=\theta^1+\pi$ and $(\omega,\theta)$ are the polar coordinates in the Fourier plane. 
Moreover, in this paper, we only consider real transforms i.e. these filters are symmetric about the origin.
This leads us to the construction of the set 
$\B^{\E\mathcal{C}}_I=\left\{\phi_1(\x),\{\psi_{nm}(\x)\}_{\substack{n=1\ldots N_s-1\\ m=1\ldots N_{\theta}}}\right\}$ but other options are possible. Indeed the previous definition of $\B^{\E\mathcal{C}}_I$ consider 
the case where the scales and angles are detected independently. It is also possible to consider that one detection depends on the other. For instance, we can detect first the scales $\Omega_{\omega}$ and then 
detect different angles sets $\Omega_{\theta}^{\omega}$ for each scale (the other case is to detect first the angles and then the scales for each angular sector, but this case will lead to some information discontinuities 
from one angular sector to its neighbours and will not be considered in this paper). In this case we then build a set 
$\B^{\E\mathcal{C}}_{II}=\left\{\phi_1(\x),\{\psi_{nm}(\x)\}_{\substack{n=1\ldots N_s-1\\ m=1\ldots N_{\theta}}}\right\}$ where the angular window is dependent of $n$ because based on $\Omega_{\theta}^{\omega}$ instead of 
$\Omega_{\theta}$ (concretely, the bounds $\theta^m$ depend on $n$: $\theta^m_n$).

\section{Frame properties}
The following theorem gives the conditions such that all each previous transform are tight frames of $L^2$.
\begin{theorem}
Let $\B^{\E\mathcal{T}},\B^{\E\mathcal{LP}},\B^{\E\mathcal{R}},\B^{\E\mathcal{C}}_{I}$ and $\B^{\E\mathcal{C}}_{II}$ be the families of wavelets described in the previous sections. If their corresponding parameters 
$\gamma$ and $\Delta\theta$ (for the curvelet transforms) are chosen in order to avoid any overlapping of the transition areas (both for scales and angles) then these families are tight frames of $L^2$.
\end{theorem}
\begin{proof}
Let us start with the tensor product transform. Because of the tensor product properties we have $\F(\psi_{nm})(\omega)=\F(\psi_n^{row})(\omega_1)\F(\psi_n^{col})(\omega_2)$. 
Then 
\begin{align}
\sum_n\sum_m|\F(\psi_{nm})(\omega)|^2&=\sum_n\sum_m|\F(\psi_n^{row})(\omega_1)|^2|\F(\psi_n^{col})(\omega_2)|^2\\ \notag
&=\left(\sum_n|\F(\psi_n^{row})(\omega_1)|^2\right)\left(\sum_m|\F(\psi_n^{col})(\omega_2)|^2\right).    
\end{align}

But each term can be viewed as a 1D empirical wavelet set which we know from \cite{Gilles2013} that if $\gamma$ is properly chosen then $\sum_n|\F(\psi_n^{row})(\omega_1)|^2=\sum_n|\F(\psi_n^{col})(\omega_2)|^2=1$. 
Finally we get $\sum_n\sum_m=|\F(\psi_{nm})(\omega)|^2=1$ which end the tensor product case.\\
The Littlewood-Paley case is straightforward because it is a radial function defined from the 1D empirical wavelet which 
directly inherits its properties.\\
The ridgelet case is also straightforward as it corresponds to take 1D EWT along lines at different angles.\\
Let us address now the curvelet cases. We start from the first option when scales and angles are independently detected. In this case, we have 
$|\F(\phi_1)(\omega)|^2+\sum_n\sum_m|\F(\psi_{nm}(\omega))|^2=|\F(\phi_1)(\omega)|^2+\sum_n\sum_m|W_n(\omega)|^2|V_m(\theta)|^2=|\F(\phi_1)(\omega)|^2+\sum_n|W_n(\omega)|^2\sum_m|V_m(\theta)|^2$. 
But by construction (it is the same concept as the low pass filters $\F(\phi_1)$) we have $\sum_m|V_m(\theta)|^2=1$ and then $|\F(\phi_1)(\omega)|^2+\sum_n|W_n(\omega)|^2$ returns us to the 
Littlewood-Paley case. The second curvelet case is almost equivalent except that now the angular tiling depends on each scale which is equivalent to study 
$|\F(\phi_1)(\omega)|^2+\sum_n\left(|W_n(\omega)|^2\sum_m|V_m^n(\theta)|^2\right)$. One more time, if we consider a fixed scale $n$, by construction we have $\sum_m|V_m^n(\theta)|^2=1$ and we conclude 
with the same argument as the first case.
\end{proof}

\section{Experiments}\label{sec:expe}
\subsection{Transform examples}
\begin{figure}[!t]
\begin{center}
\includegraphics[scale=0.5]{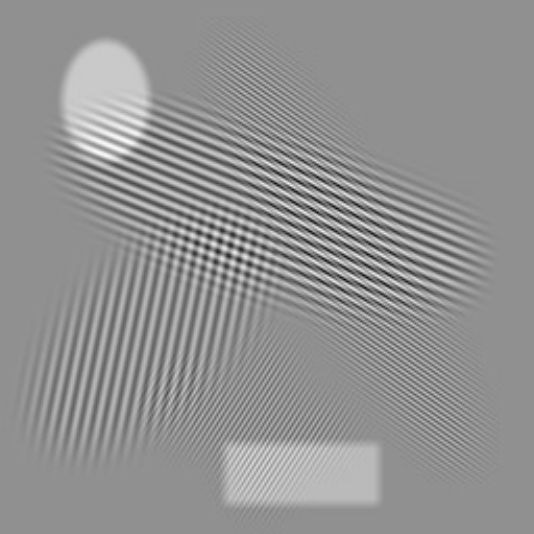}
\includegraphics[scale=0.5]{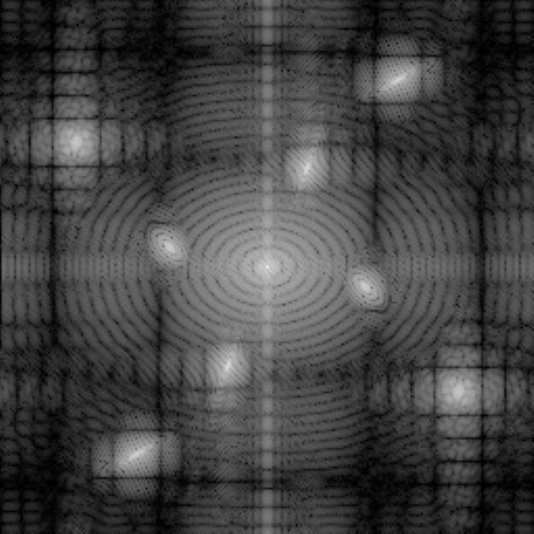}
\end{center}
\caption{Toy image and the logarithm magnitude of its Fourier spectrum.}
\label{fig:toy}
\end{figure}

\begin{figure}[!t]
\begin{center}
\includegraphics[width=0.33\textwidth]{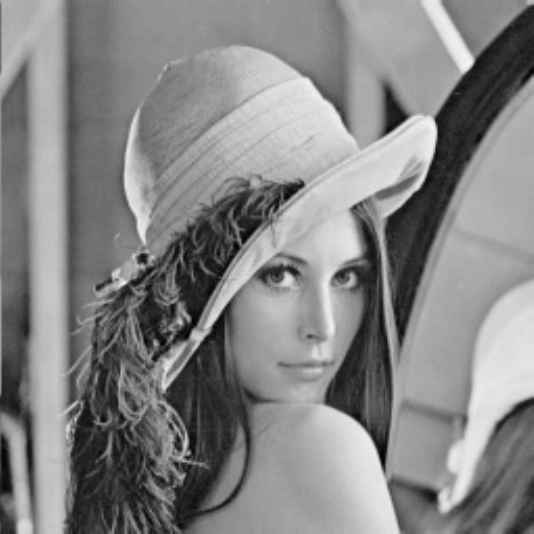}
\includegraphics[width=0.33\textwidth]{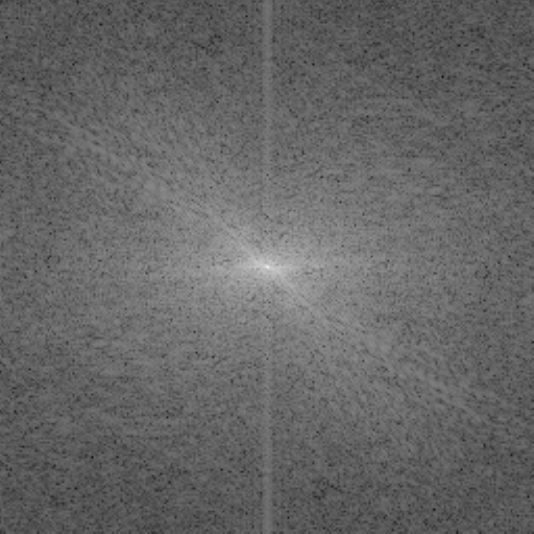}
\end{center}
\caption{Lena image and the logarithm magnitude of its Fourier spectrum.}
\label{fig:lena}
\end{figure}

\begin{figure}
\begin{tabular}{ccc}
\includegraphics[width=0.3\textwidth]{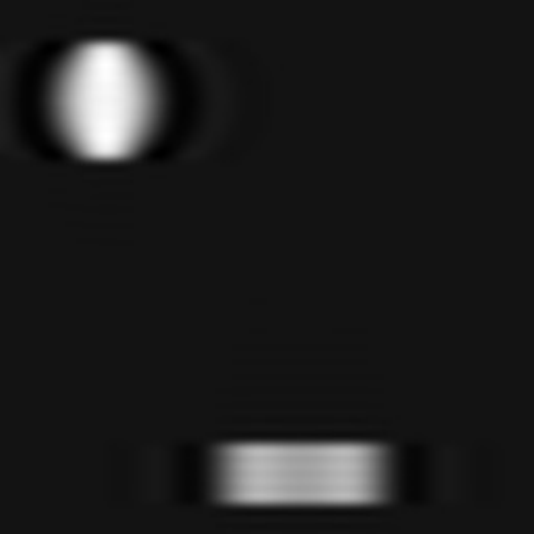} & \includegraphics[width=0.3\textwidth]{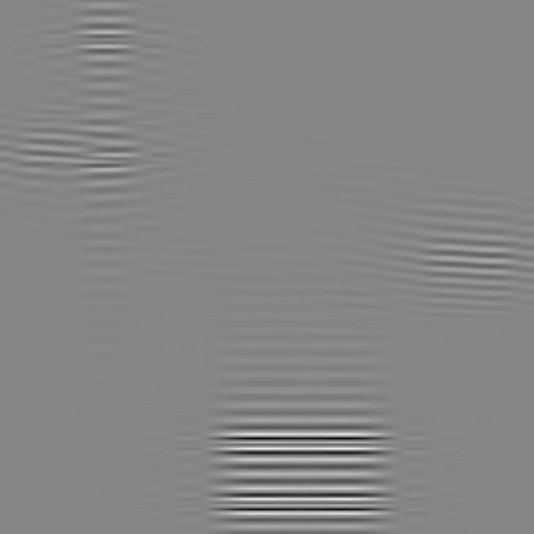} & \includegraphics[width=0.3\textwidth]{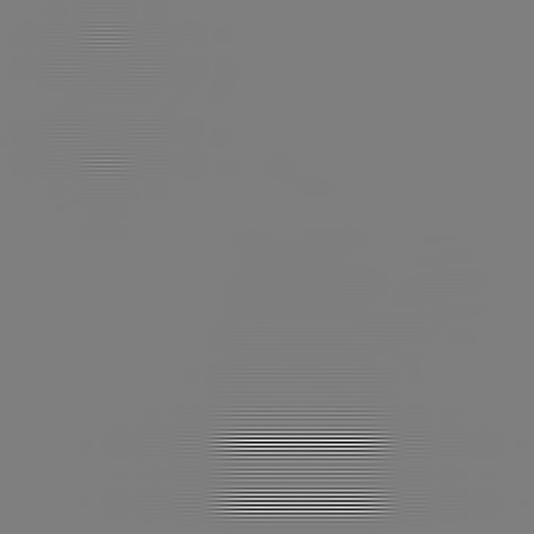}\\
\includegraphics[width=0.3\textwidth]{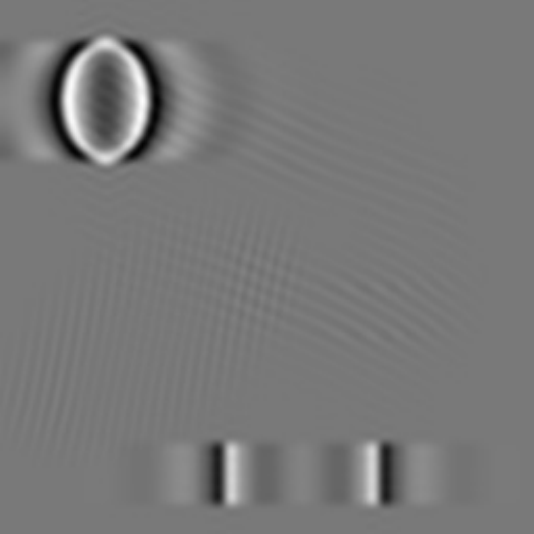} & \includegraphics[width=0.3\textwidth]{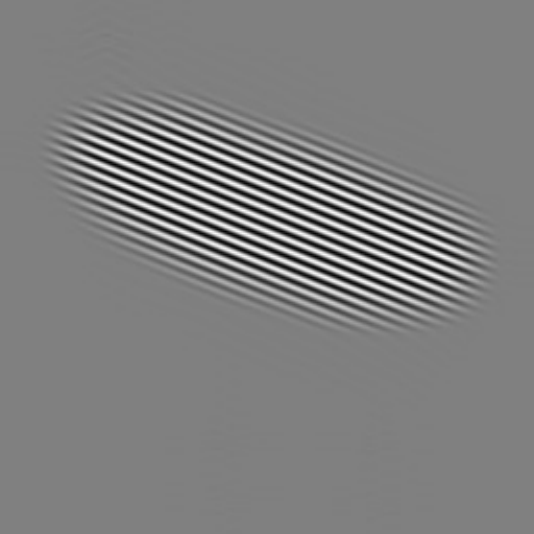} & \includegraphics[width=0.3\textwidth]{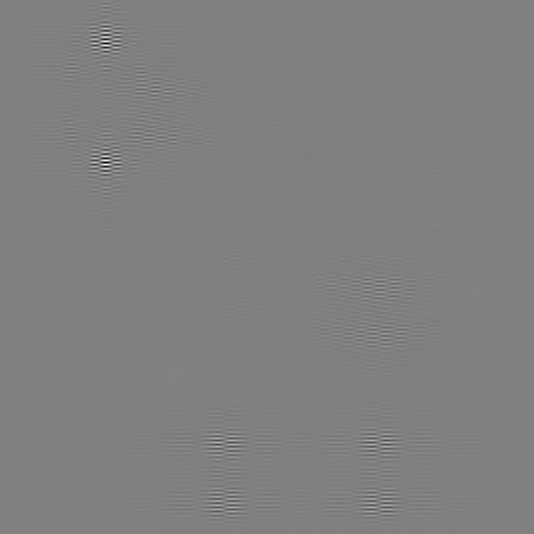}\\
\includegraphics[width=0.3\textwidth]{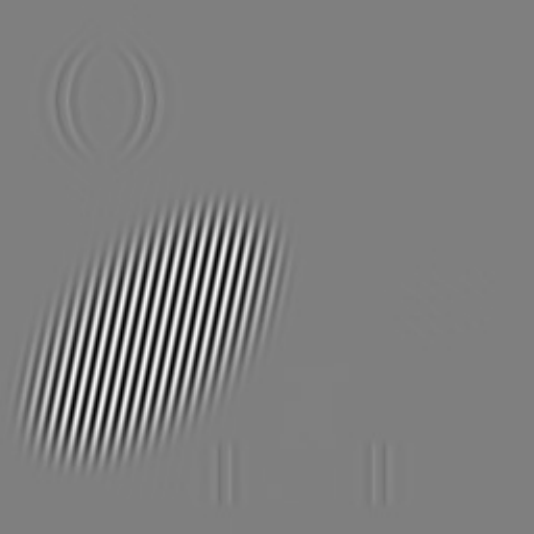} & \includegraphics[width=0.3\textwidth]{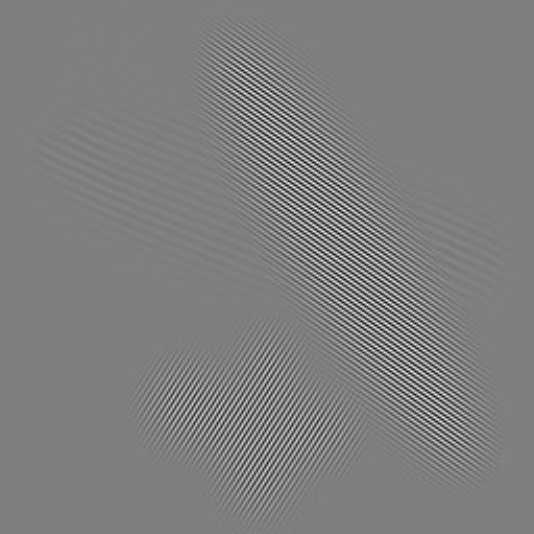} & \includegraphics[width=0.3\textwidth]{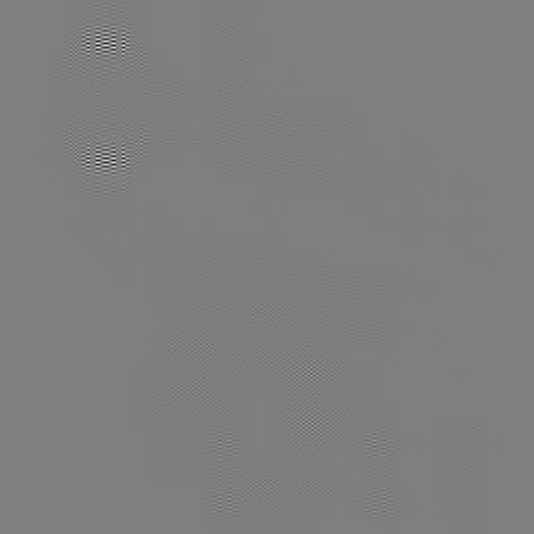}
\end{tabular}
\caption{2D Tensor EWT components of the toy image ($N_R=3, N_C=3$, the logarithm, the \textit{morpho} preprocessing and the lowest minima detection are used).}
\label{fig:tensortoy}
\end{figure}

\begin{figure}
\begin{center}
\begin{tabular}{cccc}
\includegraphics[width=0.22\textwidth]{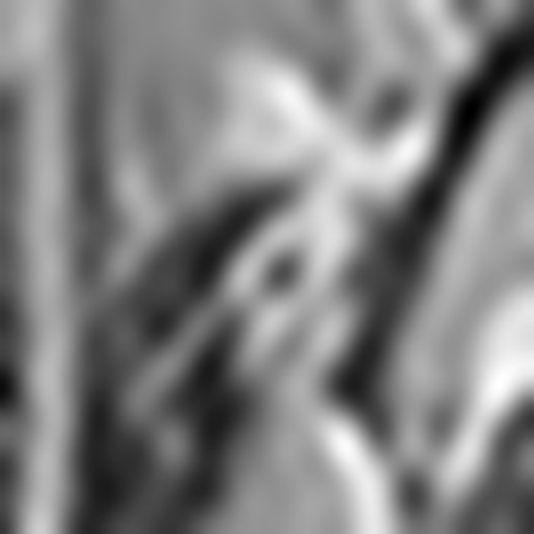} & \includegraphics[width=0.22\textwidth]{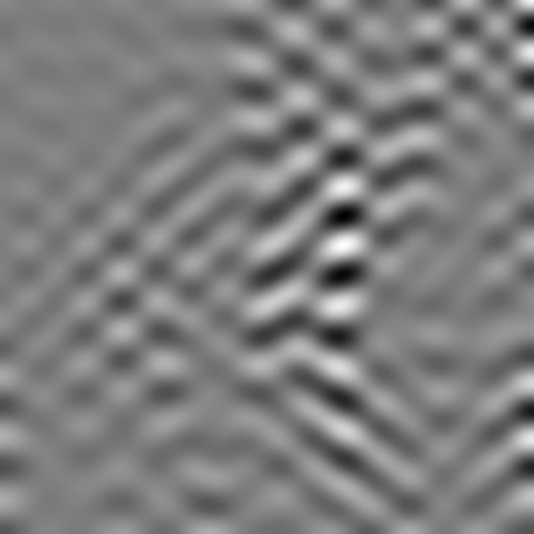} & \includegraphics[width=0.22\textwidth]{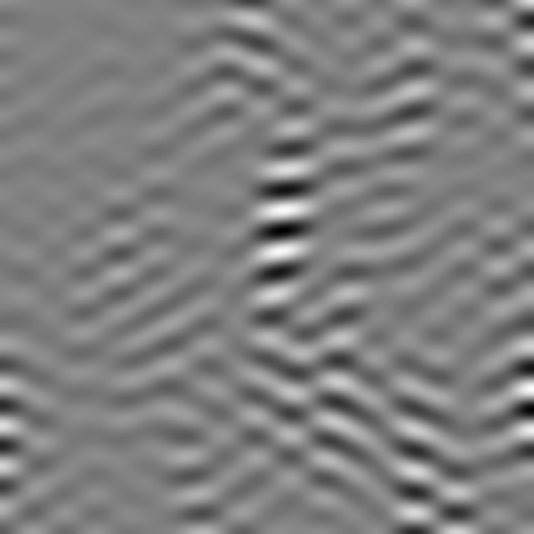} & \includegraphics[width=0.22\textwidth]{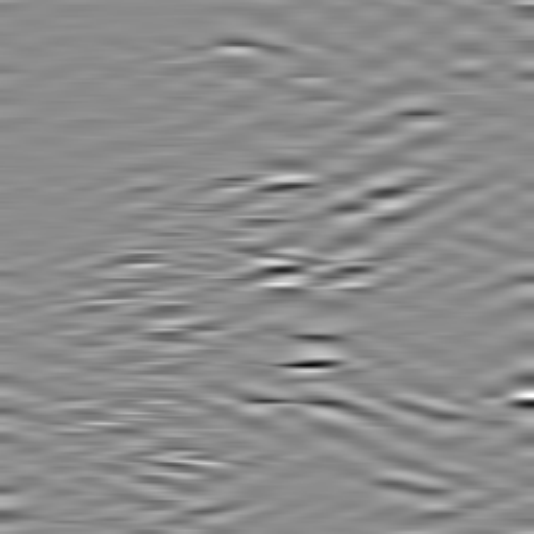}\\
\includegraphics[width=0.22\textwidth]{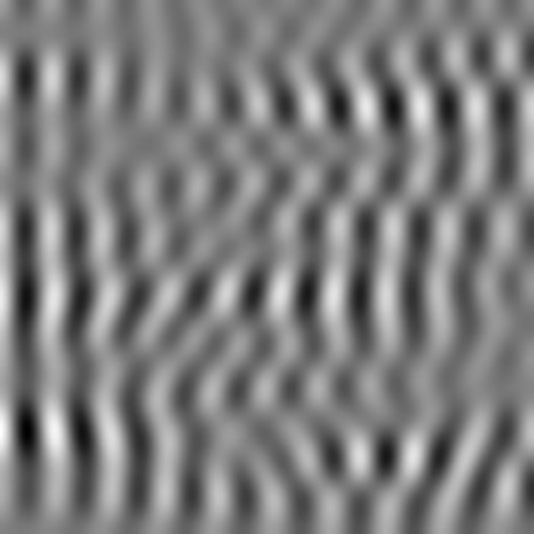} & \includegraphics[width=0.22\textwidth]{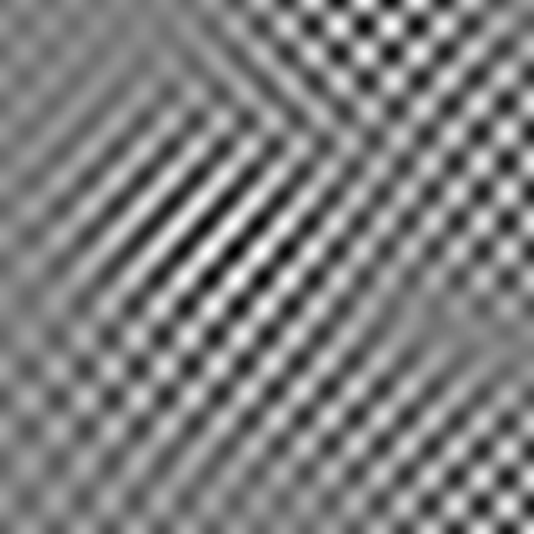} & \includegraphics[width=0.22\textwidth]{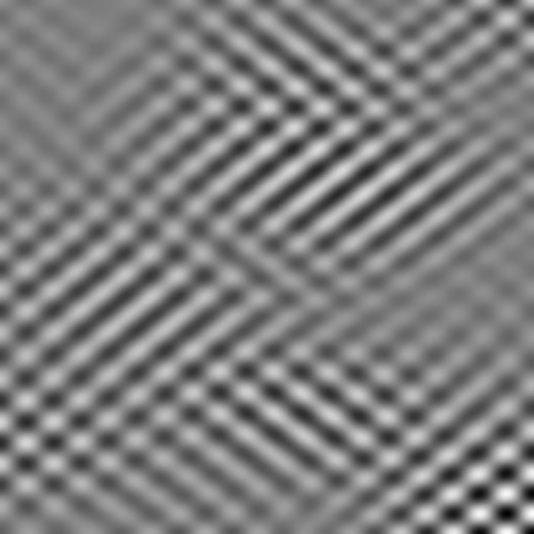} & \includegraphics[width=0.22\textwidth]{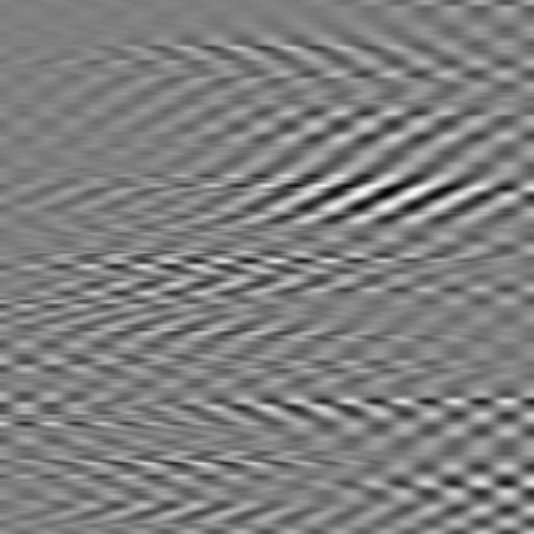}\\
\includegraphics[width=0.22\textwidth]{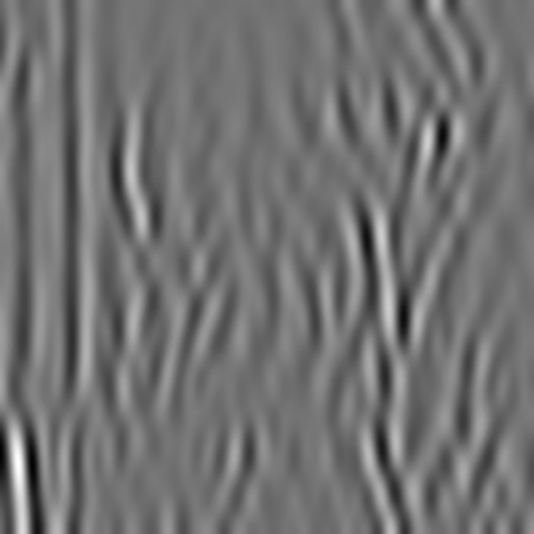} & \includegraphics[width=0.22\textwidth]{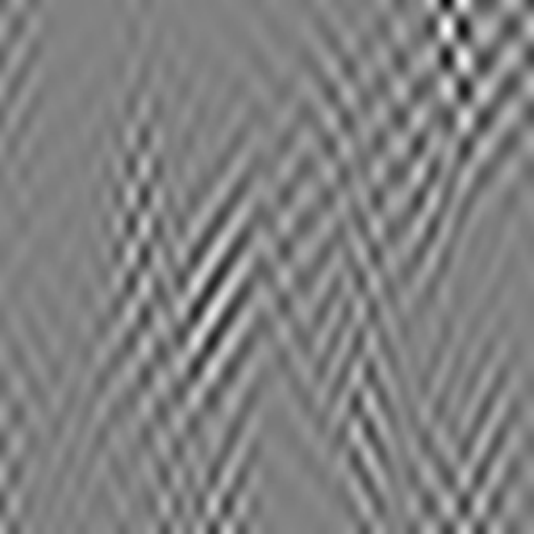} & \includegraphics[width=0.22\textwidth]{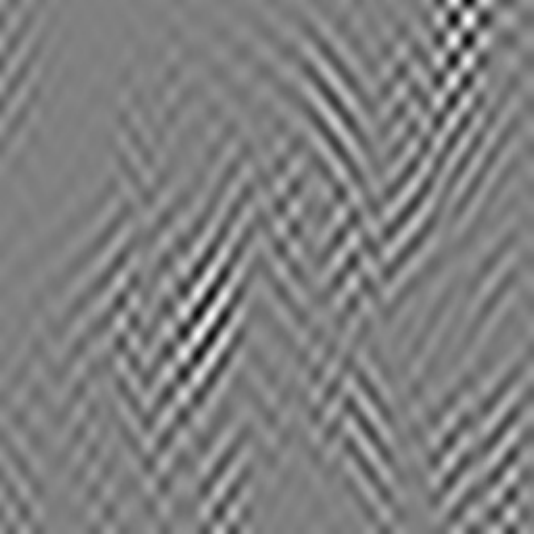} & \includegraphics[width=0.22\textwidth]{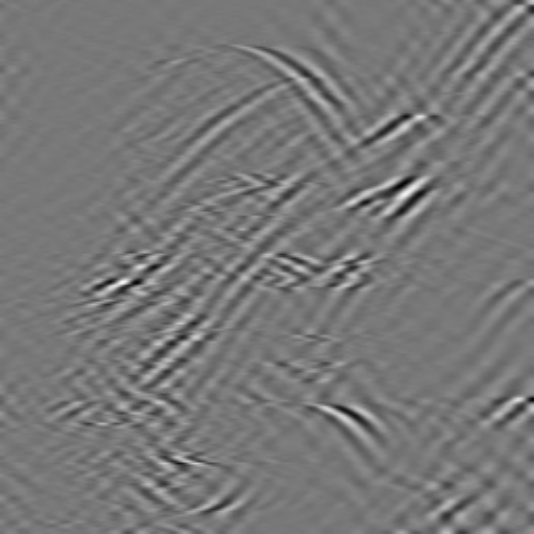} \\
\includegraphics[width=0.22\textwidth]{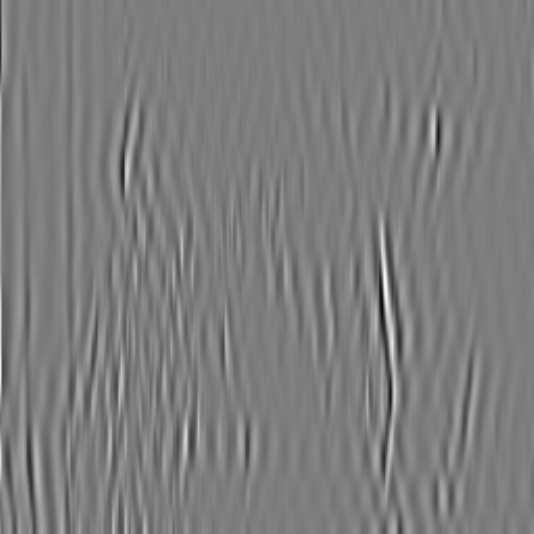} & \includegraphics[width=0.22\textwidth]{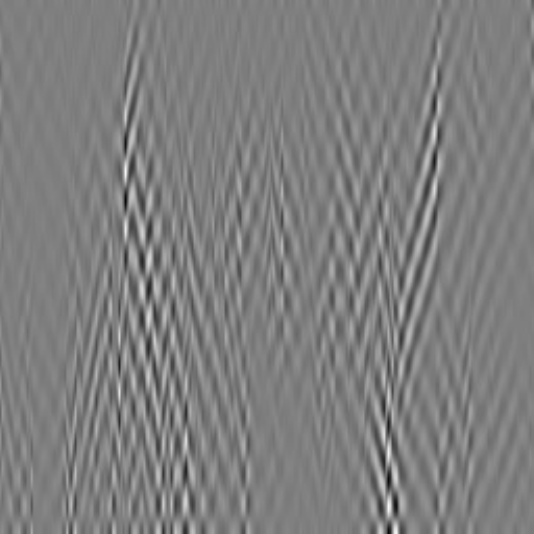} & \includegraphics[width=0.22\textwidth]{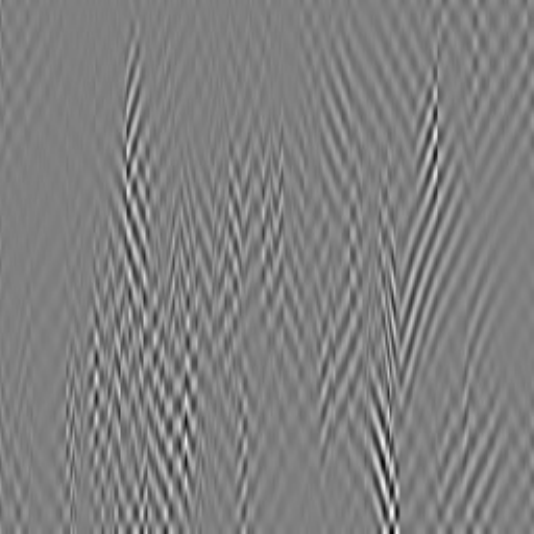} & \includegraphics[width=0.22\textwidth]{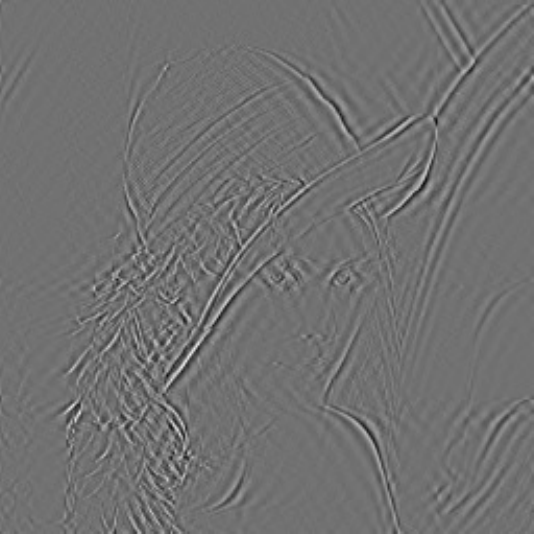} \\
\end{tabular}
\end{center}
\caption{2D Tensor EWT components of the Lena image ($N_R=4, N_C=4$, the logarithm, the \textit{morpho} preprocessing and the lowest minima detection are used).}
\label{fig:tensorlena}
\end{figure}

\begin{figure}
\begin{center}
\begin{tabular}{cc}
\includegraphics[width=0.33\textwidth]{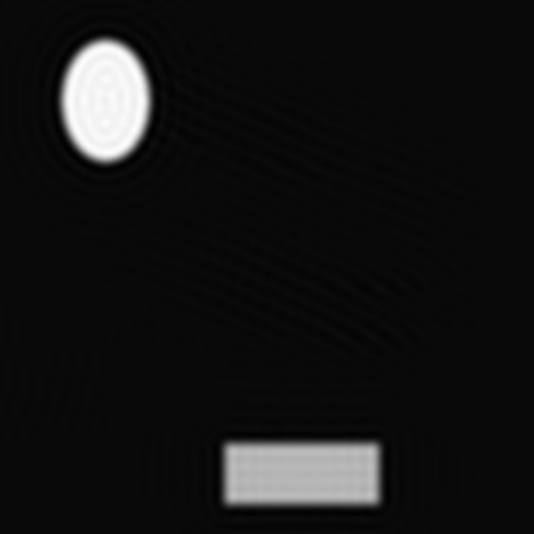} & \includegraphics[width=0.33\textwidth]{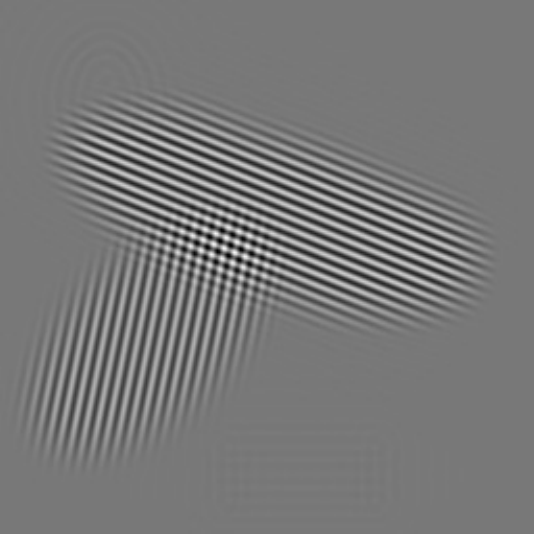} \\
\includegraphics[width=0.33\textwidth]{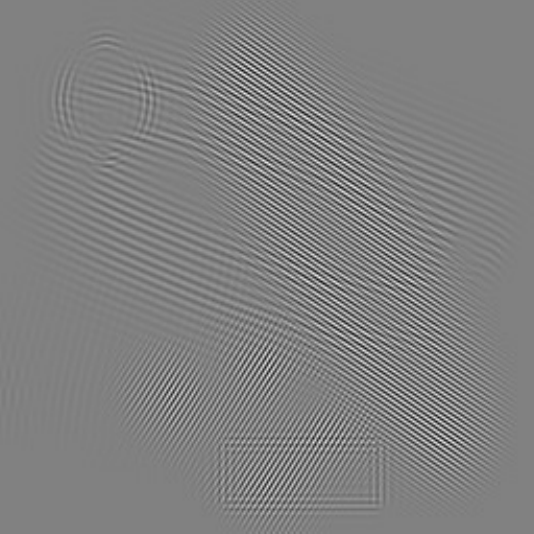} & \includegraphics[width=0.33\textwidth]{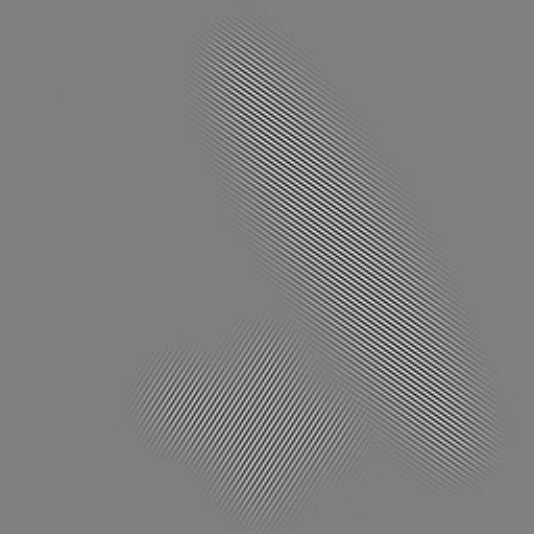} 
\end{tabular}
\end{center}
\caption{2D Littlewood-Paley EWT components of the toy image ($N=4$, the logarithm, the \textit{morpho} preprocessing and the lowest minima detection are used).}
\label{fig:lptoy}
\end{figure}

\begin{figure}
\begin{center}
\begin{tabular}{ccc}
\includegraphics[width=0.3\textwidth]{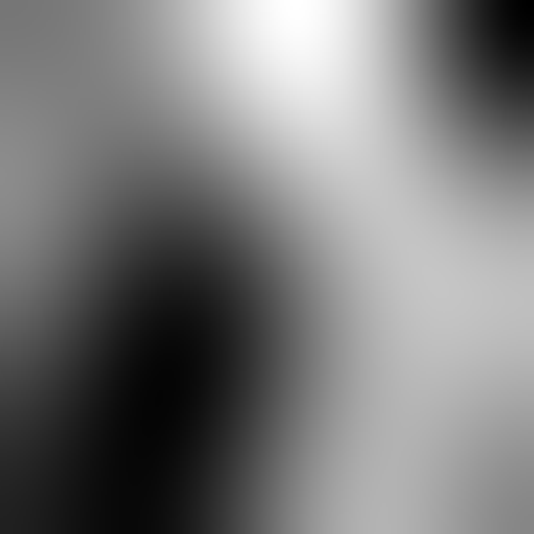} & \includegraphics[width=0.3\textwidth]{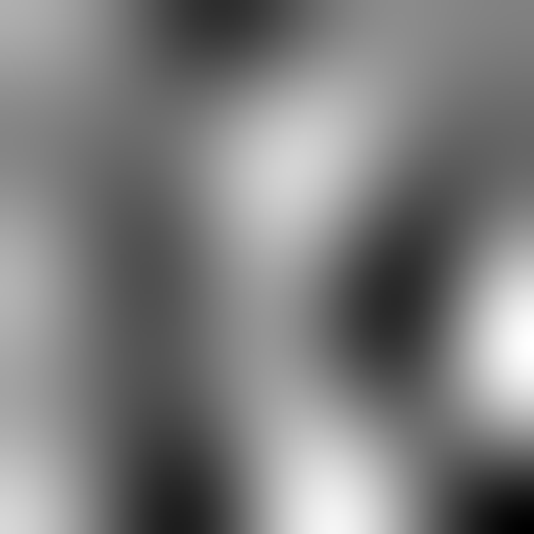} & \includegraphics[width=0.3\textwidth]{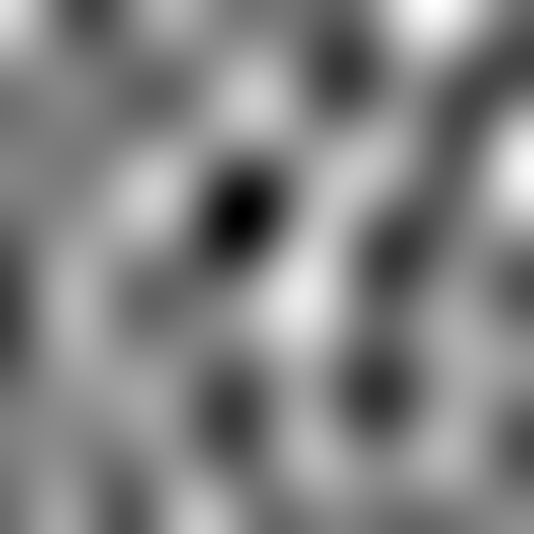}\\
\includegraphics[width=0.3\textwidth]{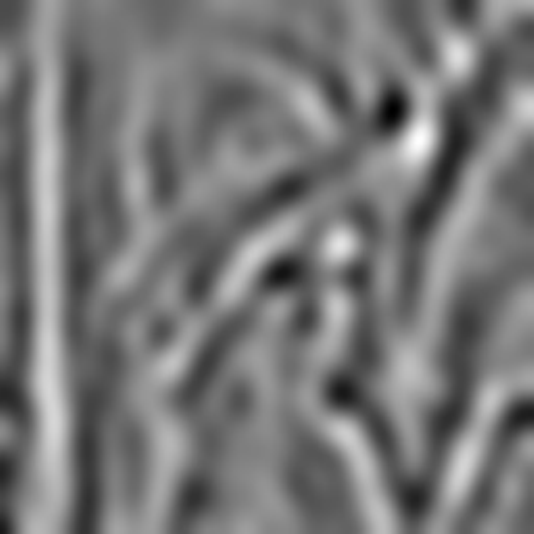} & \includegraphics[width=0.3\textwidth]{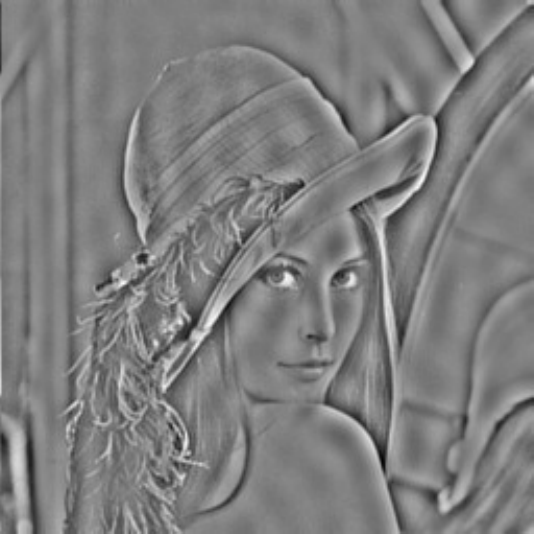} & \includegraphics[width=0.3\textwidth]{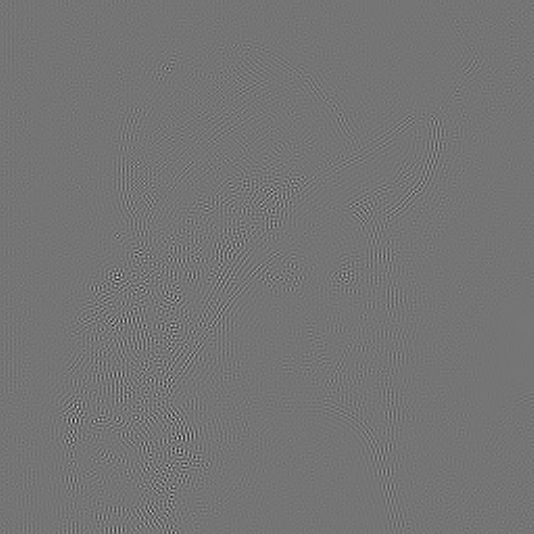}
\end{tabular}
\end{center}
\caption{2D Littlewood-Paley EWT components of the Lena image ($N=6$, the logarithm, the \textit{morpho} preprocessing and the lowest minima detection are used).}
\label{fig:lplena}
\end{figure}

\begin{figure}
\begin{center}
\begin{tabular}{cc}
\includegraphics[width=0.28\textwidth]{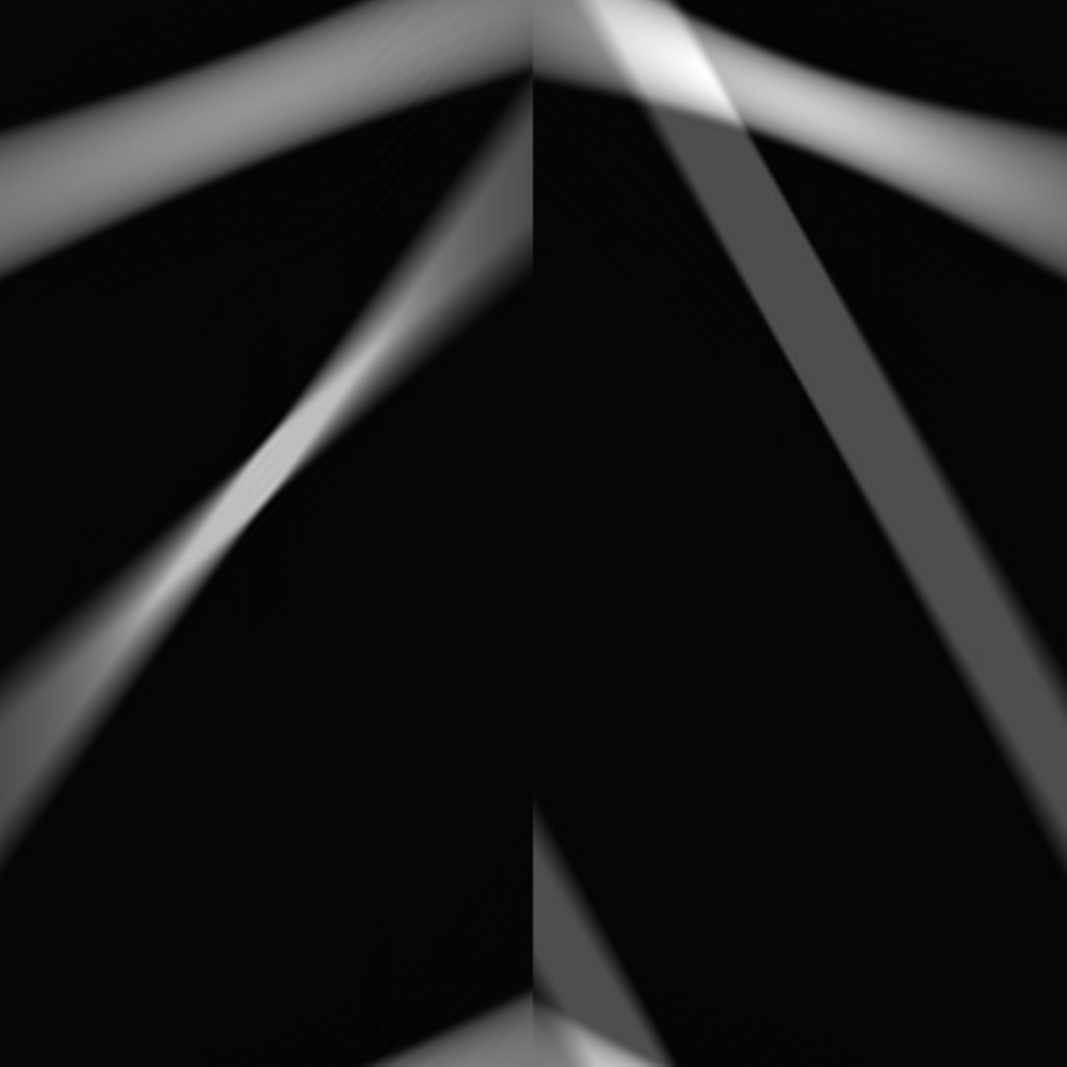} & \includegraphics[width=0.28\textwidth]{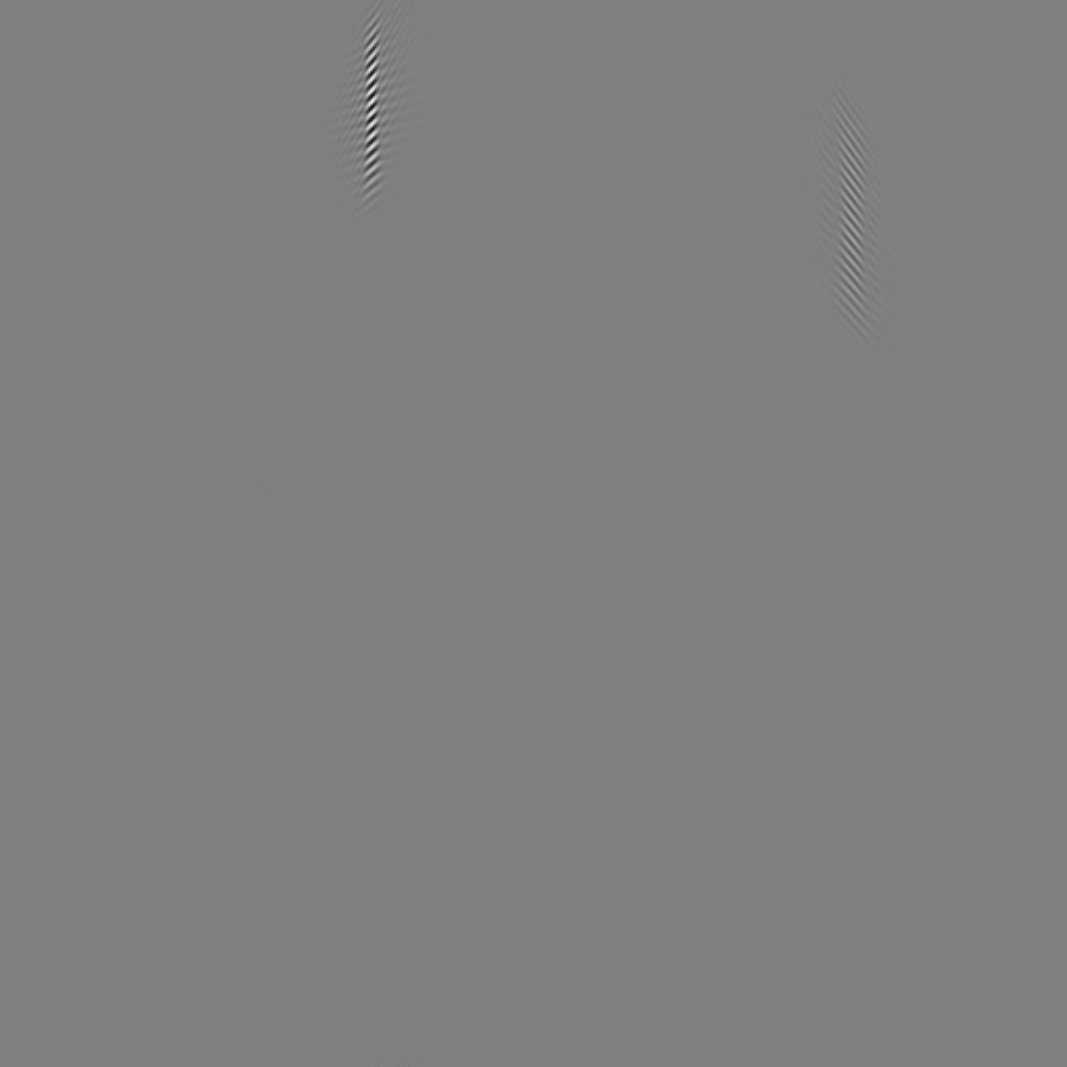} \\
\includegraphics[width=0.28\textwidth]{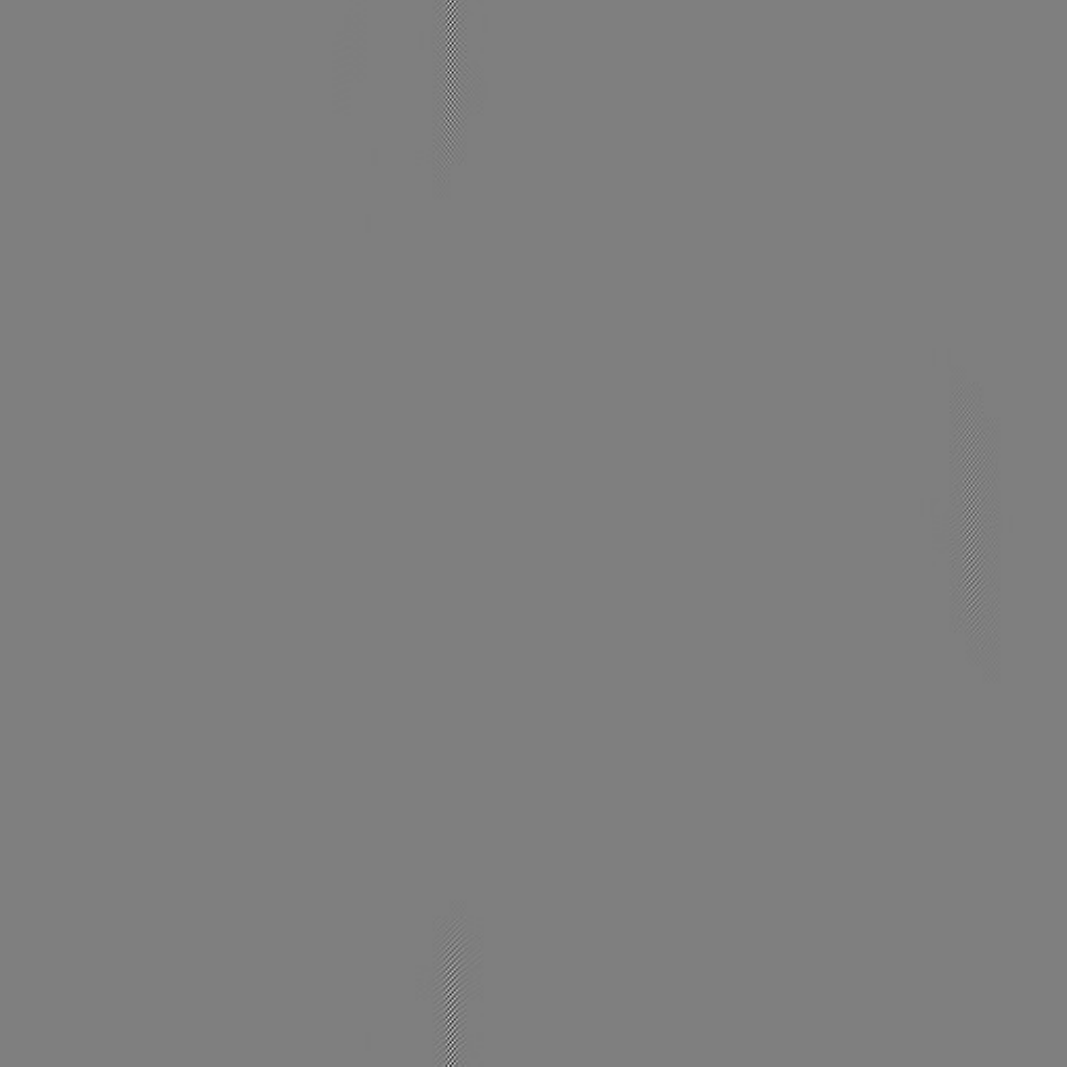} & \includegraphics[width=0.28\textwidth]{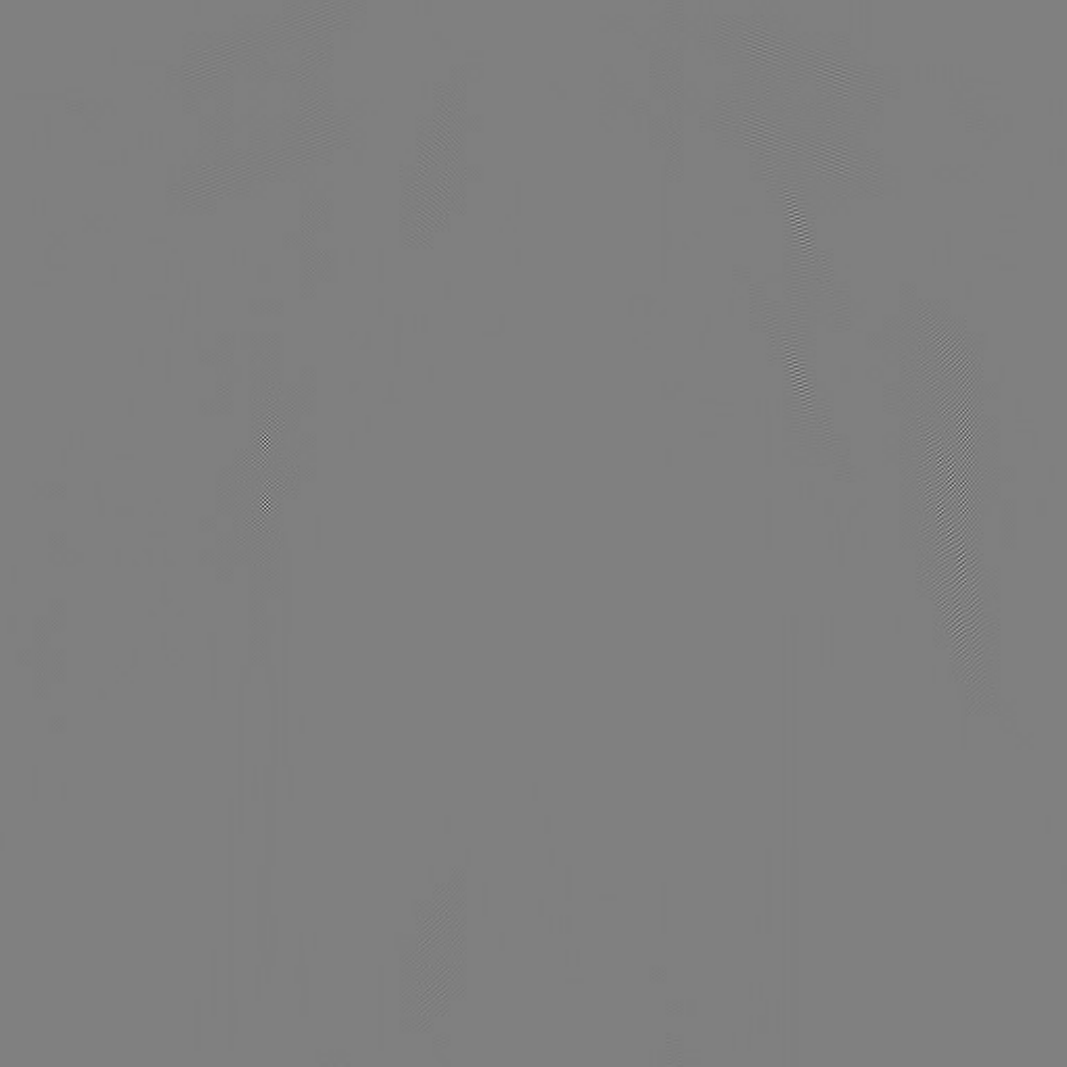} 
\end{tabular}
\end{center}
\caption{Empirical ridgelet components of the toy image ($N=4$).}
\label{fig:ridgetoy}
\end{figure}

\begin{figure}
\begin{center}
\begin{tabular}{ccc}
\includegraphics[width=0.28\textwidth]{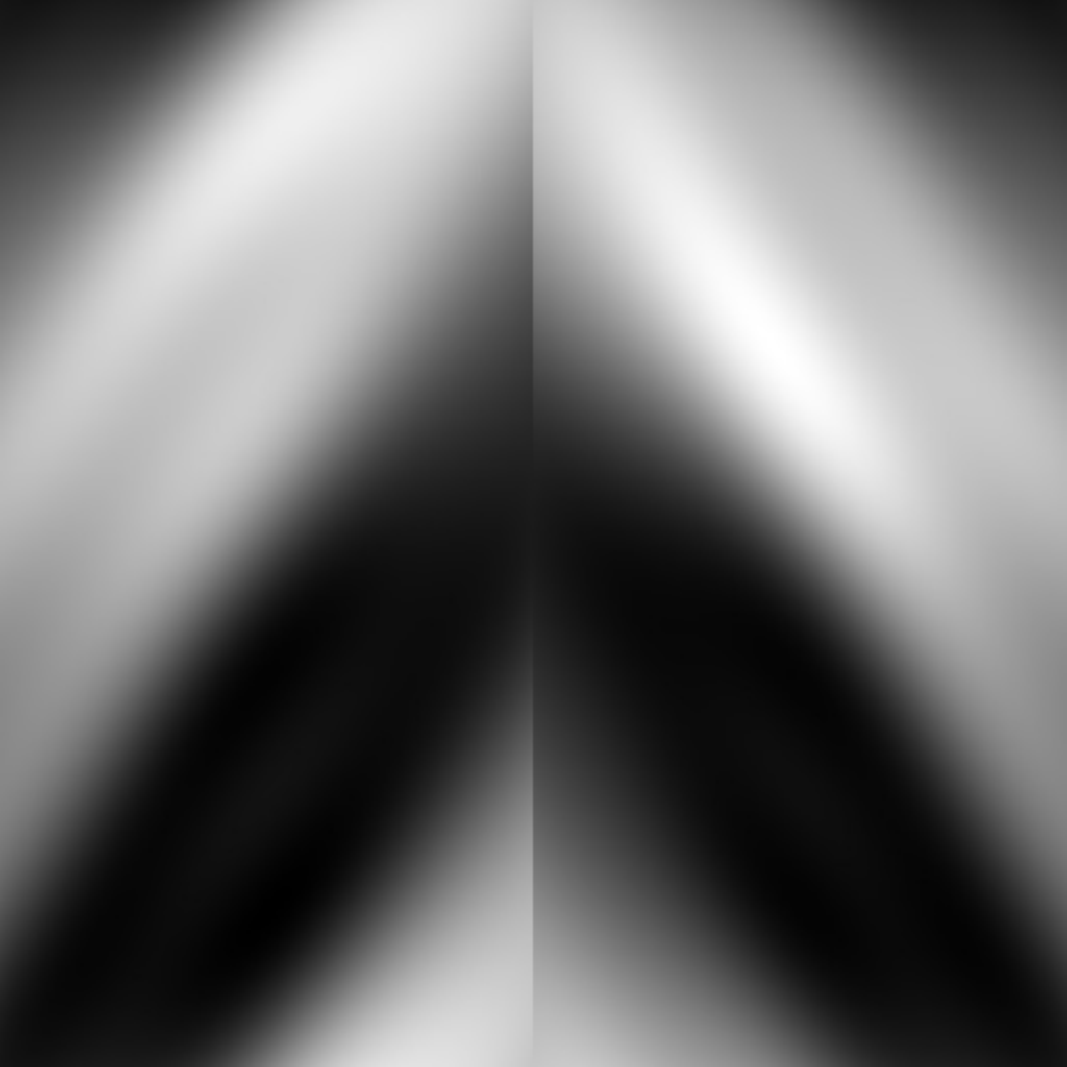} & \includegraphics[width=0.28\textwidth]{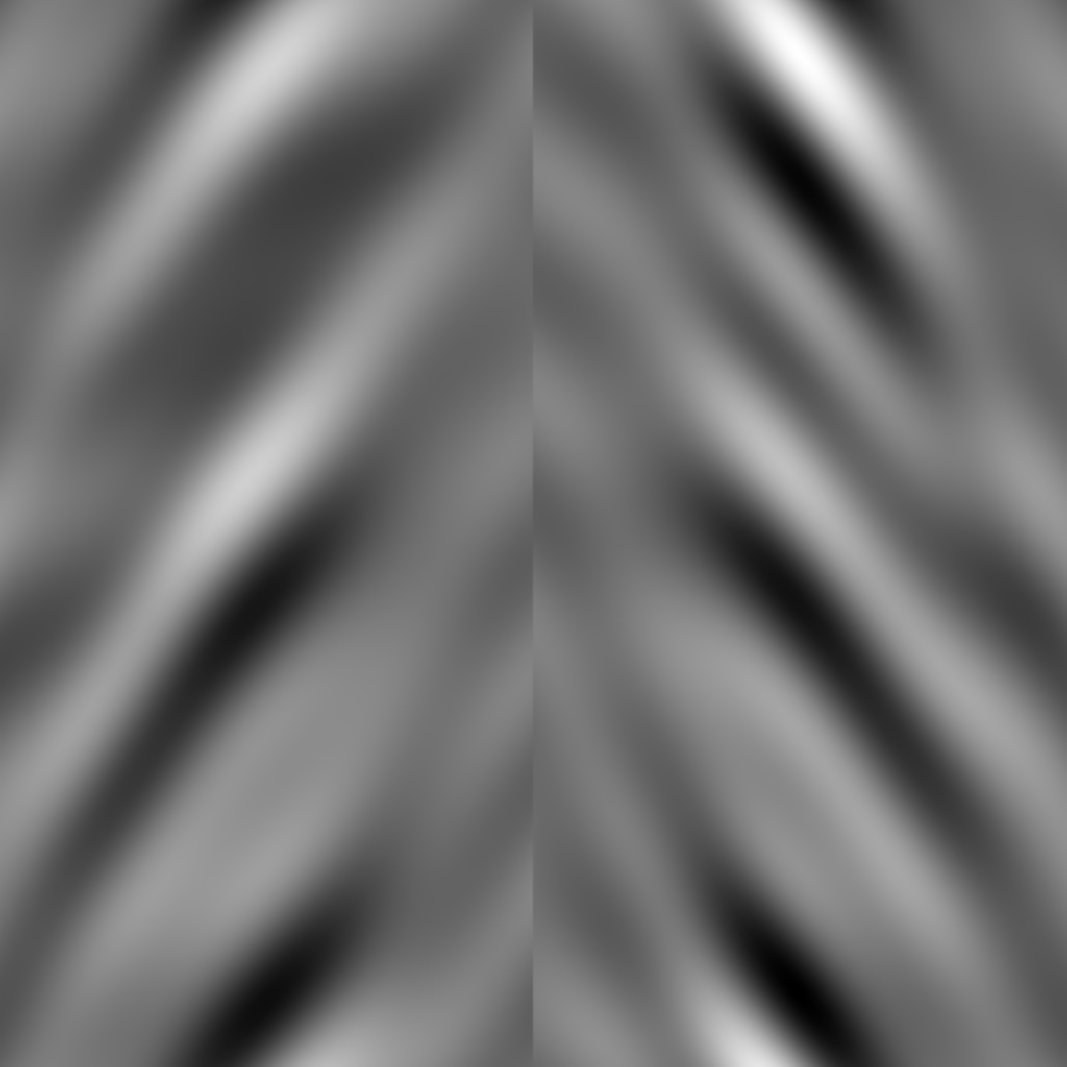} & \includegraphics[width=0.28\textwidth]{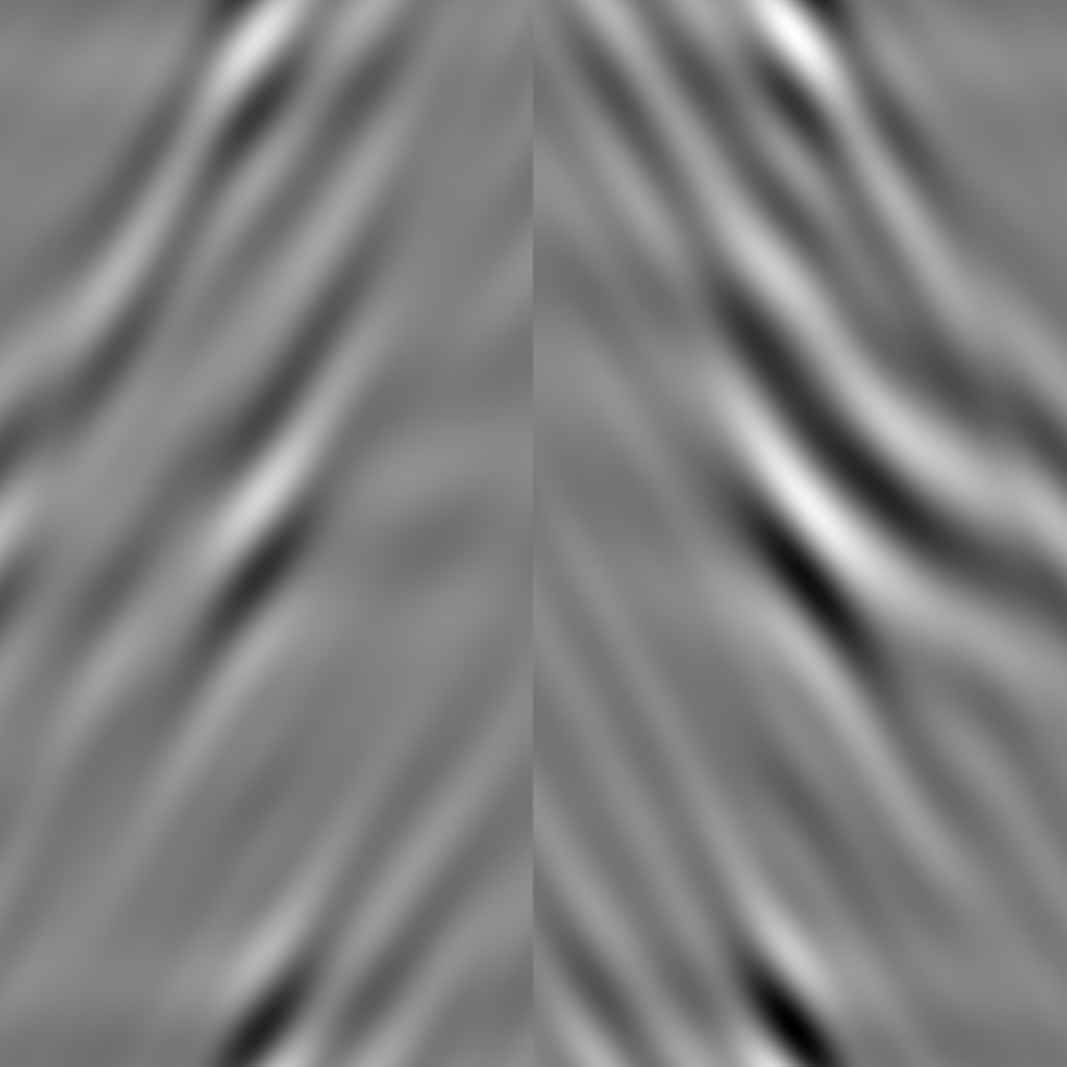}\\
\includegraphics[width=0.28\textwidth]{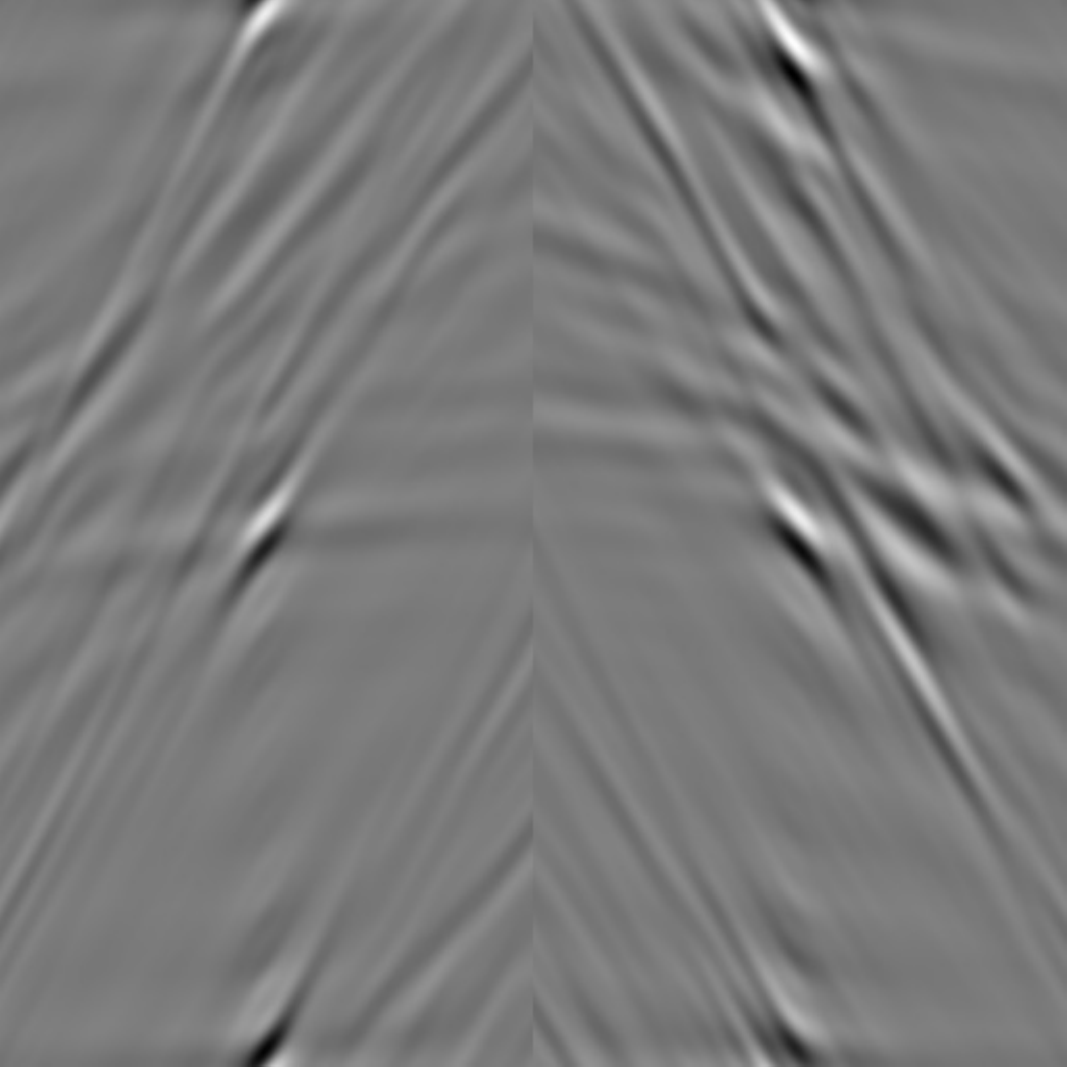} & \includegraphics[width=0.28\textwidth]{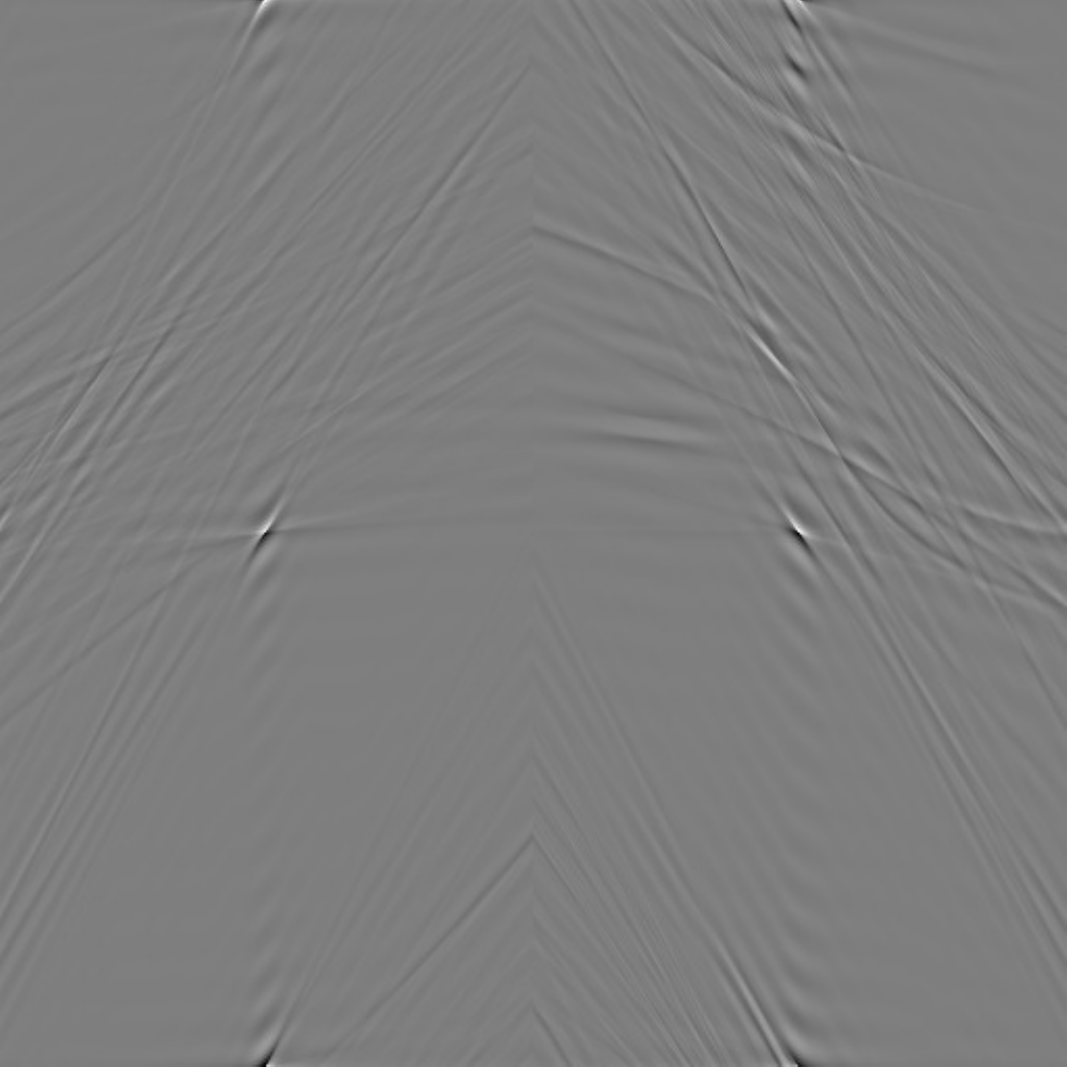} & \includegraphics[width=0.28\textwidth]{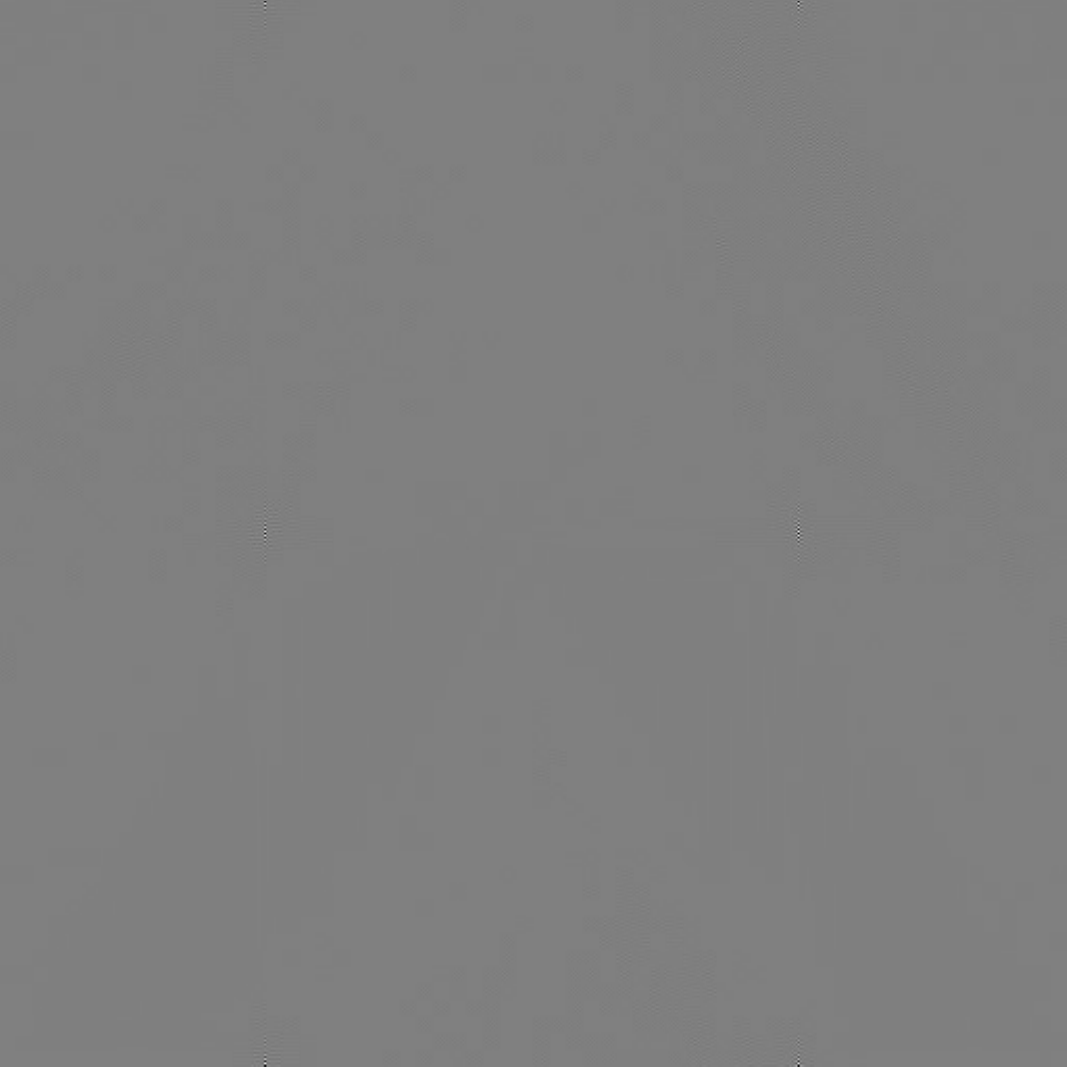}
\end{tabular}
\end{center}
\caption{Empirical ridgelet components of the Lena image ($N=6$).}
\label{fig:ridgelena}
\end{figure}

\begin{figure}[!t]
\begin{center}
\begin{tabular}{cc}
\includegraphics[width=0.27\textwidth]{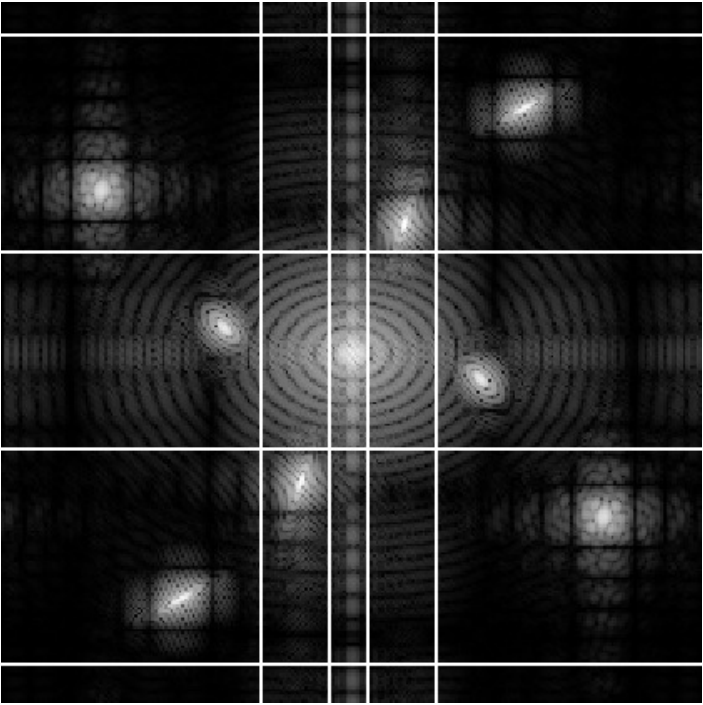} & \includegraphics[width=0.27\textwidth]{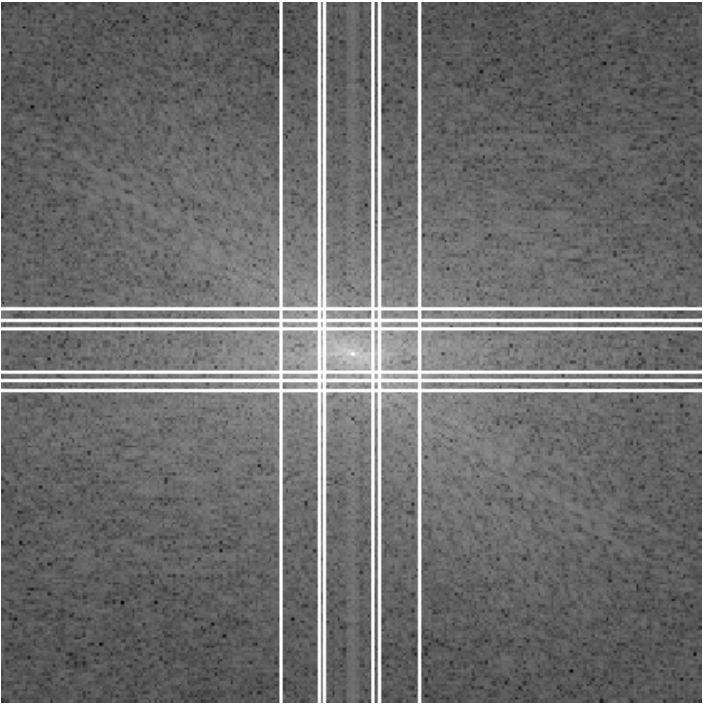} \\
\includegraphics[width=0.27\textwidth]{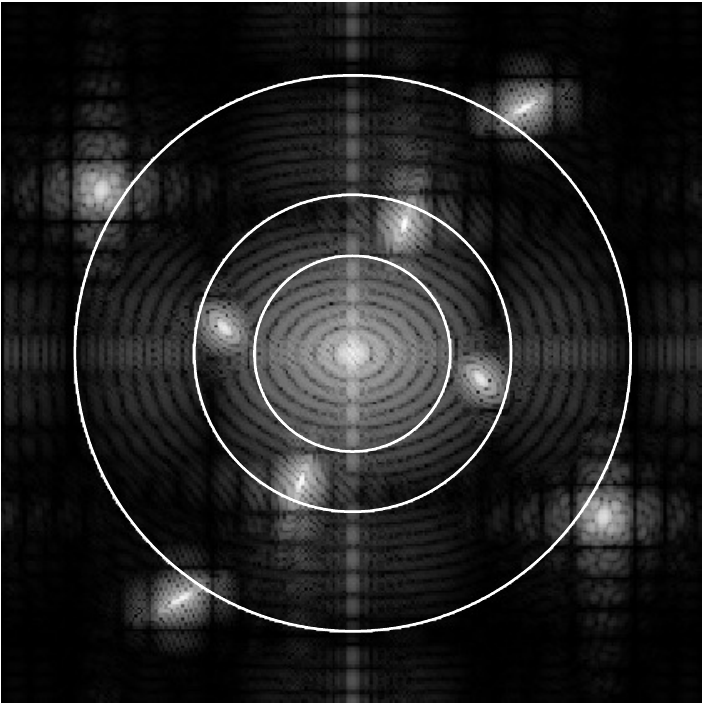} & \includegraphics[width=0.27\textwidth]{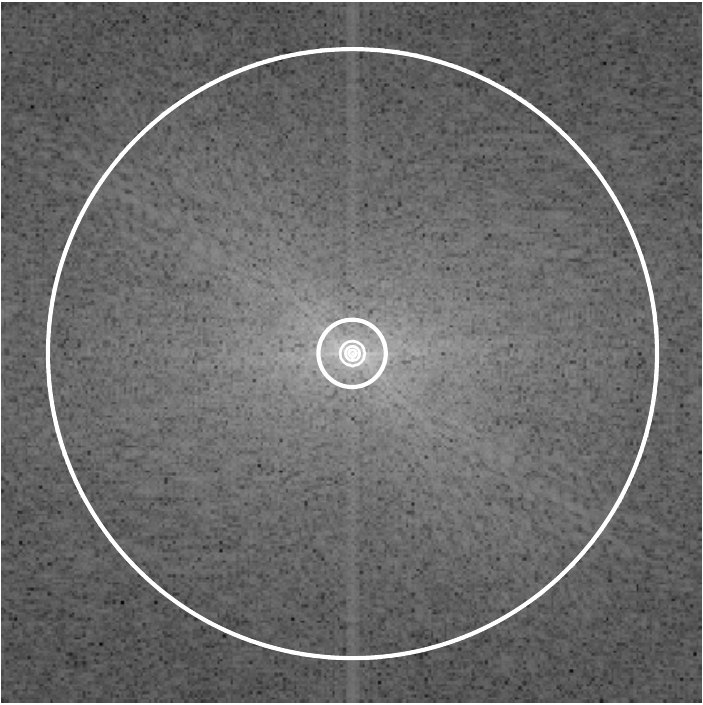}
\end{tabular}
\end{center}
\caption{Detected Fourier boundaries for each test image (toy image on left and Lena on right). Tensor EWT approach on top and Littlewood-Paley/Ridgelet approaches on bottom.}
\label{fig:Bounds}
\end{figure}

\begin{figure}
\begin{center}
\begin{tabular}{cccc}
\includegraphics[width=0.22\textwidth]{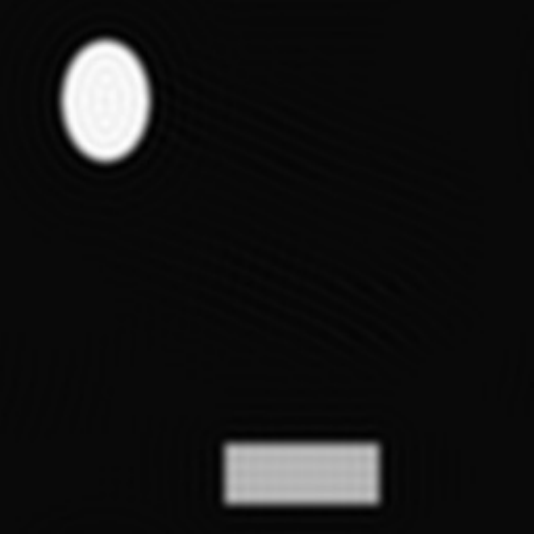} &  &  & \\
\includegraphics[width=0.22\textwidth]{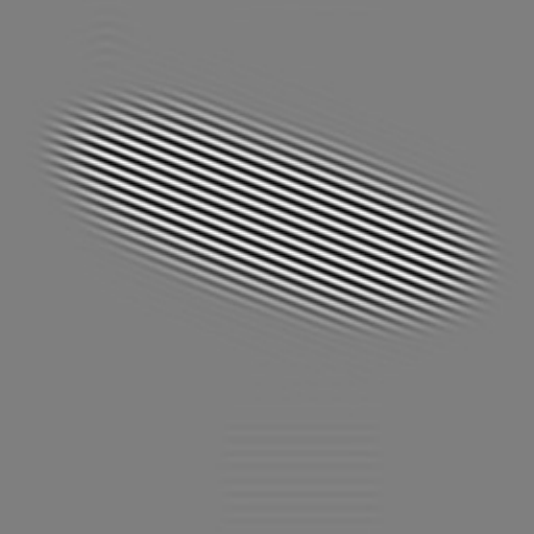} & \includegraphics[width=0.22\textwidth]{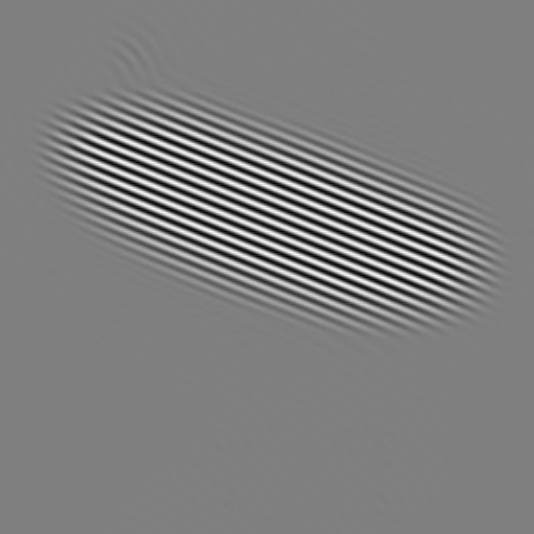} & \includegraphics[width=0.22\textwidth]{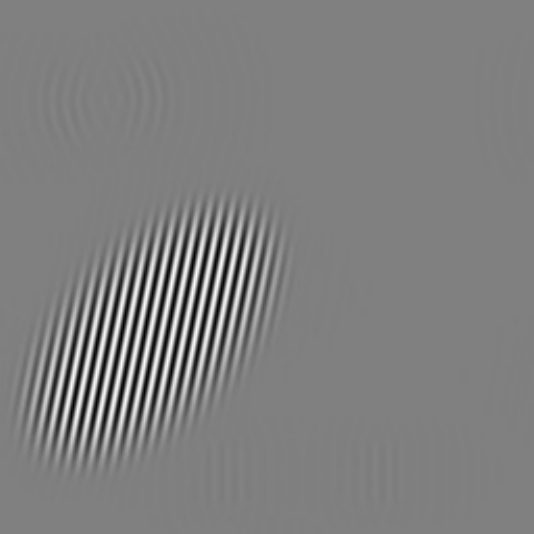} & \includegraphics[width=0.22\textwidth]{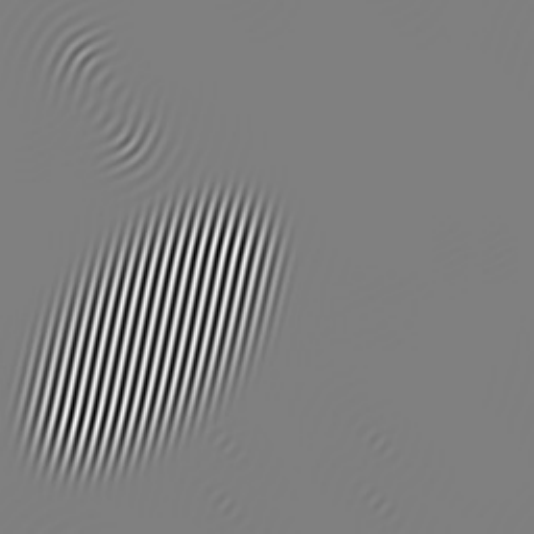}\\
\includegraphics[width=0.22\textwidth]{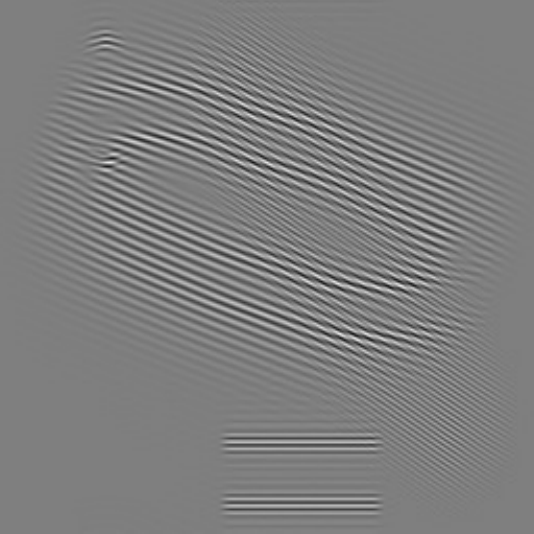} & \includegraphics[width=0.22\textwidth]{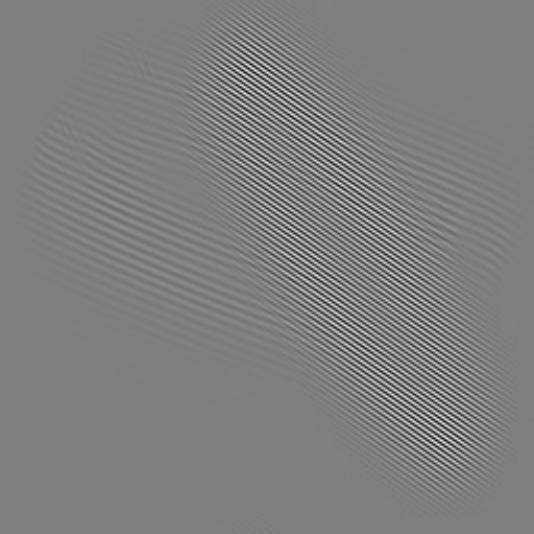} & \includegraphics[width=0.22\textwidth]{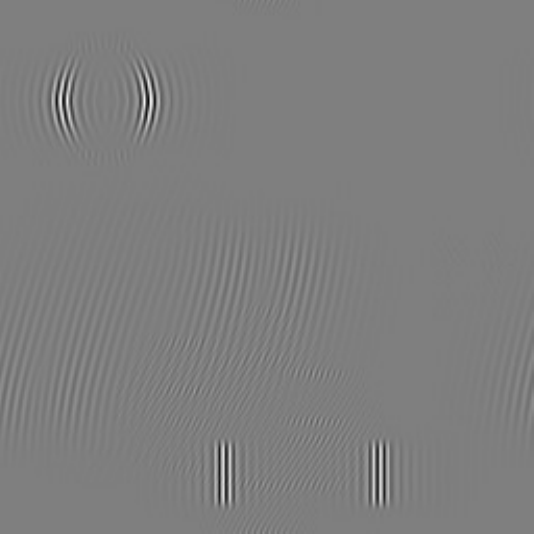} & \includegraphics[width=0.22\textwidth]{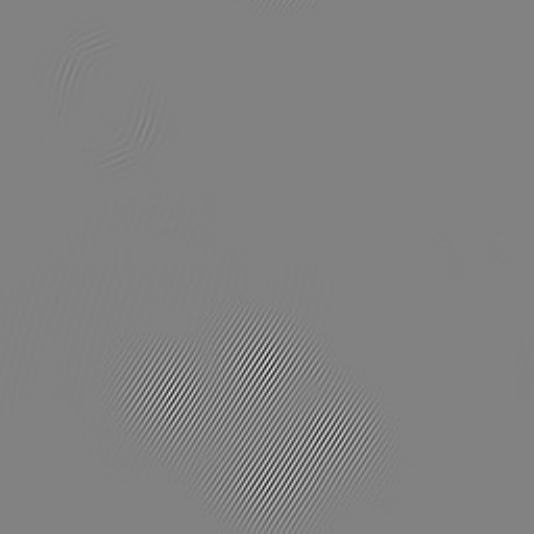} \\
\includegraphics[width=0.22\textwidth]{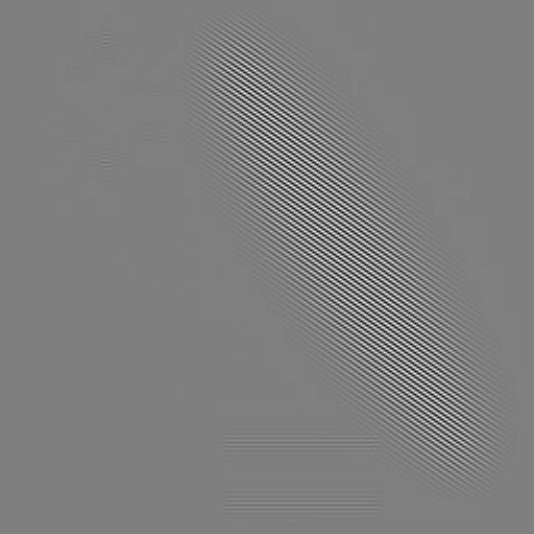} & \includegraphics[width=0.22\textwidth]{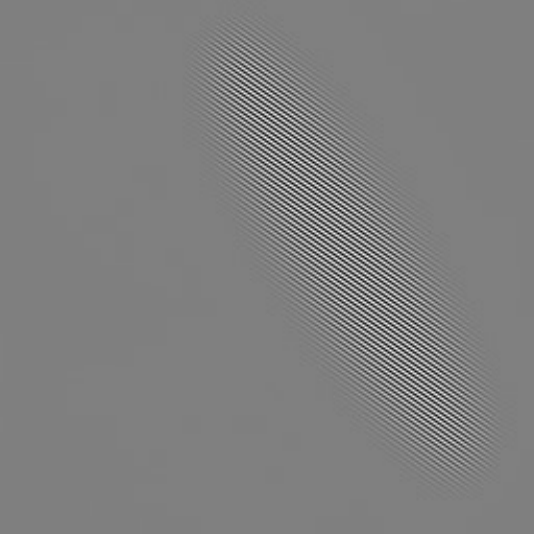} & \includegraphics[width=0.22\textwidth]{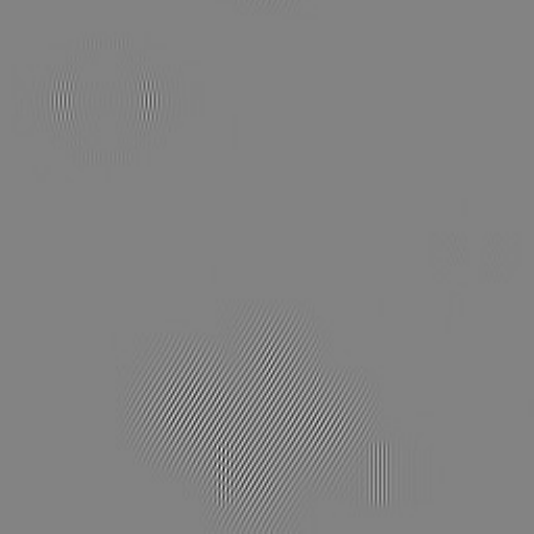} & \includegraphics[width=0.22\textwidth]{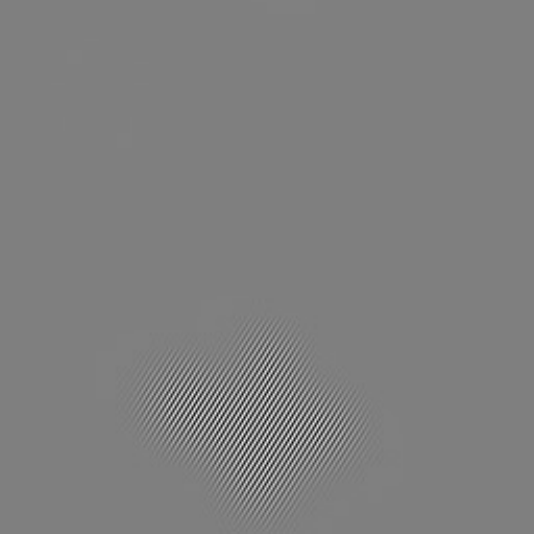} \\
\end{tabular}
\end{center}
\caption{2D Curvelet EWT-I components of the toy image ($N_s=4, N_{\theta}=4$, the logarithm, the \textit{morpho} preprocessing and the lowest minima detection are used to detect the scales and middle between local 
maxima with \textit{tophat} are used to detect the angles).}
\label{fig:cur1toy}
\end{figure}

\begin{figure}
\begin{center}
\begin{tabular}{cccc}
\includegraphics[width=0.22\textwidth]{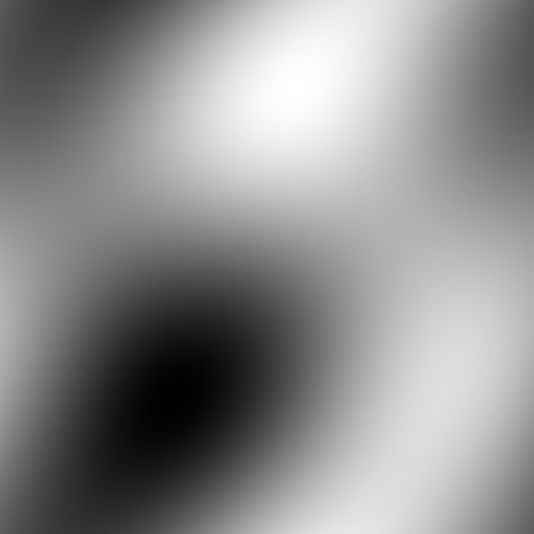} &  &  & \\
\includegraphics[width=0.22\textwidth]{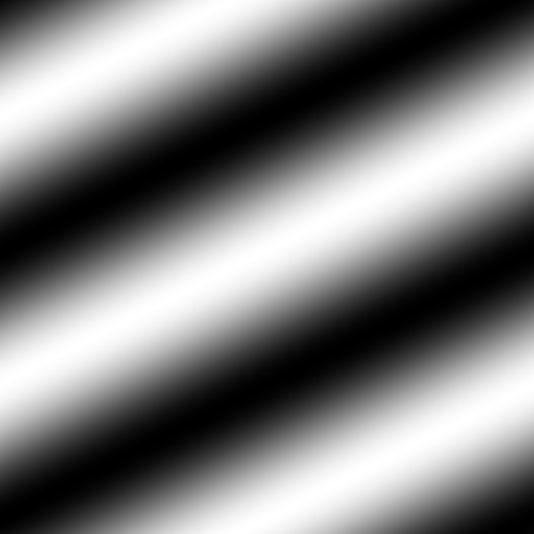} & \includegraphics[width=0.22\textwidth]{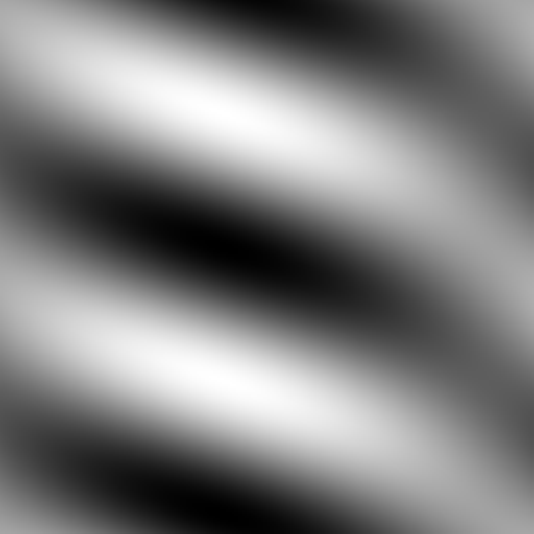} & \includegraphics[width=0.22\textwidth]{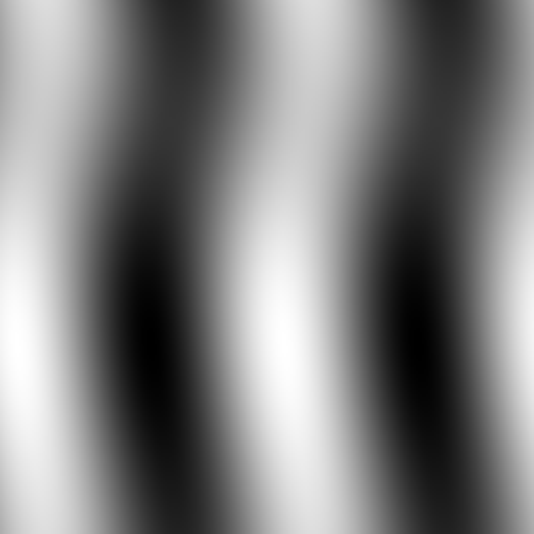} & \includegraphics[width=0.22\textwidth]{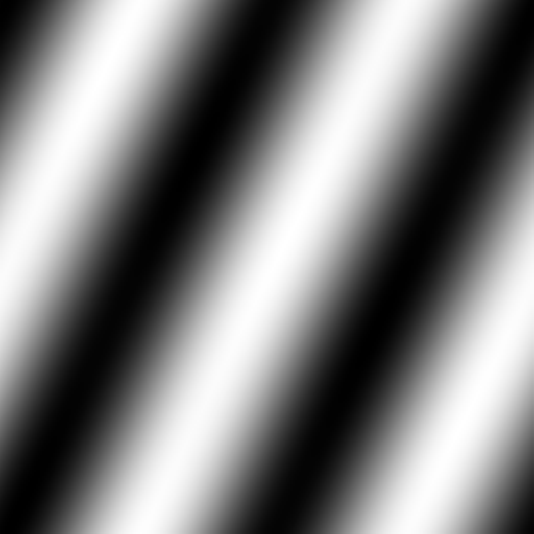}\\
\includegraphics[width=0.22\textwidth]{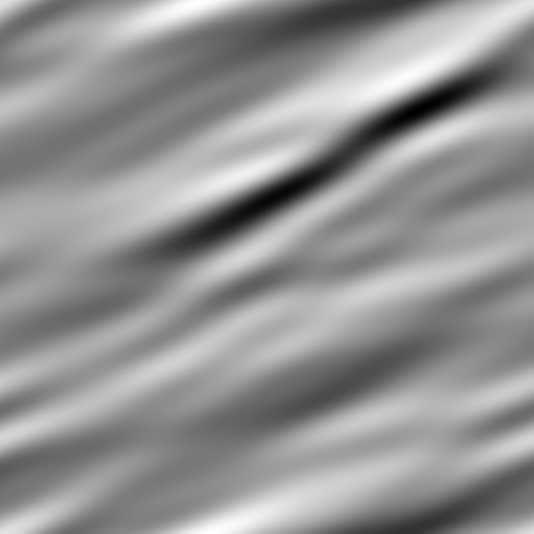} & \includegraphics[width=0.22\textwidth]{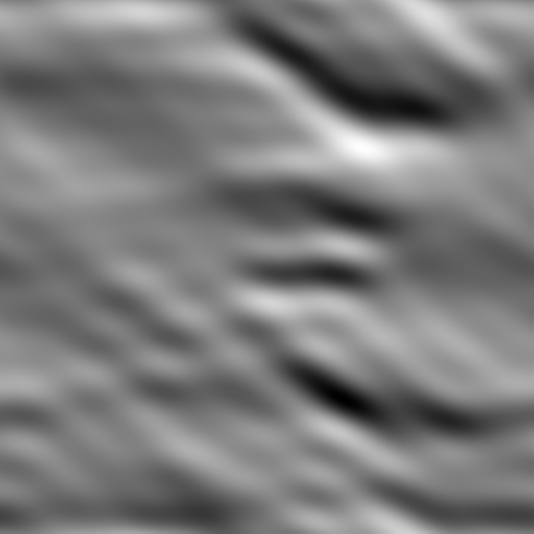} & \includegraphics[width=0.22\textwidth]{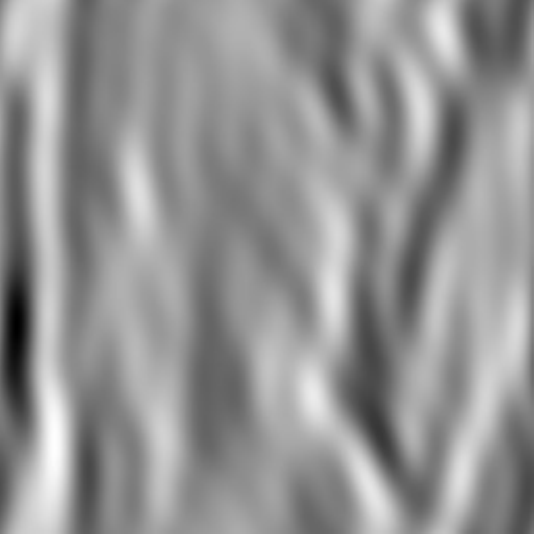} & \includegraphics[width=0.22\textwidth]{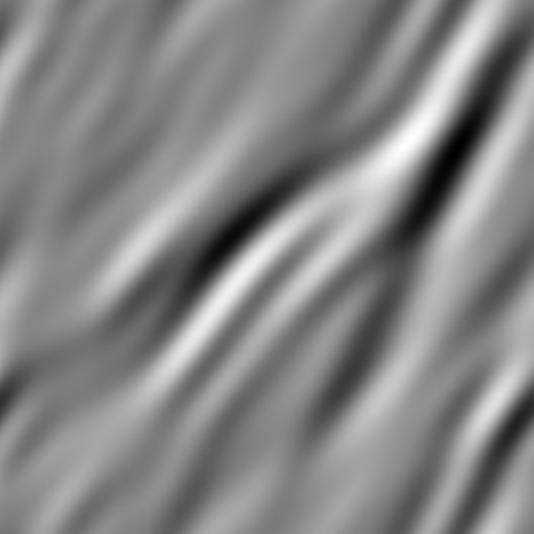} \\
\includegraphics[width=0.22\textwidth]{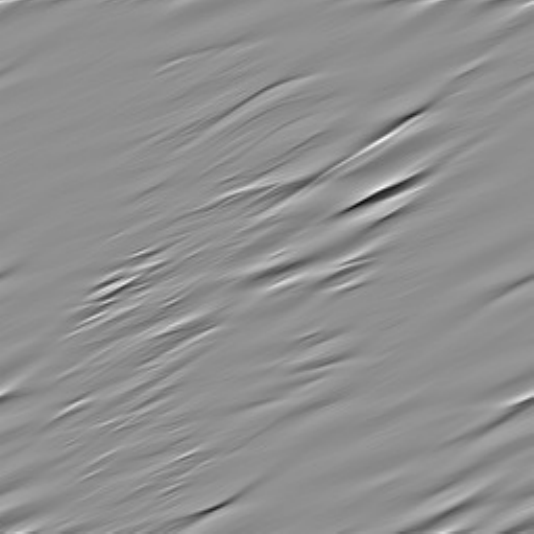} & \includegraphics[width=0.22\textwidth]{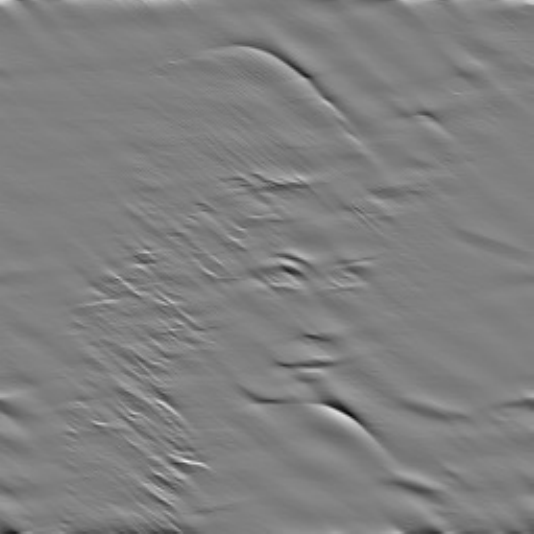} & \includegraphics[width=0.22\textwidth]{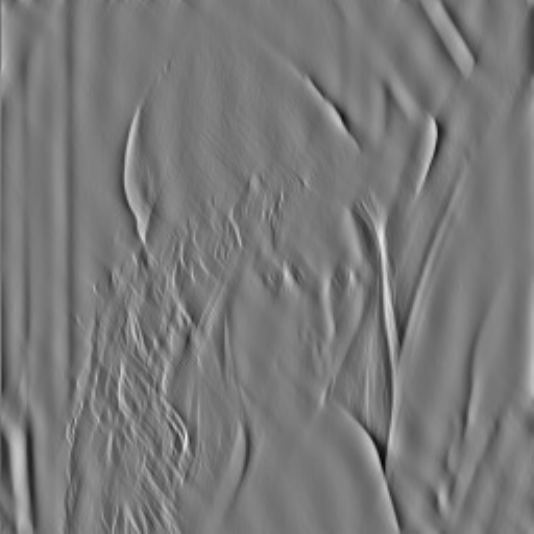} & \includegraphics[width=0.22\textwidth]{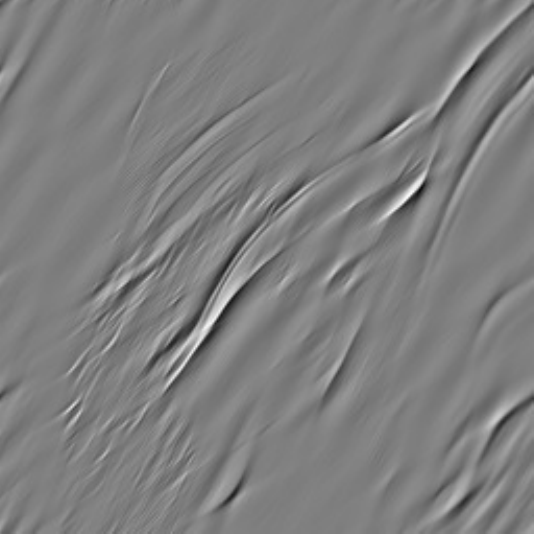} \\
\end{tabular}
\end{center}
\caption{2D Curvelet EWT-I components of the Lena image ($N_s=4, N_{\theta}=4$, the logarithm, the \textit{morpho} preprocessing and the lowest minima detection are used to detect the scales and middle between local 
maxima with \textit{tophat} are used to detect the angles).}
\label{fig:cur1lena}
\end{figure}

\begin{figure}
\begin{center}
\begin{tabular}{cccc}
\includegraphics[width=0.22\textwidth]{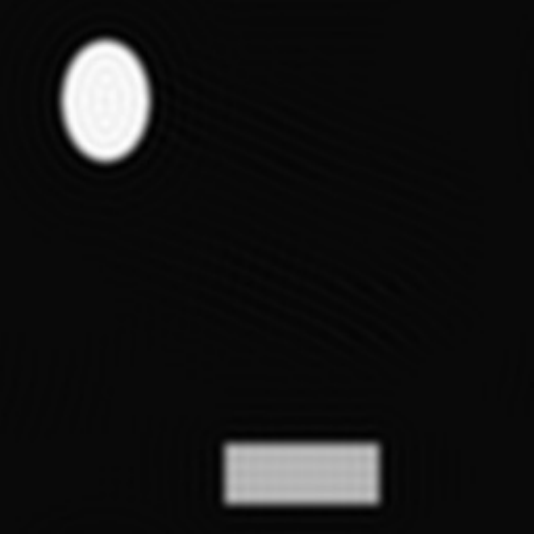} &  &  & \\
\includegraphics[width=0.22\textwidth]{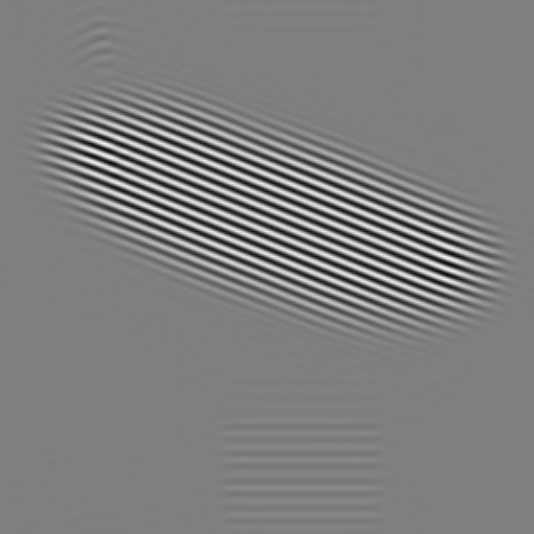} & \includegraphics[width=0.22\textwidth]{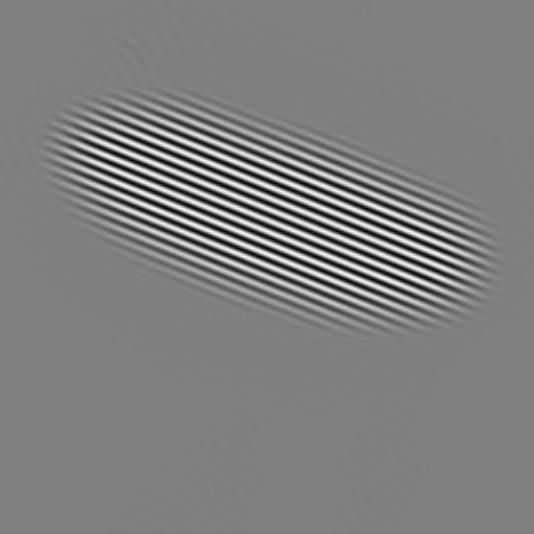} & \includegraphics[width=0.22\textwidth]{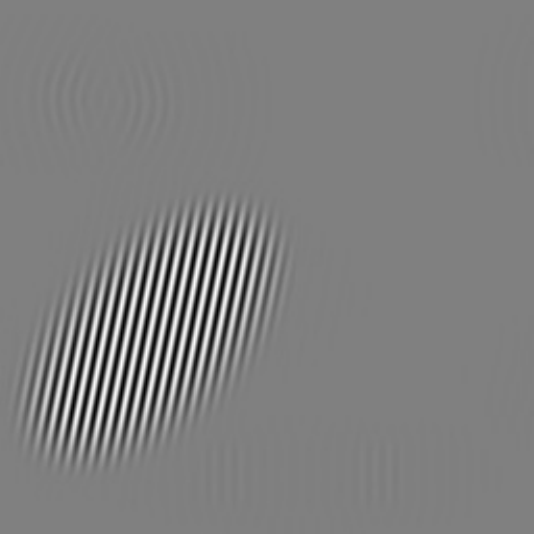} & \includegraphics[width=0.22\textwidth]{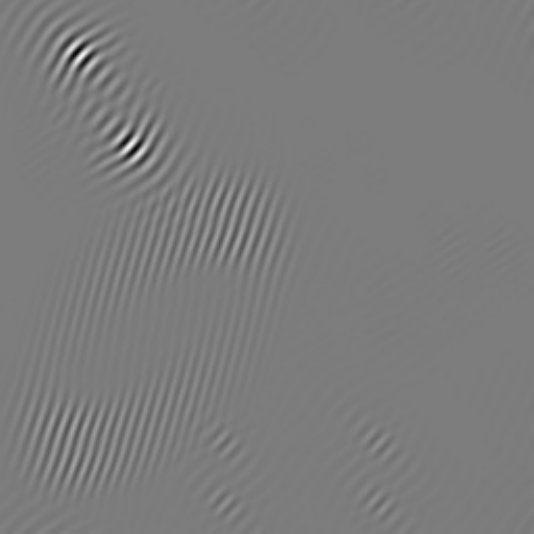}\\
\includegraphics[width=0.22\textwidth]{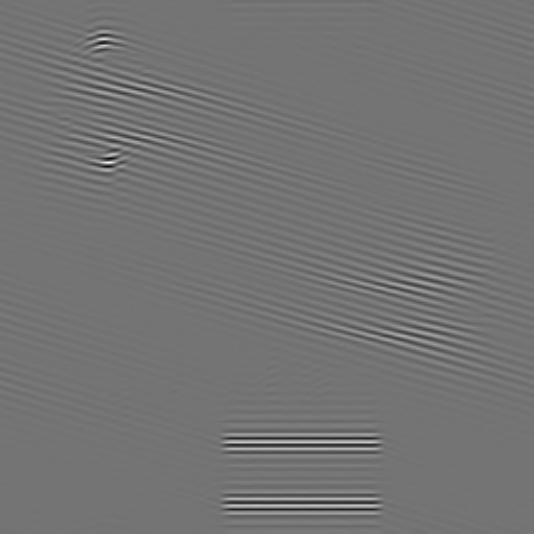} & \includegraphics[width=0.22\textwidth]{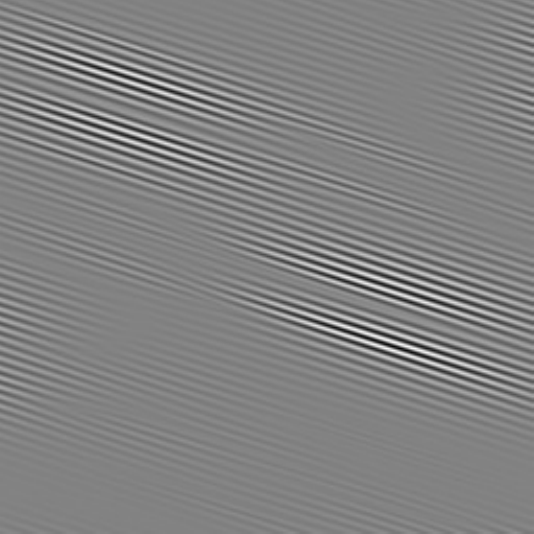} & \includegraphics[width=0.22\textwidth]{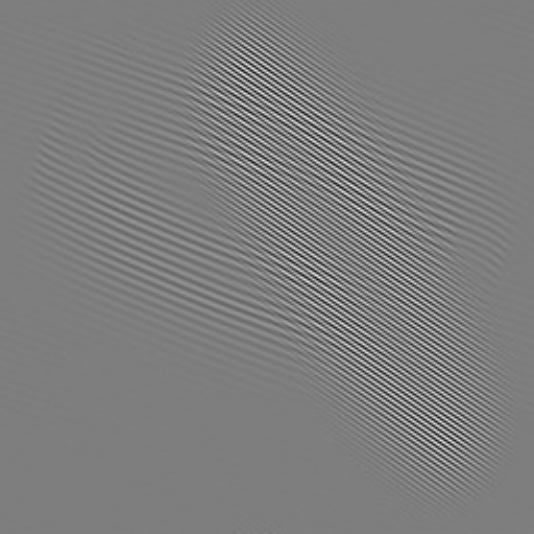} & \includegraphics[width=0.22\textwidth]{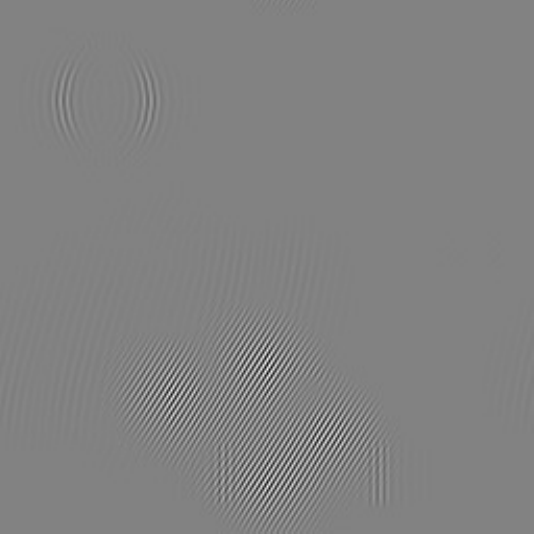}\\
\includegraphics[width=0.22\textwidth]{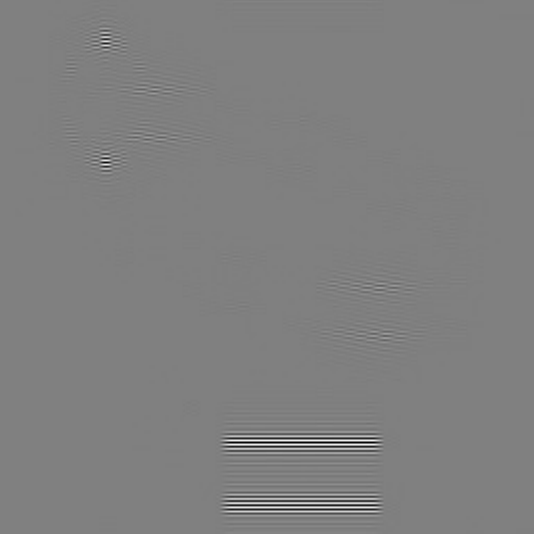} & \includegraphics[width=0.22\textwidth]{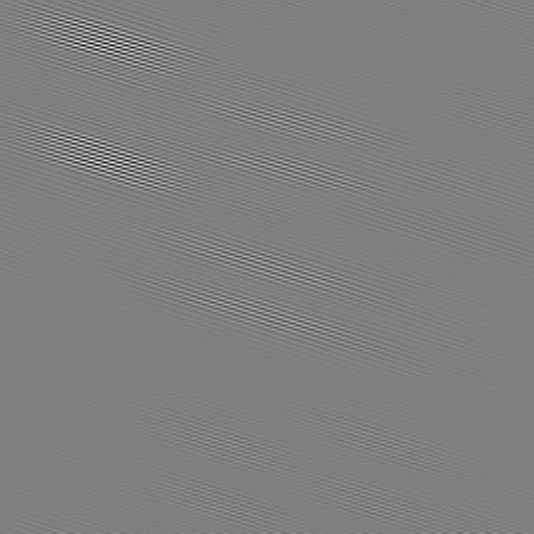} & \includegraphics[width=0.22\textwidth]{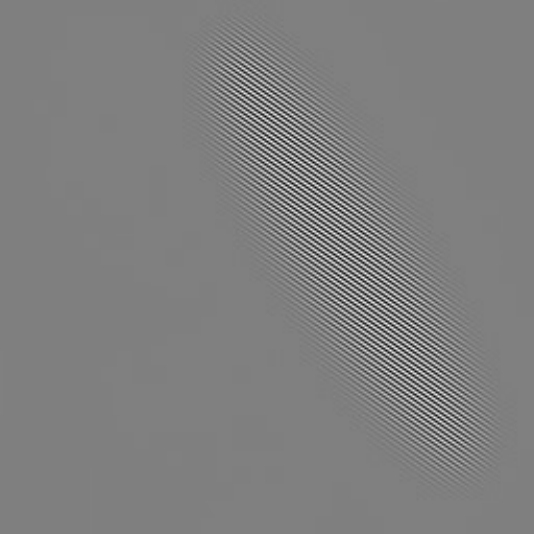} & \includegraphics[width=0.22\textwidth]{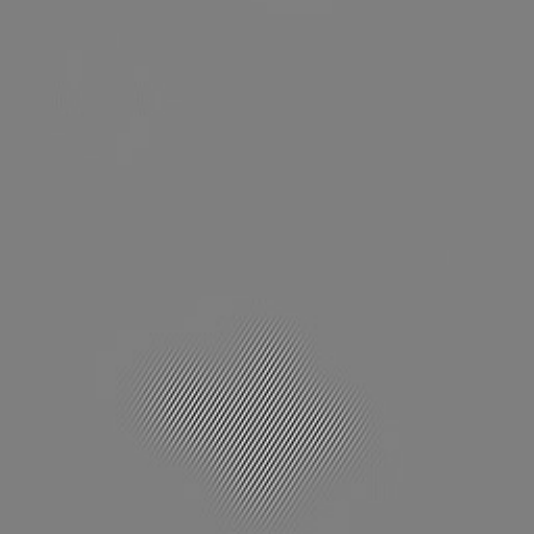}\\
\end{tabular}
\end{center}
\caption{2D Curvelet EWT-II components of the toy image ($N_s=4, N_{\theta}=4$, the logarithm, the \textit{morpho} preprocessing and the lowest minima detection are used to detect the scales and middle between local 
maxima with \textit{tophat} are used to detect the angles).}
\label{fig:cur2toy}
\end{figure}

\begin{figure}
\begin{center}
\begin{tabular}{cccc}
\includegraphics[width=0.22\textwidth]{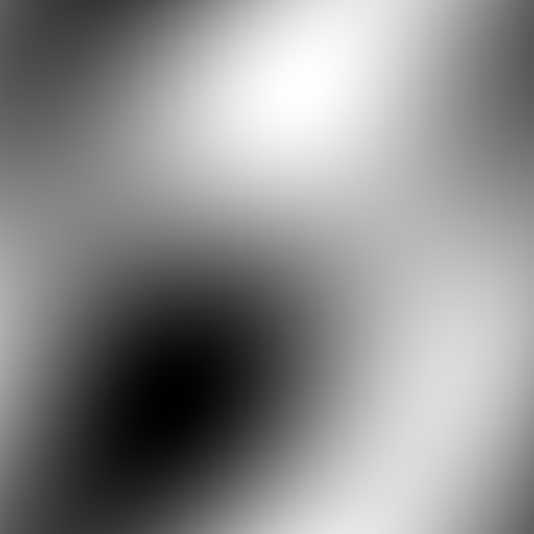} &  &  & \\
\includegraphics[width=0.22\textwidth]{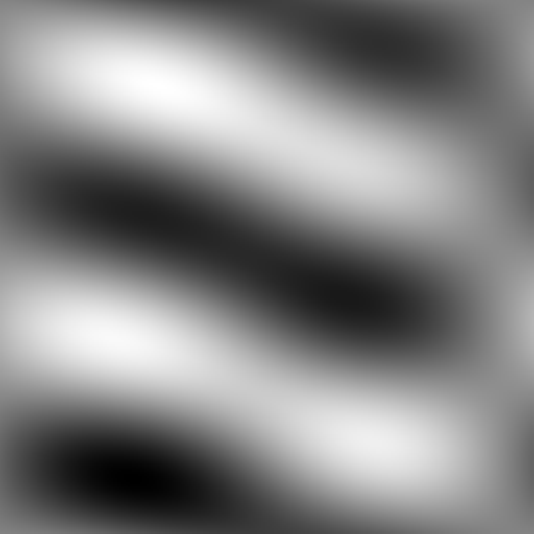} & \includegraphics[width=0.22\textwidth]{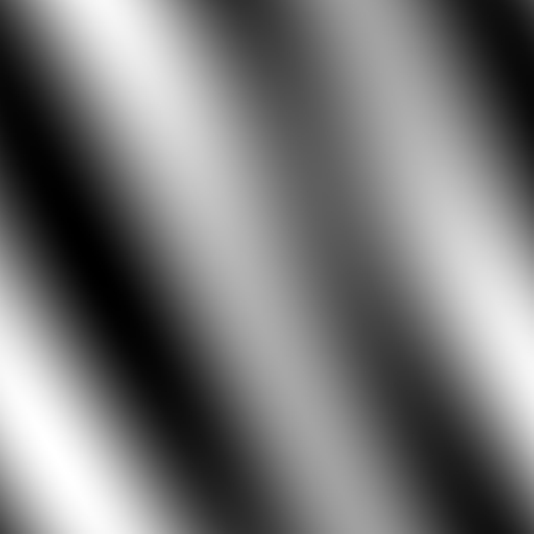} & \includegraphics[width=0.22\textwidth]{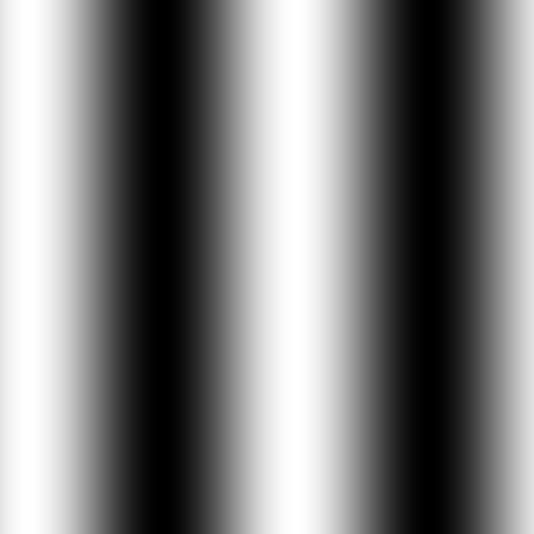} & \includegraphics[width=0.22\textwidth]{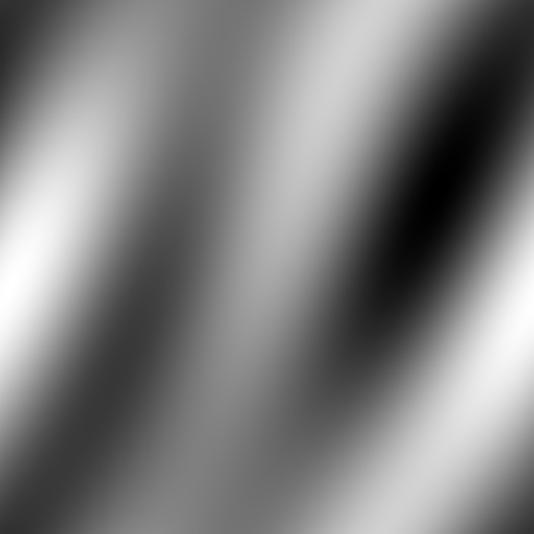}\\
\includegraphics[width=0.22\textwidth]{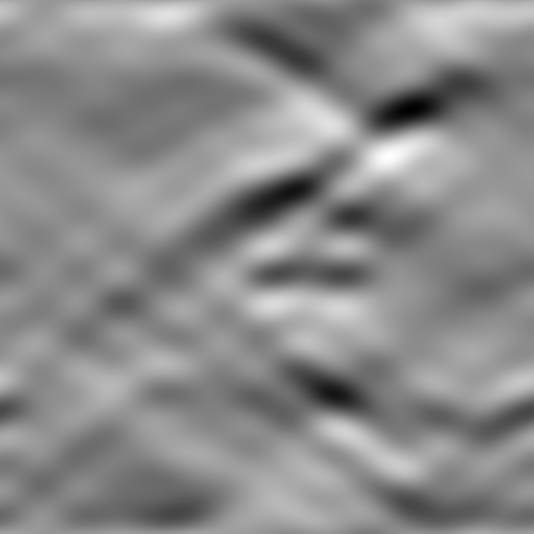} & \includegraphics[width=0.22\textwidth]{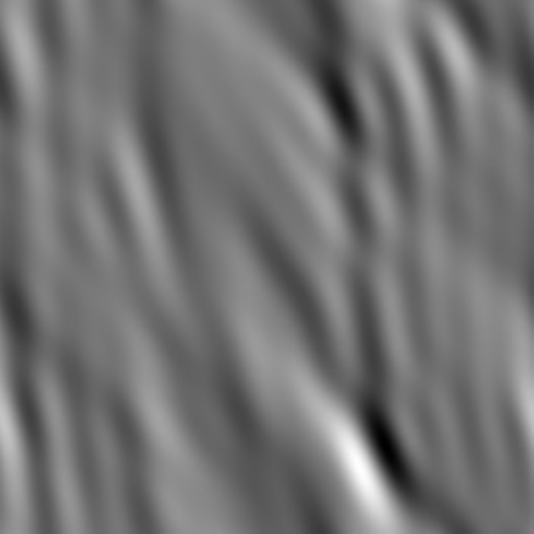} & \includegraphics[width=0.22\textwidth]{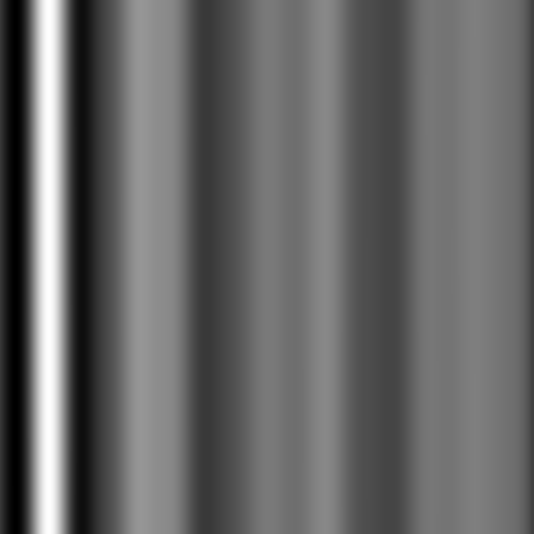} & \includegraphics[width=0.22\textwidth]{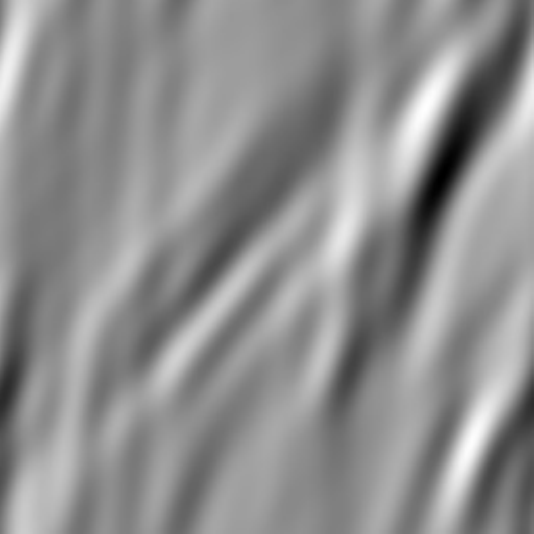} \\
\includegraphics[width=0.22\textwidth]{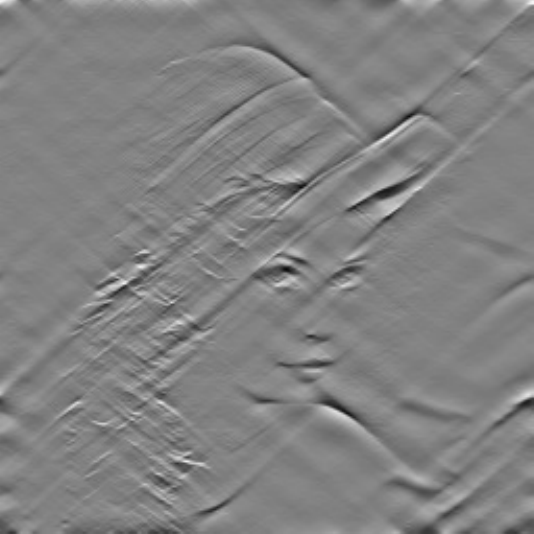} & \includegraphics[width=0.22\textwidth]{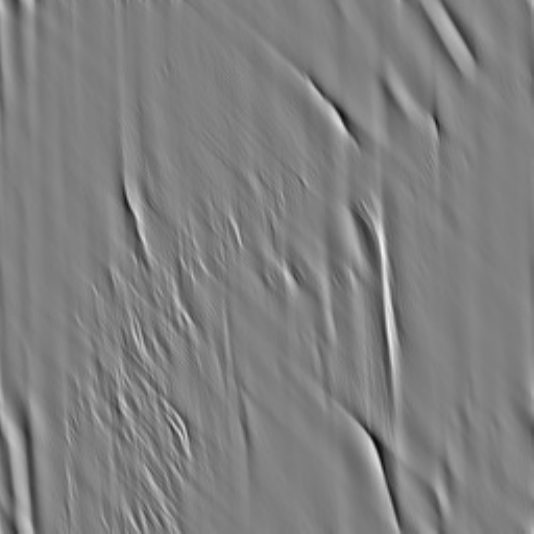} & \includegraphics[width=0.22\textwidth]{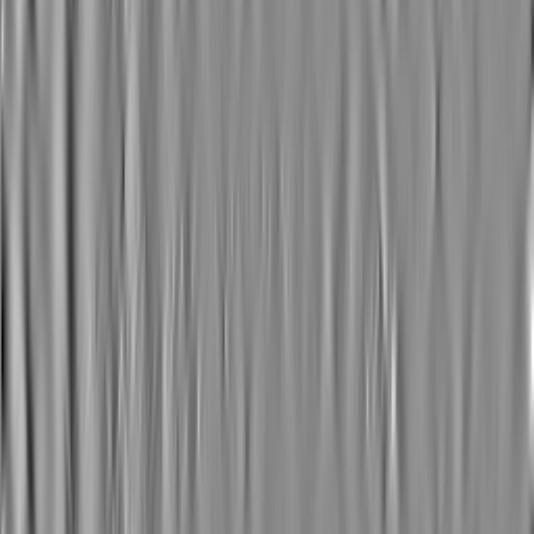} & \includegraphics[width=0.22\textwidth]{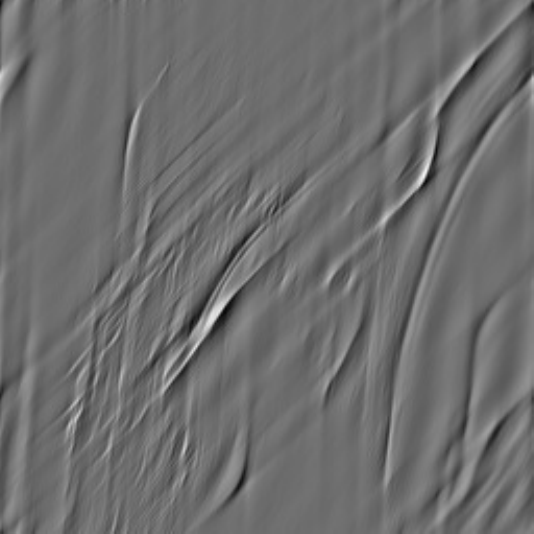} \\
\end{tabular}
\end{center}
\caption{2D Curvelet EWT-II components of the Lena image ($N_s=4, N_{\theta}=4$, the logarithm, the \textit{morpho} preprocessing and the lowest minima detection are used to detect the scales and middle between local 
maxima with \textit{tophat} are used to detect the angles).}
\label{fig:cur2lena}
\end{figure}

\begin{figure}[!t]
\begin{center}
\begin{tabular}{cc}
\includegraphics[width=0.33\textwidth]{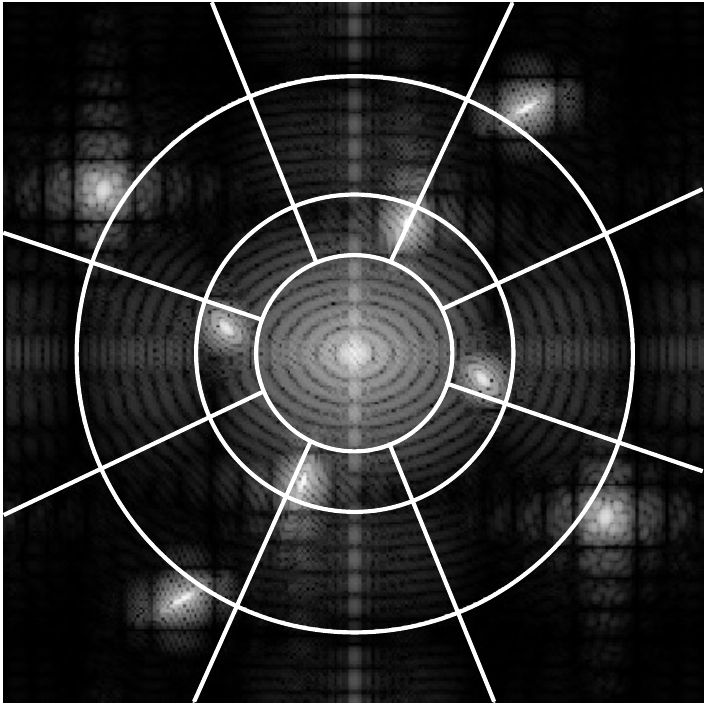} & \includegraphics[width=0.33\textwidth]{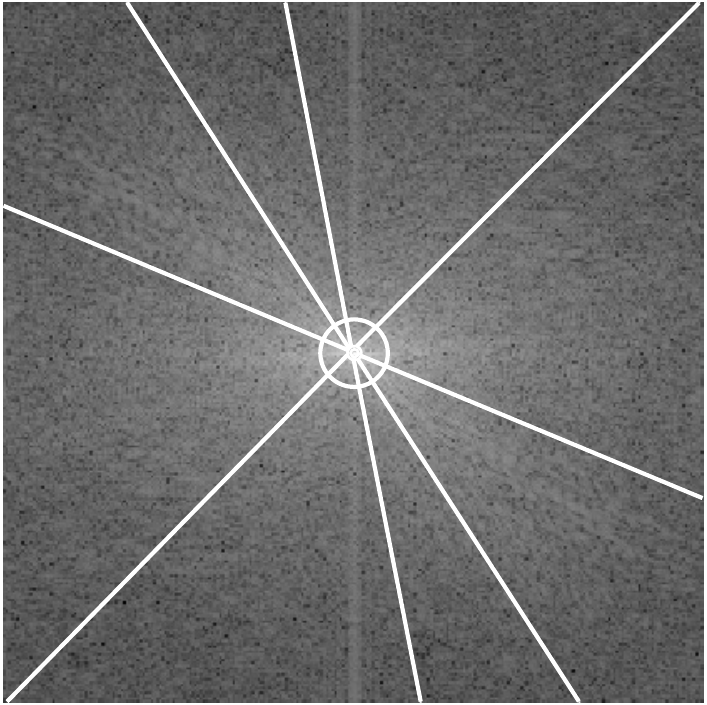} \\
\includegraphics[width=0.33\textwidth]{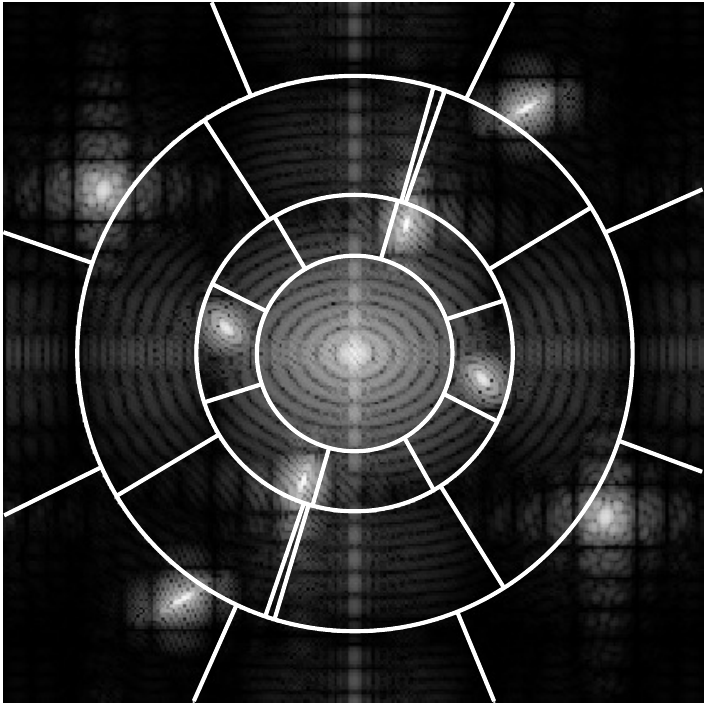} & \includegraphics[width=0.33\textwidth]{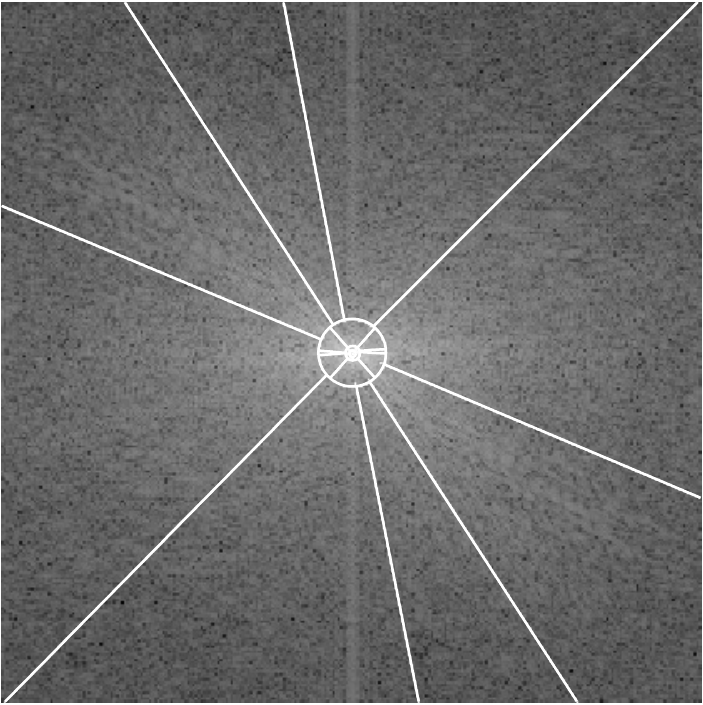}
\end{tabular}
\end{center}
\caption{Detected Fourier boundaries for each test image (toy image on left and Lena on right). Curvelet EWT-I approach on top and Curvelet EWT-II approach on bottom.}
\label{fig:CurvBounds}
\end{figure}

In this section we illustrate the output components obtained from the previously described empirical transforms. The tests are made on two different images. The first one is a toy example containing some flat objects as 
well as different textures lying at different scales and orientations, this image and its Fourier spectrum are given in figure.~\ref{fig:toy}. The second image is the classic Lena image, see figure.~\ref{fig:lena} for 
the image and its Fourier spectrum. Figures.~\ref{fig:tensortoy} and \ref{fig:tensorlena} present the coefficients obtained from the tensor transform approach. The corresponding Fourier supports are 
given in the top row of Figure.~\ref{fig:Bounds}, we can see that the detection performs relatively differently with respect to the horizontal and vertical directions. 
Coefficients obtained from the 2D Littlewood-Paley transform are shown on Figures.~\ref{fig:lptoy} and \ref{fig:lplena}, respectively, while the corresponding detected annuli are given in the bottom row 
of Figure.~\ref{fig:Bounds}. Figures~\ref{fig:ridgetoy} and \ref{fig:ridgelena} present the ridgelet coefficients (the Fourier boundaries correspond to the Littlewood-Paley case) while 
Figures.~\ref{fig:cur1toy}, \ref{fig:cur1lena}, \ref{fig:cur2toy} and \ref{fig:cur2lena} show the curvelet coefficients for each option. The corresponding tiling are given in Figure.~\ref{fig:CurvBounds}. 
For the toy example, it is clear that the Fourier supports for each transform try to separate the different modes.

%---Tensor---

%---LP---

%---Ridgelet---

%---Curvelet---

\subsection{Denoising}
In order to have some intuition on the applicability of the proposed transforms, in this section, we present some denoising 
experiments. The used denoising method consists in three simple steps: we perform the transform of a noisy image, then we apply a soft thresholding on the 
obtained coefficients and finally we perform the inverse transform to get the denoised image. This method is not the best denoising method 
in the literature but it is sufficient to get some clues on the proposed transforms. The noisy images are built from the toy and Lena images, 
mentioned in the previous section, on which we add an additive Gaussian noise. The used threshold is given by $\delta\sqrt{2\log N_p}$, where $N_P$ 
is the number of pixels and $\delta$ is a tuning parameter, optimized on each experiment to get the best denoising results. To compare the 
denoising efficiency, we compute two different metrics (the reference image is denoted $I_r$ and the denoised image $I_d$):
\begin{itemize}
\item the usual PSNR (Peak Signal to Noise Ratio), defined by $PSNR=10\log_{10}\left(\frac{Max^2}{MSE}\right)$, where $Max$ is the maximum possible value in 
the image and $MSE=\frac{1}{N_P}\|I_r-I_d\|_2^2$,
\item the SSIM (Structural Similarity Index Measure \cite{ssim}), this metric is more consistent with the human eye perception than the PSNR. This 
metric is defined by $$SSIM=\frac{(2\mu_r\mu_d+c_1)(2\sigma_{rd}+c_2)}{(\mu_r^2+\mu_d^2+c_1)(\sigma_r^2+\sigma_d^2+c_2)},$$ where 
$\mu_r,\mu_d$ are the averages of $I_r$ and $I_d$ and $\sigma_r,\sigma_d,\sigma_{rd}$ the variances and covariance of $I_r$ and $I_d$, 
respectively. The constants $c_1$ and $c_2$ are defined by $c_1=(0.01L)^2$ and $c_2=(0.03L)^2$ where $L$ is the dynamic range (see 
\cite{ssim} for full details).
\end{itemize}
In order to get a good intuition, we also perform the same denoising technique with the classic wavelet, ridgelet and curvelet transforms. 
The corresponding results are given in Table.~\ref{tab1}. Even though this table shows that the classic curvelet transform outperforms all the other 
ones, it is interesting to note that, in the case of the toy image, the empirical ridgelet and curvelet transforms provide better 
results (in terms of SSIM) than the classic wavelet and ridgelet transforms. Surprisingly, on Lena, the empirical tensor wavelet transform 
is the best empirical approach (and even beats the classic ridgelet transform). If these experiments do not give the empirical approach 
as the best ones, they show that the empirical methods can behave totally differently accordingly to the kind of processed image. These 
results lead us to think that it should be interesting to use a more advanced denoising method like considering different thresholds adapted 
to each subbands or finding the image which have the 
sparsest representation in each transform domain by solving an $L^1-L^2$ problem like in \cite{Shen2010}.

\begin{table}[!h]
\caption{Denoising results (noise variance: $\sigma=10$ for Lena and $\sigma=1$ for the toy image).}
\label{tab1}
\begin{center}
\begin{tabular}{|c|c|c||c|c|} \hline
& \multicolumn{2}{c||}{Lena} & \multicolumn{2}{c|}{Toy image} \\ \hline
& PSNR & SSIM & PSNR & SSIM\\ \hline
Noisy  & 28.096956 & 0.709367 & 48.143956 & 0.456469\\ \hline \hline
Wavelet & 31.151106 & 0.876431 & 52.365411 & 0.678816 \\ \hline
Ridgelet & 29.018199 & 0.811983 & 55.038634 & 0.708038\\ \hline
Curvelet & \textbf{31.762547} & \textbf{0.9019} & \textbf{58.977812} & \textbf{0.928533}\\ \hline
EWT Tensor & 28.956752 & 0.821412 & 52.480827 & 0.648255\\ \hline
EWT Ridgelet & 9.610751 & 0.744987 & 57.10384 & 0.821604\\ \hline
EWT Curvelet I & 16.253502 & 0.781414 & 52.79413 & 0.712689\\ \hline
EWT Curvelet II & 19.802442 & 0.769585 & 52.366085 & 0.740332\\ \hline
\end{tabular}
\end{center}
\end{table}

\section{Conclusion}\label{sec:conc}
In this paper, we extend the 1D Empirical Wavelet Transform proposed in \cite{Gilles2013} to the 2D case. Thus we revisit the tensor wavelet transform, the 2D Littlewood-Paley transform, the ridgelet 
transform and the curvelet transform and their inverses. We show that it is possible to build adaptive wavelet frames which are pretty flexible to use in practice. We also propose different 
options concerning how to perform the Fourier boundary detection. It is clear that depending on how this detection task is performed the corresponding expansions can be very different. This step 
is probably the most important step in these empirical transforms. Experiments seem to show that the algorithm used to do this detection depends on the final goal and the type of analyzed image. 
Indeed, we will not chose the same algorithm if the Fourier spectrum has or not clearly distinct modes.\\
On the practical side, these transforms provide sparse representations of images which could be used for restoration purposes (denoising, 
deblurring for instance). We present 
some denoising results based on a very simple approach, which do not give the advantage to the empirical transforms but lead us to think that more advanced 
denoising methods, like the ones presented in 
\cite{Cai2009,Cai2009a,Cai2009b}, must be investigated. Such study will be presented in a future work.\\
Another extension of this work should be the opportunity to perform a true 2D segmentation of the Fourier spectrum and then build 2D wavelets on the obtained partitioning. If this idea seems 
reasonable, it opens some challenging questions like: how to efficiently perform such segmentation? How to build such wavelet on arbitrary supports?

\section{Acknowledgement}
This work is partially founded by the following grants NSF DMS-0914856, NSF DMS-1118971, ONR N00014-08-1-1119, ONR N0014-09-1-0360, ONR MURI USC, the UC Lab Fees Research and the Keck Foundation.

%=== APPENDICES ====
\appendix
\section{Pseudo-Polar Fourier Transform}\label{app:ppfft}
\begin{figure}[!h]
\centering\includegraphics[scale=0.3]{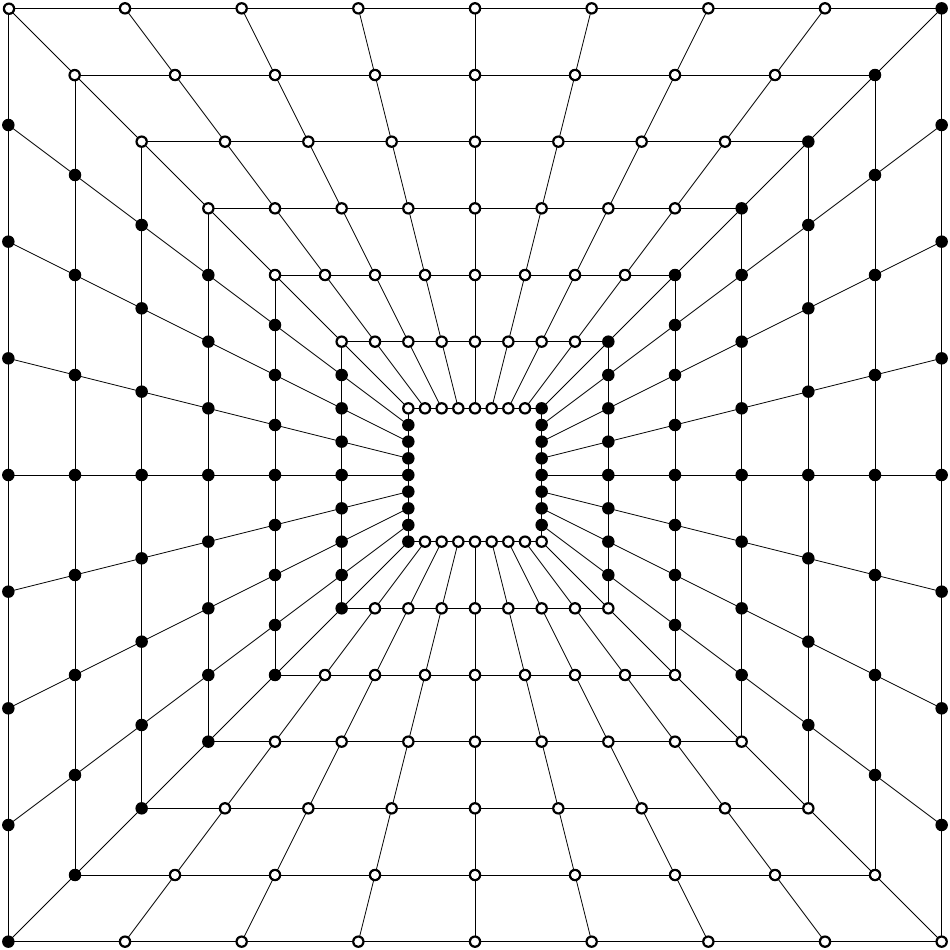}
\caption{Pseudo-polar grid.}
\label{fig:ppgrid}
\end{figure}

The polar Fourier transform is a useful tool for applications like tomography, radon transform computations or in the definition of the ridgelet transform. In \cite{PseudoPolar}, the authors propose a fast 
algorithm to compute the Pseudo-Polar Fourier transform where the polar frequency points are not defined on a regular grid but on a pseudo grid like the one depicted in Figure.~\ref{fig:ppgrid}.
Providing the frequencies $(\omega_1,\omega_2)$ sampled on this pseudo grid, and assuming that the input image $f$ has $N\times N$ pixels, the pseudo-polar Fourier transform is defined by
\begin{equation}
\F_P(f)(\omega_1,\omega_2)=\sum_{x_1=0}^{N-1}\sum_{x_2=0}^{N-1}f(x_1,x_2)exp\left(-\imath (x_1\omega_1+x_2\omega_2)\right).
\end{equation}
The inverse transform is performed by a least square minimization scheme. Providing a pseudo-polar Fourier transform $f_P$ of an image $f$, a gradient descent is used to perform  
\begin{equation}
f=\arg\min_x \|\F_P(x)-f_P\|_2^2.
\end{equation}
See \cite{PseudoPolar} for all numerical details. The authors freely provide a Matlab implementation of this algorithm.

\section{Mathematical morphology operators}\label{app:morpho}
In the boundary detection section (\ref{sec:bound}), we used some mathematical morphological operators that we recall in this appendix. The two basic operators are the dilation ($Dil$) and erosion ($Er$), 
they are respectively defined by
\begin{equation}
Dil(f,b)(x)=\sup_y(f(y)+b(x-y))\qquad;\qquad Er(f,b)(x)=\inf_y(f(y)-b(x-y)),
\end{equation}
where the function $b$ is called the structuring function. It is usual in the literature to consider flat structuring functions: $b(x)=0$ if $x\in B$ and $b(x)=-\infty$ otherwise. Then the dilation 
and erosion operators resumed to 
\begin{equation}
Dil(f,b)(x)=\sup_{y-x\in B}(f(y))\qquad;\qquad Er(f,b)(x)=\inf_{y-x\in B}(f(y)),
\end{equation}
the only remaining parameter is the size of the support $B$. From these two basic operators, we can define the opening and closing operators:
\begin{equation}
Cl(f,b)=Er(Dil(f,b),b) \qquad;\qquad Op(f,b)=Dil(Er(f,b),b).
\end{equation}
Basically, the closing operator ``fills'' holes smaller than the structuring function while the opening operator remove spikes smaller than the structuring function. 

%=== BIBLIOGRAPHY ===

\end{document}